\pdfoutput=1
\documentclass[11pt]{amsart}
\usepackage{amssymb, adjustbox, enumerate, amsbsy, stmaryrd}
\usepackage{geometry}
\usepackage{todonotes}
\usepackage{amsfonts, amssymb, amscd}
\numberwithin{equation}{section}

\usepackage[symbol]{footmisc}

\usepackage{bm}
\usepackage{verbatim}
\usepackage{mathrsfs}
\usepackage{graphicx}
\usepackage{tikz-cd}
\usepackage{subcaption}
\usepackage{listings}
\usepackage{subfiles}
\usepackage[toc,page]{appendix}
\usepackage{mathtools}
\usepackage{comment}
\usepackage{enumerate}
\usepackage{enumitem}
\usepackage[all]{xy}

\usepackage{graphicx}
\graphicspath{{images/}}

\usepackage{appendix}
\usepackage{hyperref}
\hypersetup{
    colorlinks=true,
    citecolor=red,
    linkcolor=blue,
    filecolor=magenta,      
    urlcolor=red,
}
\newcommand{\bQ}{\mathbb{Q}}

\newcommand{\bb}{\bm{b}}

\newcommand{\Qq}{\mathbb{Q}}

\newcommand{\Rr}{\mathbb{R}}

\newcommand{\vol}{\operatorname{vol}}
\newcommand{\Center}{\operatorname{center}}

\newcommand{\Exc}{\operatorname{Exc}}

\newcommand{\Bir}{\operatorname{Bir}}
\newcommand{\Aut}{\operatorname{Aut}}
\newcommand{\glct}{\operatorname{glct}}

\newcommand{\tmld}{{\operatorname{tmld}}}
\newcommand{\PA}{{\operatorname{PA}}}
\newcommand{\A}{{\operatorname{A}}}
\newcommand{\mld}{{\operatorname{mld}}}
\newcommand{\sdeg}{{\operatorname{sdeg}}}

\newcommand{\pet}{\operatorname{pet}}

\newcommand{\lcm}{\operatorname{lcm}}

\newcommand{\lct}{\operatorname{lct}}

\newcommand{\Supp}{\operatorname{Supp}}

\newcommand{\mult}{\operatorname{mult}}
\newcommand{\Rct}{\operatorname{Rct}}

\newcommand{\Div}{\operatorname{Div}}

\newcommand{\Ii}{\Gamma}

\newcommand{\cy}{\mathrm{cy}}
\newcommand{\rc}{\mathrm{rc}}

\newcommand{\ft}{\mathrm{ft}}

\newtheorem{thm}{Theorem}[section]

\newtheorem{cor}[thm]{Corollary}
\newtheorem{lem}[thm]{Lemma}
\newtheorem{prop}[thm]{Proposition}

\newtheorem{claim}[thm]{Claim}

\theoremstyle{definition}
\newtheorem{defn}[thm]{Definition}

\theoremstyle{definition}
\newtheorem{rem}[thm]{Remark}

\newtheorem{nota}[thm]{Notation}

\newtheorem{setup}[thm]{Set-up}

\theoremstyle{definition}

\begin{document}

\title{On explicit bounds of Fano threefolds}
\author{Caucher Birkar and Jihao Liu}

\address{Yau Mathematical Sciences Center, JingZhai, Tsinghua University, Hai Dian District, Beijing, China 100084}
\email{birkar@tsinghua.edu.cn}

\address{Department of Mathematics, Northwestern University, 2033 Sheridan Rd, Evanston, IL 60208, USA}
\email{jliu@northwestern.edu}

\subjclass[2020]{14J30,14J45,14E30,14C20}
\date{\today}

\begin{abstract} In this paper, we study the explicit geometry of threefolds, in particular, Fano varieties. We find an explicitly computable positive integer $N$, such that all but a bounded family of Fano threefolds have $N$-complements. This result has many applications on finding explicit bounds of algebraic invariants for threefolds. We provide explicit lower bounds for the first gap of the $\mathbb R$-complementary thresholds for threefolds, the first gap of the global lc thresholds, the smallest minimal log discrepancy of exceptional threefolds, and the volume of log threefolds with reduced boundary and ample log canonical divisor. We also provide an explicit upper bound of the anti-canonical volume of exceptional threefolds. While the bounds in this paper may not and are not expected to be optimal, they are the first explicit bounds of these invariants in dimension three.
\end{abstract}

\maketitle
\tableofcontents

\section{Introduction}

We work over an algebraically closed field of characteristic zero.

Since the proof of the Borisov-Alexeev-Borisov conjecture \cite{Bir19,Bir21}, a significant amount of work has been dedicated to the classification of Fano varieties. One natural method to classify Fano varieties is to study algebraic invariants related to them, such as the $\alpha$-invariant, the minimal log discrepancy (mld), the anti-canonical volume, the effective birationality, the Fano index, the $\Rr$-complementary threshold, and the boundedness of $N$-complements.

Numerous papers have shown that these invariants satisfy pleasing arithmetic properties. For example, Fano varieties with $\alpha$-invariant greater than $1$ belong to a bounded family \cite{Bir19}, and Fano varieties with mld bounded away from $0$ also belong to a bounded family \cite{Bir21}. In particular, these varieties have both an upper bound and a lower bound on their anti-canonical volume \cite{Bir19}, and the effective birationality property holds \cite{Bir19}. The ACC (ascending chain condition) has also been established for Fano indices \cite{HMX14} and $\Rr$-complementary thresholds \cite{HLS19,Sho20}. Finally, we have the existence and the boundedness of $N$-complements for Fano varieties \cite{Bir19,HLS19,Sho20}.

It is natural to ask whether we can leverage some explicit structure of these algebraic invariants. In particular, since exceptional Fano varieties are always K-stable \cite{Tia87,OS12}, studying the explicit structure of these algebraic invariants is not only an intriguing question but will also be beneficial for the construction of K-moduli spaces.

However, finding the optimal value of the invariants listed above is a very challenging question, and only a few results are known for surfaces and canonical threefolds. In \cite{CC08}, the authors show that the optimal lower bound of the anti-canonical volume of a canonical (weak) Fano threefold is $\frac{1}{330}$. In a very recent work of the second author with Shokurov \cite{LS23}, the authors show that the optimal $1$-gap of surface $\Rr$-complementary thresholds is equal to $\frac{1}{13}$, while the smallest mld of exceptional Fano surfaces is equal to $\frac{3}{35}$.

Nevertheless, we still hope to obtain at least some explicit bounds, if not optimal, of these algebraic invariants. Many works have been done in this direction for surfaces \cite{Kol94,AM04,Lai16,AL19,Liu23,LL23,LS23}, and canonical or terminal (weak) Fano threefolds \cite{Pro05,Pro07,Pro10,Che11,Pro13,CJ16,CJ20,Jia22,Bir23,Jia23,JZ23a,JZ23c}, but few works (\cite{Lai16,Bir23,JZ23b}) have focused on threefolds with worse than canonical singularities. Therefore, in this paper, we will study the explicit behavior of algebraic invariants of klt Fano threefolds (that are not necessarily canonical). As applications, we will also get results on the explicit behavior of algebraic invariants for (log) Calabi-Yau threefolds and varieties of (log) general type.

\medskip

\noindent\textbf{Explicit complements for non-exceptional Fano threefolds}. The first result of our paper is the explicit boundedness of $N$-complements of non-exceptional Fano threefolds. Recall that a Fano variety $X$ is called \emph{non-exceptional} if and only if $\alpha(X)\leq 1$ \cite[Theorem 1.7]{Bir21}, where $\alpha(X)$ is Tian's alpha invariant of $X$, and is called \emph{exceptional} otherwise.

\begin{thm}[Explicit boundedness of $N$-complements for non-exceptional Fano threefolds]\label{thm: explicit bdd nonexc fano threefold intro}
Any non-exceptional klt Fano threefold has an $N$-complement for some positive integer 
$$N<192\cdot (42)!\cdot 84^{128\cdot 42^5}.$$
\end{thm}

Since exceptional Fano varieties of fixed dimension form a bounded family \cite[Theorem 1.11]{Bir19}, Theorem \ref{thm: explicit bdd nonexc fano threefold intro} provides an explicit boundedness result of $N$-complements for all but a bounded family of Fano threefolds. We remark that any non-exceptional Fano surface has an $N$-complement for some integer $N\in\{1,2,3,4,6\}$ (cf. \cite[Theorem 3.1]{Sho00}).

Theorem \ref{thm: explicit bdd nonexc fano threefold intro} has many applications in explicit geometry of threefolds, particularly on estimating extremal values and gaps of algebraic invariants for threefolds.

\medskip

\noindent\textbf{Explicit gap for $\Rr$-complementary thresholds}. Recall that a projective pair $(X,B)$ is called \emph{$\Rr$-complementary} if there exists $B^+\geq B$ such that $(X,B^+)$ is lc and $K_X+B^+\sim_{\mathbb R}0$. The $\Rr$-complementary threshold ($\Rct$ for short) measures how far a pair is away from being not $\Rr$-complementary. It is known that the set of $\Rr$-complementary thresholds for Fano type varieties satisfies the ACC (\cite[Theorem 8.20]{HLS19}, \cite[Theorem 21]{Sho20}). 

In particular, there should be a gap between (non-trivial) $\Rr$-complementary thresholds and $1$ in any dimension. For surfaces, the gap is known to be $\frac{1}{13}$ (\cite[Theorem 1.3]{LS23}). For threefolds, the question is much more difficult and we are not able to find the optimal gap. Nevertheless, we are able to find an explicit gap in this paper by applying Theorem \ref{thm: explicit bdd nonexc fano threefold intro}:

\begin{thm}[Explicit gap of $\Rr$-complementary thresholds]\label{thm: rct 1gap threefold}
Let $I_0:=2\cdot 84^{256\cdot 42^5+338}$. For any projective Fano threefold $X$ and any effective $\Qq$-Cartier Weil divisor $S$ on $X$, let
$$\Rct(X,0;S):=\sup\{t\geq 0\mid (X,tS)\text{ is }\Rr\text{-complementary}\}.$$
Then $\Rct(X,0;S)=1$ or $\Rct(X,0;S)<1-\frac{1}{I_0}$. 
\end{thm}

\medskip

\noindent\textbf{Explicit gap of global lc threshold}. In this paper, global lc thresholds stand for values $b$, such that there exist lc Calabi-Yau pairs $(X,bS)$ for some non-zero effective Weil divisor $S$. (In some references, Tian's $\alpha$-invariant is also called the ``global lc threshold". We will not use this notation in this paper.) Theorem \ref{thm: rct 1gap threefold} immediately implies the following theorem on the explicit $1$-gap of global lc thresholds:

\begin{thm}\label{thm: explicit 1gap glct dim 3}
  Let $I_0:=2\cdot 84^{256\cdot 42^5+338}$. Let $(X,bS)$ be an lc Calabi-Yau pair such that $S$ is a non-zero effective Weil divisor and $b<1$. Then $b<1-\frac{1}{I_0}$.
\end{thm}

\medskip

\noindent\textbf{Explicit bounds for exceptional varieties}. A normal projective variety $X$ is called \emph{exceptional} if $|-K_X|_{\mathbb Q}\not=\emptyset$, and for any $D\in |-K_X|_{\mathbb Q}$, $(X,D)$ is klt. As explained above, exceptional Fano varieties are exactly Fano varieties with Tian's alpha invariant strictly larger than $1$, and form a bounded family. As applications of Theorems \ref{thm: explicit bdd nonexc fano threefold intro} and \ref{thm: rct 1gap threefold}, we can find several explicit bounds of algebraic invariants of exceptional threefolds. 

First, we get an explicit lower bound of the mlds for exceptional threefolds:

\begin{thm}[Explicit lower bound of exceptional threefold mld]\label{thm: exc mld 3fold}
Let $I_0:=2\cdot 84^{256\cdot 42^5+338}$. Then any exceptional threefold is $\frac{1}{I_0}$-klt.
\end{thm}

By abundance (cf. \cite[V.4.6 Theorem]{Nak04}), any  klt Calabi-Yau variety is exceptional. Therefore, Theorem \ref  {thm: exc mld 3fold} immediately implies the following result:

\begin{thm}\label{thm: cy mld 3fold}
 Let $I_0:=2\cdot 84^{256\cdot 42^5+338}$. Then any klt Calabi-Yau threefold is $\frac{1}{I_0}$-klt.   
\end{thm}

By Theorem \ref{thm: exc mld 3fold} and \cite[Theorem 1.1]{JZ23b}, we get the following upper bound of the volume of exceptional threefolds: 

\begin{cor}[Explicit upper bound of exceptional threefold volume]\label{cor: exc upper vol 3fold}
Let
$$V_0:=3200\cdot 84^{1024\cdot 42^5+1352}\approx 10^{10^{11.4}}.$$
Then for any exceptional projective threefold $X$, $\vol(-K_X)<V_0$.
\end{cor}

We remark that we can also apply \cite[Theorem 1.2]{Bir23} to get an explicit upper bound, but that bound will be significantly larger (around $10^{10^{10^{10^{11.3}}}}$).

\medskip

\noindent\textbf{Explicit bounds for volumes of log threefolds}. In moduli theory, when considering degeneration of stable varieties, it is often natural to consider the volume of lc pairs $(X,S)$ where $S$ is a non-zero reduced divisor and $K_X+S$ is ample. In particular, the lower bound of the volume of $K_X+S$ is closely related to the upper bound of the number of irreducible components of the special fiber of the stable degeneration (cf. \cite[Section 4]{Kol94}).

In this paper, we have found a lower bound of such volume for (log) threefolds:

\begin{thm}\label{thm: lower bound k+s volume threefold}
    Let $(X,S)$ be a projective threefold lc pair such that $S$ is a non-zero reduced divisor and $K_X+S$ is ample. Then $$\vol(K_X+S)>\frac{1}{84^{384\cdot 42^5+507}}.$$
\end{thm}

Although the bounds in the previous results may not, and are not expected to be optimal, they represent the first known explicit bounds for threefolds in the literature. Furthermore, we anticipate that the methods employed in this paper will continue to be useful when seeking explicit bounds for algebraic invariants in higher dimensions. We note that the most significant hurdles in finding explicit bounds in high dimensions are the effective $1$-gap of minimal log discrepancies (which is unknown in dimensions $\geq 4$) and the effective finiteness of {\textbf{B}}-representations (which is unknown in dimensions $\geq 2$).

Another approach to studying the explicit behavior of algebraic invariants is to construct examples with conjectural maximal/minimal algebraic invariants. Although these results focus only on constructing examples without extensive proofs, many of the constructed examples are quite inspiring, as they demonstrate that extremal values of algebraic invariants may not appear on ``nice" varieties, such as toric varieties or quasi-smooth general hypersurfaces of weighted projective spaces. Additionally, these examples show that if we can find high-dimensional bounds in the same way as the proof of Theorems \ref{thm: explicit bdd nonexc fano threefold intro}, \ref{thm: rct 1gap threefold}, \ref{thm: exc mld 3fold}, and \ref{cor: exc upper vol 3fold}, the results may not be be very far from being optimal. A sequence of examples, constructed by Totaro \cite{Tot23}, shows that most of these algebraic invariants grow (decrease) at least (at most) double-exponentially, i.e., they grow similarly to $2^{2^{\dim X}}$ $\left(\frac{1}{2^{2^{\dim X}}}\right)$. We will not delve into this approach in detail in this paper, but interested readers may refer to \cite{ETW22,ET23,ETW23,Tot23,TW23} for further information.

\medskip

\noindent\textit{Structure of the proof}. We will present the proof of Theorem \ref{thm: explicit bdd nonexc fano threefold intro} first as it is crucial to prove the other main theorems. Since we are in the non-exceptional case, we can reduce the problem to showing the boundedness of $N$-complements for $\Qq$-factorial dlt threefold pairs $(X,S)$ such that $0\not=S$ is reduced and $-(K_X+S)$ is semi-ample. First we may assume that $(X,S)$ is plt. We consider the following four cases:
\begin{enumerate}
    \item $-(K_X+S)$ is big. In this case, we may use lifting complements \cite[Proposition 6.2]{PS01} from $S$. Thus we only need to show the explicit boundedness of $N$-complements for surface pairs $(S,B_S)$ with coefficients in the standard set $$\left\{1-\frac{1}{n}\Biggm| n\in\mathbb N^+\right\}\cup\{1\}.$$
    \item $\kappa(-(K_X+S))=2$. In this case, we have a contraction $f: X\rightarrow Z$ induced by $-(K_X+S)$ with $\dim Z=2$. As $(X,S)$ is plt, the Stein degree $\sdeg(T/Z)\leq 2$ for any component $T$ of $S$ which dominates $Z$. Thus we also only need to find the explicit boundedness of $N$-complements for surface pairs $(T,B_T)$ with coefficients in the standard set, because we can descend a complement from $T$ to $Z$ and then lift to $X$.
    \item $\kappa(-(K_X+S))=1$. In this case, we have a contraction $f: X\rightarrow Z$ induced by $-(K_X+S)$ with $\dim Z=1$. We use the canonical bundle formula
    $$K_X+B\sim_{\mathbb Q}f^*(K_Z+B_Z+M_Z)$$
    to lift complements from $Z$. There are two remaining tasks at this point:
    \begin{enumerate}
\item The first is to consider the explicit boundedness of monotonic $N$-complements for the generalized pair $(Z,B_Z+M_Z)$, which is relatively easy since $\dim Z=1$.
\item The second is to control the coefficients of $B_Z$ and $M_Z$. Proposition \ref{prop: cbf explicit boundary set} implies that we only need to have boundedness of relative $N$-complements of $(X/Z\ni z,S)$ for any point $z\in Z$, but this will involve complicated arguments. However, we can simplify the situation by letting $K_S+B_S:=(K_X+B)|_S$, so that we also have $$K_S+B_S\sim_{\mathbb Q}f|_S^*(K_Z+B_Z+M_Z).$$
Since $(X,S)$ is plt, the Stein degree $\sdeg(T/Z)$ is $\leq 2$ for a component $T$ of $S$ which dominates $Z$. Thus we may reduce the construction of a complement on $Z$ to the case when $T\rightarrow Z$ is a contraction. Now, by \cite[3.1 Theorem]{Sho00}, there exists a relative $N$-complement of $(T/Z\ni z,B_T)$ for some $N\in\{1,2,3,4,6\}$, and we may control the coefficients of $B_Z$ and the Weil index of $M_Z$ by applying Proposition \ref{prop: cbf explicit boundary set}.
\end{enumerate}
\item $\kappa(-(K_X+S))=0$. In this case, we have $K_X+S\equiv 0$. We may run a $K_X$-MMP and assume that there exists a Mori fiber space $X\rightarrow Z'$. If $\dim Z'>0$, we can apply similar arguments as in the $\kappa(-(K_X+B))=2$ or $1$ case. If $\dim Z'=0$, then $\rho(X)=1$, and by applying Kawamata-Viehweg vanishing theorem, we can lift complements from $S$.
 \end{enumerate}
Now let's assume that $(X,S)$ is not plt. There are two ways to deal with this case: one way is to directly apply the results in the unpublished preprint \cite[Theroems 4,5]{FFMP22}. The other way is to, again, consider the following four cases:
\begin{enumerate}
    \item $-(K_X+S)$ is big. In this case, we may apply \cite[Lemma 2.7]{Bir21} and reduce to the case when $(X,S)$ is plt.
    \item $\kappa(-(K_X+S))=2$. In this case, the argument is the same as the plt case as we still have $\sdeg(T/Z)\leq 2$ for any component $T$ of $S$.
    \item $\kappa(-(K_X+S))=1$. In this case, we can either reduce to the plt case, or we have $2(K_X+S)\sim 0$ over the generic point of $Z$. Possibly passing to a dlt model and running a $K_X$-MMP, we may assume that $K_X$ is $\frac{1}{2}$-lc over the generic point of $Z$ and $X\rightarrow Z$ is a $K_X$-Mori fiber space. Thus any general fiber $F$ of $X\rightarrow Z$ is a $\frac{1}{2}$-lc del Pezzo surface. By using the $\gamma$-invariant in \cite{LS23}, we have $\rho(F)\leq 9$, so we have $\sdeg(T/Z)\leq 11$ for any component $T$ of $S$. The explicit bound is now controlled by lifting complements from $T$.
    \item $\kappa(-(K_X+S))=0$. In this case, similar to the plt case, we may reduce to the case when we have a $K_X$-Mori fiber space $X\rightarrow Z'$, and by applying similar arguments as above, we may assume that $\dim Z'=0$ and $\rho(X)=1$. Now $S$ is a semi-lc surface and we can lift complement from $S$ by applying the Kawamata-Viehweg vanishing theorem. Finally, to provide an explicit boundedness of complements on $S$, we can apply results in \cite{Xu20}.
\end{enumerate}
 
To summarize, we have reduced the non-exceptional case to the question of the explicit boundedness of monotonic $N$-complements for surface pairs with standard coefficients, which are not necessarily of Fano type. For this question, $N$-complements are known to exist for non-exceptional surfaces for some $N\in\{1,2,3,4,6\}$. On the other hand, for exceptional surfaces with standard coefficients, we can establish the existence of $N$-complements by finding an explicit bound on the Cartier indices (Lemma \ref{lem: index epsilonlc fano surface}) and applying the base-point-freeness theorem.

Now we turn to the proof of Theorem \ref{thm: explicit 1gap glct dim 3}. We only need to find a positive real number $\delta$, such that for any threefold lc log Calabi-Yau pair $(X,B)$ where $$B\geq (1-\delta)\Supp B\not=0,$$ we have $B=\Supp B$. By applying known results for surfaces (cf. \cite{Kol94}) we may reduce to the case when $X$ is a $\Qq$-factorial Fano threefold of Picard number $1$ and $(X,\Supp B)$ is lc. There are two possibilities:
\begin{enumerate}
    \item $\vol(-K_X)$ is bounded from above by an explicit real number. In this case, if $B\not=\Supp B$, then we may pick a component $S$ of $B$ and consider the volume of $$K_{S^\nu}+B_{S^\nu}:=(K_X+\Supp B)|_{S^\nu},$$ 
    where $S^\nu$ is the normalization of $S$. By assumption, $(S^\nu,B_{S^\nu})$ is an lc pair with standard coefficients and $K_{S^\nu}+B_{S^\nu}$ is ample, so the volume of $K_{S^\nu}+B_{S^\nu}$ has an explicit lower bound \cite[Section 10]{AL19}. We obtain a contradiction by comparing the volumes.
    \item $\vol(-K_X)$ is bounded away from $0$ by an explicit real number. In this case we can find an explicitly computable real number $\delta_1$ such that $(X,(1-\delta_1)B)$ is not exceptional. Then, we can follow the same lines of the proof of the non-exceptional case to show that $(X,(1-\delta_1)B)$ has a monotonic $N$-complement for some explicit positive integer $N$. Such a monotonic $N$-complement must also be an $N$-complement of $(X,\Supp B)$, hence $B=\Supp B$ and we are done. In practice, to avoid repetitive work, we combine these arguments with the proof of Theorem \ref{thm: explicit bdd nonexc fano threefold intro}.
\end{enumerate}
The rest of the main theorems immediately follow from Theorem \ref{thm: explicit 1gap glct dim 3}. More precisely, by Proposition \ref{prop: compare exc mld and glct} below, Theorems \ref{thm: rct 1gap threefold} and \ref{thm: explicit 1gap glct dim 3} are known to be equivalent, and Theorem \ref{thm: exc mld 3fold} is a consequence of Theorem \ref{thm: rct 1gap threefold}. Theorem \ref{thm: cy mld 3fold} is a direct consequence of Theorem \ref{thm: exc mld 3fold}. Corollary \ref{cor: exc upper vol 3fold} is a direct consequence of Theorem \ref{thm: exc mld 3fold} and \cite[Theorem 1.1]{JZ23b}. Finally, Theorem \ref{thm: lower bound k+s volume threefold} is a consequence of Theorem \ref{thm: exc mld 3fold} and \cite[Section 10]{AL19}.

\medskip

\noindent\textit{Structure of the paper}. Section \ref{sec: preliminaries} introduces preliminary results and notations that will be used throughout the paper. In Section \ref{sec: gaps of algebraic invariants}, we study the relationship between extremal values of algebraic invariants related to the structure of Fano varieties and complements. Section \ref{sec: curve} estimates explicit bounds for the natural invariants introduced inT Section \ref{sec: gaps of algebraic invariants} for curves, while Section \ref{sec: surface standard coefficient} estimates explicit bounds for surfaces. In Section \ref{sec: non-exceptional case}, we prove Theorem \ref{thm: explicit bdd nonexc fano threefold intro}. In Section \ref{sec: proof some main theorems}, we prove all the rest of the main theorems. Finally, in Section \ref{sec: Non-standard coefficient case}, we discuss and provide some explicit bounds for pairs with more general coefficients in low dimensions.

\medskip

\noindent\textbf{Acknowledgements}. The first author was supported by a grant from Tsinghua University and a grant of the National Program of Overseas High Level Talent. We would like to thank Jingjun Han, Yuchen Liu, and Lingyao Xie for useful comments, Vyacheslav V. Shokurov for discussions on related topics, and Stefano Filipazzi and Joaqu\'in Moraga for discussions on Subsection \ref{subsec: non-plt}.

\section{Preliminaries}\label{sec: preliminaries}

We adopt the standard notation and definitions in \cite{Sho92,KM98,BCHM10} and will freely use them. 

\subsection{Sets} In this paper, $\mathbb N$ stands for the set $\{0,1,2,\dots\}$, and $\mathbb N^+$ stands for the set $\{1,2,3,\dots\}$.

\begin{defn}\label{defn: DCC and ACC}
Let $\Ii\subset\Rr$ be a set. We say that $\Ii$ satisfies the \emph{descending chain condition} (DCC) if any decreasing sequence in $\Ii$ stabilizes. We say that $\Ii$ satisfies the \emph{ascending chain condition} (ACC) if any increasing sequence in $\Ii$ stabilizes. We say that $\Ii$ is a \emph{hyperstandard} set if there exists a finite set $\Ii_0\subset\mathbb R_{\geq 0}$ such that $0,1\in\Ii_0$, and 
$$\Ii=\left\{1-\frac{\gamma}{n}\Biggm| n\in\mathbb N^+,\gamma\in\Ii_0\right\}\cap [0,1],$$ 
and we also denote $\Ii$ by $\Phi(\Ii_0)$. We denote $$\Phi_p=\Phi\left(\left\{\frac{k}{p}\Biggm| 0\leq k\leq p\right\}\right)$$ for any positive integer $p$.  We call $\Phi_1$ the \emph{standard set}. For any $\Rr$-divisor $D$, we write $D\in\Ii$ if the non-zero coefficients of $D$ belong to $\Ii$. In particular, $0\in\Ii$ for any $\Ii$.
\end{defn}

We will frequently use the following two facts in this paper, which can be checked easily:
\begin{lem}\label{lem: sum in hyperstandard set}
 Let $p$ and $m$ be two positive integers. 
 \begin{enumerate}
     \item  For any $\gamma_1,\dots,\gamma_m\in\Ii_p$ such that $\sum_{i=1}^m\gamma_i\in [0,1]$, we have $\sum_{i=1}^m\gamma_i\in\Ii_p$.
     \item For any $\gamma\in\Ii_p$ and positive integer $n$, we have
     $$\frac{n-1+\gamma}{n}\in\Ii_p.$$
 \end{enumerate}
\end{lem}
\begin{proof}
    (1) We may apply induction on $m$ and reduce to the case when $m=2$. We write $\gamma_i=1-\frac{k_i}{pn_i}$ for $i=1,2$, such that $k_i\in\mathbb N,n_i\in\mathbb N^+$, and $k_i\leq p$ for $i=1,2$. If $n_i\geq 2$ for each $i$, then $$1\geq \gamma_1+\gamma_2\geq\left(1-\frac{p}{2p}\right)+\left(1-\frac{p}{2p}\right)=1,$$
    so $\gamma_1+\gamma_2=1\in\Phi_p$ and we are done. Thus we may assume that $n_i=1$ for some $i$. Possibly reordering indices, we may assume that $n_1=1$. Then
    $$1-\frac{k_2-n_2(p-k_1)}{pn_2}=\gamma_1+\gamma_2\leq 1.$$
    Let $k':=k_2-n_2(p-k_1)$, then $k'$ is an integer, and
    $$p\geq k_2\geq k'\geq 0.$$
    Thus $$\gamma_1+\gamma_2=1-\frac{k'}{pn_2}\in\Phi_p$$
    and we are done.

    (2) We write $\gamma=1-\frac{k}{pl}$, where $k\in\mathbb N,l\in\mathbb N^+$, and $k\leq p$. Then 
    $$\frac{n-1+\gamma}{n}=\frac{npl-k}{npl}=1-\frac{k}{p(nl)}\in\Phi_p.$$
\end{proof}

\subsection{Pairs and singularities}

\begin{defn}
A \emph{contraction} is a projective morphism $f: X\rightarrow Z$ such that $f_*\mathcal{O}_X=\mathcal{O}_Z$. 
\end{defn}

\begin{defn}\label{defn sing}
A \emph{pair} $(X/Z\ni z, B)$ consists of a contraction $\pi: X\rightarrow Z$, a (not necessarily closed) point $z\in Z$, and an $\mathbb{R}$-divisor $B\geq 0$ on $X$, such that $K_X+B$ is $\Rr$-Cartier over a neighborhood of $z$. If $\pi$ is the identity map and $z=x$, then we may use $(X\ni x, B)$ instead of $(X/Z\ni z,B)$. In addition, if $B=0$, then we use $X\ni x$ instead of $(X\ni x,0)$. If $(X\ni x,B)$ is a pair for any codimension $\geq 1$ point $x\in X$, then we call $(X,B)$ a pair.
\end{defn}

\begin{defn}[Singularities of pairs]\label{defn: relative mld}
 Let $(X\ni x,B)$ be a pair and $E$ a prime divisor over $X$ such that $x\in \Center_XE$. Let $f: Y\rightarrow X$ be a log resolution of $(X,B)$ such that $\Center_Y E$ is a divisor, and suppose that $K_Y+B_Y=f^*(K_X+B)$ over a neighborhood of $x$. We define $a(E,X,B):=1-\mult_EB_Y$ to be the \emph{log discrepancy} of $E$ with respect to $(X,B)$.
 
 For any prime divisor $E$ over $X$, we say that $E$ is \emph{over} $X\ni x$ if $\Center_XE=\bar x$. We define
 $$\mld(X\ni x,B):=\inf\{a(E,X,B)\mid E\text{ is over }X\ni x\}$$
 to be the \emph{minimal log discrepancy} (\emph{mld}) of $(X\ni x,B)$. We define $$\mld(X,B):=\inf\{a(E,X,B)\mid E\text{ is exceptional over }X\}.$$
 We define
 $$\tmld(X,B):=\inf\{a(E,X,B)\mid E\text{ is over }X\}$$
to be the \emph{total minimal log discrepancy (tmld)} of $(X,B)$.
 
 Let $\epsilon$ be a non-negative real number. We say that $(X\ni x,B)$ is lc (resp. klt, $\epsilon$-lc, $\epsilon$-klt) if $\mld(X\ni x,B)\geq 0$ (resp. $>0$, $\geq\epsilon$, $>\epsilon$). We say that $(X,B)$ is lc (resp. klt, $\epsilon$-lc, $\epsilon$-klt) if $\tmld(X,B)\geq 0$ (resp. $>0$, $\geq\epsilon$, $>\epsilon$). We say that $(X,B)$ is \emph{canonical} (resp. \emph{terminal}, \emph{plt}) if $\mld(X,B)\geq 1$ (resp. $>1,>0$).% We say that $(X/Z\ni z,B)$ is \emph{canonical} (resp. \emph{terminal}) if $(X/Z\ni z,B)$ is $1$-lc (resp. $1$-klt).
 
 Let $(X,B)$ be a pair. We say that $(X,B)$ is \emph{dlt} if there exists a log resolution $f: Y\rightarrow X$ of $(X,B)$, such that for any prime $f$-exceptional divisor $E$, $a(E,X,B)>0$.
 \end{defn}

\begin{defn}\label{defn: alct local}  Let $a$ be a non-negative real number, $(X\ni x,B)$ (resp. $(X,B)$) an lc pair, and $D\geq 0$ an $\Rr$-Cartier $\Rr$-divisor on $X$. We define
$$\lct(X\ni x,B;D):=\sup\{t\mid t\geq 0, (X\ni x,B+tD)\text{ is lc}\}$$
$$\text{(resp. }\lct(X,B;D):=\sup\{t\mid t\geq 0, (X,B+tD)\text{ is lc}\} \text{)}$$
to be the \emph{lc threshold} of $D$ with respect to $(X\ni x,B)$ (resp. $(X,B)$). 

For any positive integer $d$ and set $\Ii\subset [0,1]$, we define
$$\lct(d,\Ii):=\{\lct(X,B;D)\mid \dim X=d,(X,B)\text{ is lc}, B\in\Ii,D\in\mathbb N^+\}.$$
\end{defn}

\begin{defn}
Let $(X,B)$ be an lc pair. A \emph{dlt modification} of $(X,B)$ is a projective birational morphism $f: Y\rightarrow X$, such that
\begin{enumerate}
    \item $(Y,B_Y)$ is $\Qq$-factorial dlt, where $K_Y+B_Y:=f^*(K_X+B)$.
    \item for any prime $f$-exceptional divisor $E$, $a(E,X,B)=0$.
\end{enumerate}
We say that $(Y,B_Y)$ is a \emph{dlt model} of $(X,B)$. We will freely use the fact that dlt models exist for any lc pair (cf. \cite[Theorem 3.1]{AH12}).
\end{defn}

\begin{defn}
Let $X\rightarrow Z$ be a contraction. We say that $X$ is \emph{of Fano type} over $Z$ if there exists $B\geq 0$ on $X$ such that $(X,B)$ is klt and $-(K_X+B)$ is ample$/Z$. 
\end{defn}

\begin{thm}[{\cite[Corollary 1.3.2]{BCHM10}}]\label{thm: fano type mds}
Fano type varieties are Mori dream spaces. In particular, let $X\rightarrow Z$ be a contraction such that $X$ is of Fano type over $Z$, and let $D$ be an $\Rr$-Cartier $\Rr$-divisor on $X$. Then we may run a $D$-MMP$/Z$, and any sequence of $D$-MMP$/Z$ terminates with either a good minimal model$/Z$ or a Mori fiber space$/Z$.
\end{thm}

\subsection{Complements}

\begin{defn}\label{defn: complement}
Let $N$ be a positive integer, and $(X/Z\ni z,B)$ and $(X/Z\ni z,B^+)$ two pairs. We say that $(X/Z\ni z,B^+)$ is an \emph{$\Rr$-complement} of $(X/Z\ni z,B)$ if 
\begin{itemize}
    \item $(X/Z\ni z,B^+)$ is lc,
    \item $B^+\geq B$, and
    \item $K_X+B^+\sim_{\Rr}0$ over a neighborhood of $z$.
\end{itemize}
We say that $(X/Z\ni z,B^+)$ is an \emph{$N$-complement} of $(X/Z\ni z,B)$ if
\begin{itemize}
\item $(X/Z\ni z,B^+)$ is lc,
\item $nB^+\geq \lfloor (n+1)\{B\}\rfloor+n\lfloor B\rfloor$, and
\item $n(K_X+B^+)\sim 0$ over a neighborhood of $z$.
\end{itemize}
Here $\{B\}$ means the fractional part of $B$. We say that $(X/Z\ni z,B)$ is $\Rr$-complementary if $(X/Z\ni z,B)$ has an $\Rr$-complement. We say that $(X/Z\ni z,B^+)$ is a \emph{monotonic $N$-complement} of $(X/Z\ni z,B)$ if $(X/Z\ni z,B^+)$ is an $N$-complement of $(X/Z\ni z,B)$ and $B^+\geq B$.
\end{defn}

 \begin{defn}
 Let $(X,B)$ be a projective pair.

 $(X,B)$ is called \emph{log Calabi-Yau} if $K_X+B\equiv 0$. $X$ is called \emph{Calabi-Yau} if $(X,0)$ is log Calabi-Yau. $(X,B)$ is called \emph{exceptional} if $|-(K_X+B)|_{\mathbb R}\not=\emptyset$, and for any $D\in |-(K_X+B)|_{\mathbb R}$, $(X,B+D)$ is klt. $(X,B)$ is called \emph{log Fano} if $-(K_X+B)$ is ample and $(X,B)$ is klt. $X$ is called \emph{of Fano type} if $(X,\Delta)$ is log Fano for some $\Delta$.  $X$ is called \emph{Fano} if $-K_X$ is ample and $X$ is klt. A Fano surface is also called a \emph{del Pezzo} surface. 
 
 We denote by $\rho(X)$ the Picard number of $X$.
 \end{defn}

\begin{lem}\label{lem: antimmp does not contract lc place}
Let $(X/Z,B)$ be an $\Rr$-complementary pair. Then for any steps of a $-(K_X+B)$-MMP$/Z$, no component of $\lfloor B\rfloor$ is contracted. 
\end{lem}
\begin{proof}
Let $X\dashrightarrow X'$ be steps of a $-(K_X+B)$-MMP$/Z$ and let $B'$ be the image of $B$ on $X'$. Let $E$ be a component of $\lfloor B\rfloor$. Then $(X'/Z,B')$ is $\Rr$-complementary, so
$$0=a(E,X,B)\geq a(E,X',B')\geq 0.$$
Thus $a(E,X,B)=a(E,X',B')$, so $E$ is not contracted by $X\dashrightarrow X'$.
\end{proof}

\subsection{Stein degree}

\begin{defn}[Stein degree, {\cite[Page 8]{Bir22}}]\label{defn: stein degree}
Let $f: S\rightarrow Z$ be a projective morphism between varieties and let $S\xrightarrow{\pi} V\xrightarrow{\tau} Z$ be the Stein factorization of $f$. We define the Stein degree of $S$ over $Z$ to be $$\sdeg(S/Z):=\deg(V/Z):=\deg\tau.$$ If $f$ is not surjective, then we define $\sdeg(S/Z):=0$. 
\end{defn}

\begin{lem}\label{lem: plt stein degree 2}
Let $(X,B)$ be a pair and $S$ a component of $\lfloor B\rfloor$. Let $X\rightarrow Z$ be a contraction such that $K_X+B\sim_{\mathbb R,Z}0$ and $S$ is horizontal$/Z$. Suppose that
\begin{enumerate}
    \item either $\dim X-\dim Z=1$, or
    \item $(X,B)$ is plt near $S$ over the generic point of $Z$.
\end{enumerate}
Then $\sdeg(S/Z)\in\{1,2\}$.
\end{lem}
\begin{proof}
Let $F$ be a general fiber of $X\rightarrow Z$, $S_F:=S|_F$, and $B_F:=B|_F$. Then $K_F+B_F\equiv 0$. Moreover, since either $\dim F=1$ or $(X,B)$ is plt near $S$ over the generic point of $Z$, $(F,B_F)$ is plt near $S_F$. By \cite[Proposition 5.1]{KK10}, the non-klt locus of $(F,B_F)$ has at most $2$ connected components, so $S_F$ has at most $2$ connected components. Since $(F,S_F)$ is plt near $S_F$, $S_F$ has at most $2$ irreducible components. Thus $\sdeg(S/Z)\in\{1,2\}$.
\end{proof}

\subsection{Generalized pairs}

Generalized pairs are a natural structure in birational geometry that has origions in Kodaira’s canonical bundle formula for elliptic fibrations. They were formally introduced in \cite{BZ16} as part of the study of effective Iitaka fibrations and have since become a central topic in modern-day birational geometry. In this paper, for simplicity, we will only use projective generalized pairs. 

\begin{defn}
A \emph{projective generalized pair} (\emph{projective g-pair} for short) $(X,B+M)$ consists of a normal projective variety $X$, and $\Rr$-divisor $B\geq 0$ on $X$, and a nef $\Rr$-divisor $M'$ over a model $X'$ of $X$, such that $M$  is the image of $M'$ on $X$ and $K_X+B+M$ is $\Rr$-Cartier. Here we think of $M'$ as a $\bb$-divisor in the sense that when replacing $X'$ by a higher birational model, we replace $M'$ with its pullback. We say that $(X,B+M)$ is NQC if $M'$ can be written as an $\Rr_{\geq 0}$-combination of nef Cartier divisors.

We say that $(X,B+M)$ is \emph{lc} if, possibly replacing the morphism $f: X'\rightarrow X$ by a log resolution of $X$, we have
$$K_{X'}+B'+M'=f^*(K_X+B+M)$$
for some $\Rr$-divisor $B'$, such that the coefficients of $B'$ are $\leq 1$. An $\Rr$-complement (resp. monotonic $N$-complement) of $(X,B+M)$ is an lc g-pair $(X,B^++M)$ such that $B^+\geq B$ and $K_X+B^++M\sim_{\mathbb R}0$ (resp. $N(K_X+B^++M)\sim 0$ and $NM'$ is $\bb$-Cartier).
\end{defn}

\section{Explicit bounds of invariants and their general behavior}\label{sec: gaps of algebraic invariants}

\subsection{Some bounds of algebraic invariants}
\begin{defn}
Let $d$ and $p$ be two positive integers. We define the following sets:
\begin{itemize}
    \item 
    $$\mathcal{S}(d,p):=\{(X,B)\mid (X,B)\text{ is a projective lc pair}, \dim X=d, B\in\Phi_p\}$$
    \item $$\mathcal{S}_{\ft}(d,p):=\{(X,B)\mid (X,B)\in \mathcal{S}(d,p),X\text{ is of Fano type}\}.$$
    \item $$\mathcal{S}_{\rc}(d,p):=\{(X,B)\mid (X,B)\in \mathcal{S}(d,p),(X,B)\text{ is }\Rr\text{-complementary}\}.$$
    \item  $$\mathcal{S}_{\cy}(d,p):=\{(X,B)\mid (X,B)\in \mathcal{S}(d,p),K_X+B\equiv 0\}.$$
    \item $$\mathcal{S}_{\rc,\ft}(d,p):=\mathcal{S}_{\rc}(d,p)\cap\mathcal{S}_{\ft}(d,p),\ \mathcal{S}_{\cy,\ft}(d,p):=\mathcal{S}_{\cy}(d,p)\cap\mathcal{S}_{\ft}(d,p),$$
     \item \begin{align*}
   S_{g}(d,p):=\Bigg\{(X,B+M)\Biggm|
    \begin{array}{r@{}l}
    (X,B+M)\text{ is a projective lc g-pair}\\
    \text{with nef part }
    M',\dim X=d,\\
    X\text{ is of Fano type}, B\in\Phi_p, pM'\text{ is }\bb\text{-Cartier}
    \end{array}\Bigg\}.
    \end{align*}
      \item $$\mathcal{S}_{g,\rc}(d,p):=\{(X,B+M)\mid (X,B+M)\in \mathcal{S}_g(d,p),(X,B+M)\text{ is }\Rr\text{-complementary}\}.$$
\end{itemize}
\end{defn}

\begin{defn}\label{defn: special explicit values}
Let $d$ and $p$ be two positive integers. In the rest of the paper, we adopt the following notations:
\begin{itemize}
    \item $\delta_{\lct}(d,p)$: the $1$-gap of lc thresholds in dimension $d$  with coefficients in the hyperstandard set $\Phi_p$, i.e.
    $$\delta_{\lct}(d,p):=\inf\{1-t\mid t<1, t\in\lct(d,\Phi_p)\}.$$
    \item $\delta_{\glct}(d,p)$: the $1$-gap of global lc thresholds in dimension $d$ with coefficients in the hyperstandard set $\Phi_p$, i.e.
\begin{align*}
   \delta_{\glct}(d,p):=\inf\Bigg\{1-t\Biggm|
    \begin{array}{r@{}l}t<1,
    \text{there exists a projective lc pair }(X,B+tS)\text{ of dimension } d\\\text {such that }K_X+B+tS\equiv 0,
     B\in\Phi_p,\text{and }0\not=S\in\mathbb N^+
    \end{array}\Bigg\}.
    \end{align*}
    \item $I(d,p)$: the maximal index of lc log Calabi-Yau pairs in dimension $d$ with coefficients in the hyperstandard set $\Phi_p$, i.e.
    $$I(d,p):=\inf\{I_0\mid \text{For any }(X,B)\in\mathcal{S}_{\cy}(d,p), I(K_X+B)\sim 0\text{ for some } I\leq I_0\}.$$
     \item $N(d,p)$: the maximal monotonically complementary index of $\Rr$-complementary pairs in dimension $d$ with coefficients in the hyperstandard set $\Phi_p$, i.e.
    \begin{align*}
   N(d,p):=\inf\Bigg\{N_0\Biggm|
    \begin{array}{r@{}l}
    \text{For any }(X,B)\in\mathcal{S}_{\rc}(d,p), (X,B) \text{ has a monotonic }\\ N\text{-complement for some }N\leq N_0
    \end{array}\Bigg\}.
    \end{align*}
     \item $N_{\ft}(d,p)$: the maximal monotonically complementary index of $\Rr$-complementary Fano type pairs in dimension $d$ with coefficients in the hyperstandard set $\Phi_p$, i.e.
    \begin{align*}
   N_{\ft}(d,p):=\inf\Bigg\{N_0\Biggm|
    \begin{array}{r@{}l}
    \text{For any }(X,B)\in\mathcal{S}_{\rc,\ft}(d,p), (X,B) \text{ has a monotonic }\\ N\text{-complement for some }N\leq N_0
    \end{array}\Bigg\}.
    \end{align*}
      \item $N_{g}(d,p)$: the maximal monotonically complementary index of $\Rr$-complementary Fano type g-pairs in dimension $d$ with boundary coefficients in the hyperstandard set $\Phi_p$ and with $p$ as a $\bb$-Cartier index of the moduli part i.e.,
    \begin{align*}
   N_{g}(d,p):=\inf\left\{N_0\Biggm|
    \begin{array}{r@{}l}
    \text{For any }(X,B+M)\in\mathcal{S}_{g,\rc}(d,p), (X,B+M) \text{ has a monotonic }\\
    N\text{-complement for some }N\leq N_0
    \end{array}\right\}.
    \end{align*}
    \item $\delta_{\mld}(d,p)$: the smallest mld of exceptional pairs of dimension $d$ with coefficients in the hyperstandard set $\Phi_p$, i.e.
        $$\delta_{\mld}(d,p):=\inf\{\tmld(X,B)\mid (X,B)\in\mathcal{S}_{\rc}(d,p), (X,B)\text{ is exceptional}\}.$$
    \item $\delta_{\mld,\ft}(d,p)$: the smallest mld of exceptional Fano type pairs of dimension $d$ with coefficients in the hyperstandard set $\Phi_p$, i.e.
        $$\delta_{\mld,\ft}(d,p):=\inf\{\tmld(X,B)\mid (X,B)\in\mathcal{S}_{\rc,\ft}(d,p), (X,B)\text{ is exceptional}\}.$$
\end{itemize}
\end{defn}

\begin{lem}\label{lem: elementary comparison of explicit bound}
Let $d,p$ be positive integers. Then \begin{enumerate}
\item $\mathcal{S}_{\cy}(d,p)\subset\mathcal{S}_{\rc}(d,p)$ and $\mathcal{S}_{\cy,\ft}(d,p)\subset\mathcal{S}_{\rc,\ft}(d,p)$.
\item $I(d,p)\leq N(d,p)$.
\item $N_{\ft}(d,p)\leq N(d,p)$ and $N_{\ft}(d,p)\leq N_g(d,p)$.
\end{enumerate}
\end{lem}
\begin{proof}
(1) follows from the abundance for numerically trivial pairs  (cf. \cite[V.4.6 Theorem]{Nak04}) and divisors on Fano type varieties. (2) follows from (1). (3) follows from the definitions.
\end{proof}

\begin{lem}\label{lem: low dimension epsilon is larger}
Let $d,p$ be two positive integers. Then
$$\delta_{\lct}(d,p)\geq\delta_{\lct}(d+1,p)$$ and $$\delta_{\glct}(d,p)\geq\delta_{\glct}(d+1,p).$$
\end{lem}
\begin{proof}
By \cite[Lemma 11.2(1)]{HMX14}, $\delta_{\lct}(d,p)\geq\delta_{\lct}(d+1,p)$. For any $t\in (0,1)$ and a pair $(W,\Delta+t\Psi)$ of dimension $d$ such that $(W,\Delta+t\Psi)$ is lc,  $K_W+\Delta+t\Psi\equiv 0,\Delta\in\Ii$, and $0\not=\Psi\in\mathbb N^+$, $(W\times E,(\Delta\times E)+t(\Psi\times E))$ is of dimension $d+1$ such that $(W\times E,(\Delta\times E)+t(\Psi\times E))$ is lc,  $K_{W\times E}+(\Delta\times E)+t(\Psi\times E)\equiv 0,(\Delta\times E)\in\Ii$, and $0\not=(\Psi\times E)\in\mathbb N^+$, where $E$ is a smooth elliptic curve. Thus $\delta_{\glct}(d,p)\geq\delta_{\glct}(d+1,p)$. 
\end{proof}

\begin{lem}\label{lem: alternative defn epsilon12}
Let $d,p$ be positive integers. Let $(X,B+B')$ be an lc pair of dimension $d$ such that $B\in\Phi_p$. Let $S:=\Supp B'$.
\begin{enumerate}
    \item Suppose that all components of $B'$ are $\Qq$-Cartier and $B'\in (1-\delta_{\lct}(d,p),1)$. Then $(X,B+S)$ is lc.
    \item Suppose that $X$ is projective, $B'\not=0$, and $B'\in (1-\delta_{\glct}(d,p),1)$. Then $$K_X+B+B'\not\equiv 0.$$
    \item Suppose that $X$ is projective, $B'\in [1-\frac{1}{N(d,p)+1},1)$, and $(X,B+B')$ is $\Rr$-complementary. Then $(X,B+S)$ has a monotonic $N$-complement for some $$N\leq N(d,p).$$
    \item Suppose that $X$ is of Fano type, $B'\in [1-\frac{1}{N(d,p)+1},1)$, and $(X,B+B')$ is $\Rr$-complementary.  Then $(X,B+S)$ has a monotonic $N$-complement for some $$N\leq N_{\ft}(d,p).$$
\end{enumerate}
\end{lem}
\begin{proof}
(1) Since $(X,B+B')$ is lc, $(X,B+(1-\delta_{\lct}(d,p))S)$ is lc. Thus $(X,B+S)$ is lc.

(2) Suppose that $K_X+B+B'\equiv 0$. Possibly replacing $(X,B+B')$ with a dlt model, we may assume that $(X,B+B')$ is $\Qq$-factorial dlt. Let $T$ be a component of $B'$. We run a $$(K_X+B+B'-\epsilon T)\text{-MMP}$$ 
for some $0<\epsilon\ll 1$, which terminates with a Mori fiber space $Y\rightarrow Z$. Let $B_Y,B_Y',T_Y,S_Y$ be the strict transforms of $B,B',T,S$ on $Y$ respectively. Let $F$ be a general fiber of $Y\rightarrow Z$, $B_F:=B_Y|_F,B_F':=B_Y'|_F,T_F:=T_Y|_F$, and $S_F:=S_Y|_F$. Then $T_F\not=0$, so $B_F'\not=0$. Therefore, $(F,B_F+B'_F)$ is a projective lc pair of dimension $\dim F$, $K_F+B_F+B'_F\equiv 0$, $B_F\in\Phi_p$, $B'\not=0$, and $B'\in (1-\delta_{\glct}(d,p),1)$. By Lemma \ref{lem: low dimension epsilon is larger}, $$B'_F\in (1-\delta_{\glct}(\dim F,p),1).$$ If $\dim F<\dim Y$, then by induction on dimension, $K_F+B_F+B'_F\not\equiv 0$, a contradiction. Thus $F=Y$ and $\rho(Y)=1$. 

By construction, $S_Y=\Supp B_Y'$. Since $\rho(Y)=1$, $K_Y+B_Y+S_Y$ is ample, and $K_Y+B_Y+(1-\delta_{\glct}(d,p))S_Y$ is anti-ample, so there exists $t\in (1-\delta_{\glct}(d,p),1)$ such that $K_Y+B_Y+tS_Y\equiv 0$. By the definition of $\delta_{\glct}(d,p)$, $(Y,B_Y+tS_Y)$ is not lc. Thus $(Y,B_Y+S_Y)$ is not lc. We let $t':=\lct(Y,B_Y;S_Y)$. Since $(Y,B_Y+B_Y')$ is lc, $$t'>1-\delta_{\glct}(d,p).$$ We let $E$ be a prime divisor over $Y$ such that $a(E,Y,B_Y+t'S_Y)=0$ and $\mult_ES_Y>0$. Let $g: W\rightarrow Y$ be a dlt modification of $(Y,B_Y+t'S_Y)$ such that $E$ is on $W$, and let $H$ be a $g$-exceptional prime divisor such that $H\cap g^{-1}_*S_Y\not=\emptyset$. We let $G$ be a general fiber of $H\rightarrow\Center_YH$, and let $$K_G+\tilde B_G:=g^*(K_Y+B_Y+t'S_Y)|_G.$$ 
By Lemma \ref{lem: sum in hyperstandard set}, we have the following.
\begin{itemize}
    \item $K_G+\tilde B_G\equiv 0$, $(G,\tilde B_G)$ is lc, and $\dim G<d$,
    \item For any component $D$ of $\tilde B_G$, $$\mult_D\tilde B_G=\frac{n_D-1+\gamma_D+k_Dt'}{n_D}$$ for some $\gamma_D\in\Phi_p$, positive integer $n_D$, and non-negative integer $k_D$, and
    \item There exists a component $C$ of $\tilde B_G$ such that $k_C>0$. In particular,
    $$\mult_C\tilde B_G\geq t'>1-\delta_{\glct}(d,p).$$
\end{itemize}
We may let $$B_G=\sum_{D\mid k_D=0}\frac{n_D-1+\gamma_D}{n_D}D$$ and $$B_G':=\sum_{D\mid k_D>0}\frac{n_D-1+\gamma_D+k_Dt'}{n_D}D.$$ 
By Lemma \ref{lem: sum in hyperstandard set}, $\frac{n-1+\gamma}{n}\in\Phi_p$ for any positive integer $n$ and $\gamma\in\Phi_p$. Thus $(G,B_G+B'_G)$ is an lc pair of dimension $<d$ such that $B_G\in\Phi_p$, $B'_G\not=0$, $B'_G\in (1-\delta_{\glct}(d,p),1)$, and $K_G+B_G+B'_G\equiv 0$. We get a contradiction by induction on dimension.

(3)(4) Possibly replacing $(X,B+B')$ with a dlt model, we may assume that $(X,B+B')$ is $\Qq$-factorial dlt. Then $$\left(X,B+\left(1-\frac{1}{N(d,p)+1}\right)S\right)\text{ (resp. }\left(X,B+\left(1-\frac{1}{N_{\ft}(d,p)+1}\right)S\right)\text{)}$$has a monotonic $N$-complement $(X,B^+)$ for some $N\leq N(d,p)$ (resp. $N\leq N_{\ft}(d,p)$). Then $(X,B^+)$ is a monotonic $N$-complement of $(X,B+S)$.
\end{proof}

\begin{lem}\label{lem: alternative defn ngdp}
Let $d$ and $p$ be two positive integers. Let $(X,B+B'+M)$ be a projective $\Rr$-complementary g-pair of dimension $d$ with nef part $M'$, such that $X$ is of Fano type, $B\in\Phi_p$, $B'\in [1-\frac{1}{N_{g}(d,p)+1},1)$, and $pM'$ is $\bb$-Cartier. Let $S:=\Supp B'$. Then $(X,B+S+M)$ has a monotonic $N$-complement for some $N\leq N_{g}(d,p)$.
\end{lem}
\begin{proof}
Possibly replacing $(X,B+B'+M)$ with a dlt model (cf. \cite[Proposition 3.10]{HL22}), we may assume that $(X,B+B'+M)$ is $\Qq$-factorial. Since $(X,B+B'+M)$ is $\Rr$-complementary and $B'\geq (1-\frac{1}{N_g(d,p)+1})S$, $$\left(X,B+\left(1-\frac{1}{N_g(d,p)+1}\right)S+M\right)$$ is $\Rr$-complementary. Then $$\left(X,B+\left(1-\frac{1}{N_g(d,p)+1}\right)S+M\right)$$ has a monotonic $N$-complement $(X,B^++M)$ for some $N\leq N_g(d,p)$. $(X,B^++M)$ is a monotonic $N$-complement of $(X,B+S+M)$.
\end{proof}

\begin{lem}\label{lem: equality glct d and lct d+1}
Let $d,p$ be two positive integers. Then $\delta_{\lct}(d+1,p)=\delta_{\glct}(d,p)$. In particular, $\delta_{\lct}(d+1,p)=\delta_{\glct}(d,p)\leq\delta_{\lct}(d,p)$.
\end{lem}
\begin{proof}
First we show that $\delta_{\lct}(d+1,p)\geq\delta_{\glct}(d,p)$. The proof is similar to the proof of Lemma \ref{lem: alternative defn epsilon12}(2). By ACC for lc thresholds \cite[Theorem 1.1]{HMX14}, $\delta_{\lct}(d+1,p)>0$, and there exists a pair $(X,B)$ of dimension $d+1$ and a $\Qq$-Cartier Weil divisor $S\geq 0$ on $X$, such that $B\in\Phi_p$ and $$t:=\lct(X,B;S)=1-\delta_{\lct}(d+1,p).$$ 
We let $E$ be a prime divisor over $X$ such that $a(E,X,B+tS)=0$ and $\mult_ES>0$. Let $f: Y\rightarrow X$ a dlt modification of $(X,B+tS)$ such that $E$ is on $Y$, and let $F$ be an $f$-exceptional prime divisor such that $F\cap f^{-1}_*S\not=\emptyset$. We let $G$ be a general fiber of $F\rightarrow\Center_XF$, and let $$K_G+B_G:=f^*(K_X+B+tS)|_G.$$ Then by Lemma \ref{lem: sum in hyperstandard set}, we have the following.
\begin{itemize}
    \item $K_G+B_G\equiv 0$, $(G,B_G)$ is lc, and $\dim G\leq d$,
    \item for any component $D$ of $B_G$, $$\mult_DB_G=\frac{n_D-1+\gamma_D+k_Dt}{n_D}$$ for some $\gamma_D\in\Phi_p$, positive integer $n_D$, and non-negative integer $k_D$, and
    \item there exists a component $C$ of $B_G$, such that $k_C>0$. In particular, $$\mult_CB_G\geq t=1-\delta_{\lct}(d+1,p).$$
\end{itemize}
By Lemma \ref{lem: sum in hyperstandard set}, $\frac{n-1+\gamma}{n}\in\Phi_p$ for any positive integer $n$ and $\gamma\in\Phi_p$. By Lemma \ref{lem: alternative defn epsilon12}(2) and Lemma \ref{lem: low dimension epsilon is larger}, $\delta_{\lct}(d+1,p)\geq\delta_{\glct}(\dim G,p)\geq\delta_{\glct}(d,p)$.

Next we show that $\delta_{\lct}(d+1,p)\leq\delta_{\glct}(d,p)$. Let $t:=1-\delta_{\glct}(d,p)$.  By Lemma \ref{lem: low dimension epsilon is larger}, there exists an lc pair $(X,B+tS)$ of dimension $\leq d$ such that $K_X+B+tS\equiv 0$, $B\in\Phi_p$, $S$ is a non-zero effective Weil divisor, and either $\dim X=1$ or $\delta_{\glct}(\dim X-1,p)>\delta_{\glct}(\dim X,p)$. Possibly replacing $(X,B+tS)$ with a dlt model, we may assume that $X$ is $\Qq$-factorial klt. We run a $(K_X+B)$-MMP which terminates with a Mori fiber space $X'\rightarrow Z$. Let $B',S'$ be the images of $B,S$ on $X'$ respectively and let $F$ be a general fiber of $X'\rightarrow Z$. Then $S'|_F\not=0$. Since either $\dim X=1$ or $\delta_{\glct}(\dim X-1,p)>\delta_{\glct}(\dim X,p)$, $F=X'$. Possibly replacing $X,B,S$ with $X',B',S'$, we may assume that $\rho(X)=1$. We let $\tilde X:=C(X,-K_X)$ be the cone of $X$ and let $\tilde B,\tilde S$ be the divisors induced by $B,S$ on $\tilde X$ respectively. Then $$\lct(\tilde X,\tilde B;\tilde S)=t=1-\delta_{\glct}(d,p),$$ so $$\delta_{\lct}(\dim X+1,p)\leq\delta_{\glct}(d,p).$$ By Lemma \ref{lem: low dimension epsilon is larger}, $$\delta_{\lct}(d+1,p)\leq\delta_{\lct}(\dim X+1,p)\leq\delta_{\glct}(d,p).$$
Thus $\delta_{\lct}(d+1,p)=\delta_{\glct}(d,p)$. The in particular part immediately follows from Lemma \ref{lem: low dimension epsilon is larger}.
\end{proof}

\subsection{Global lc threshold and mld of exceptional pairs}

\begin{defn}
Let $(X/Z\ni z,B)$ be an $\Rr$-complementary pair and let $D\geq 0$ be an $\Rr$-Cartier $\Rr$-divisor on $X$. We define
$$\Rct(X/Z\ni z,B;D):=\sup\{t\geq 0\mid (X/Z\ni z,B+tD)\text{ is }\Rr\text{-complementary}\}.$$
\end{defn}

\begin{lem}\label{lem: rct gap is glct gap}
Let $d,p$ be two positive integers. Let $(X,B)\in\mathcal{S}_{\rc}(d,p)$ be a pair and $S\geq 0$ a non-zero $\Qq$-Cartier Weil divisor on $X$. Assume that
\begin{itemize}
    \item either $X$ is of Fano type, or
    \item the existence of good minimal models holds in dimension $d$.
\end{itemize}
Then $\Rct(X,B;S)=1$ or $\Rct(X,B;S)\leq 1-\delta_{\glct}(d,p)$.
\end{lem}
\begin{proof}
Suppose that the proposition does not hold. Then there exists a pair $(X,B)\in\mathcal{S}_{\rc}(d,p)$ and a Weil divisor $S\geq 0$ on $X$, such that $$t:=\Rct(X,B;S)\in (1-\delta_{\glct}(d,p),1).$$ 
Since $t>\frac{1}{2}$, $S$ is reduced. 

We let $(X,B+tS+G)$ be an $\Rr$-complement of $(X,B+tS)$ and let $f: Y\rightarrow X$ be a $\Qq$-factorial dlt modification of $(X,B+tS+G)$. We let $R:=\lfloor B+tS+G\rfloor$, $G_Y$ the strict transform of $G-G\wedge R$ on $Y$, $S_Y$ the strict transform of $S-S\wedge R$ on $Y$, and
$$K_Y+B_Y+tS_Y+G_Y:=f^*(K_X+B+tS+G).$$
Then $B_Y\in\Phi_p$ and $\lfloor B_Y+tS_Y+G_Y\rfloor=\lfloor B_Y\rfloor$. In particular, there exists $\tau\in (0,t)$ such that $(Y,B_Y+(1+\tau)(tS_Y+G_Y))$ is dlt, and
$$-\tau(K_Y+B_Y+S_Y)\sim_{\mathbb R}K_Y+B_Y+(1+\tau)(tS_Y+G_Y)-\tau S_Y.$$
 Since $\tau\in (0,t)$, $(1+\tau)t-\tau>0$, so $$(Y,\tilde B_Y:=B_Y+(1+\tau)(tS_Y+G_Y)-\tau S_Y)$$ is dlt. 
 
 There are two cases.

\medskip

\noindent\textbf{Case 1}. $-(K_Y+B_Y+S_Y)$ is not pseudo-effective. Since 
$$-\tau(K_Y+B_Y+S_Y)\sim_{\mathbb R}K_Y+\tilde B_Y,$$ 
we may run a $(K_Y+\tilde B_Y)$-MMP with scaling of an ample divisor, which terminates with a Mori fiber space $\pi: Y'\rightarrow Z$. We let $B_{Y'},S_{Y'}$ be the images of $B_Y,S_Y$ on $Y'$ respectively, then $-(K_{Y'}+B_{Y'}+S_{Y'})$ is anti-ample$/Z$. Since $(Y,B_Y+tS_Y)$ is $\Rr$-complementary, $(Y',B_{Y'}+tS_{Y'})$ is $\Rr$-complementary. In particular,  $(Y',B_{Y'}+tS_{Y'})$ is lc. Since $t\in (1-\delta_{\glct}(d,p),1)$, by Lemma \ref{lem: equality glct d and lct d+1}, $(Y',B_{Y'}+S_{Y'})$ is lc. Moreover, $-(K_{Y'}+B_{Y'}+tS_{Y'})$ is pseudo-effective, hence $-(K_{Y'}+B_{Y'}+tS_{Y'})$ is nef$/Z$. Thus there exists $t'\in [t,1)$ such that $(Y',B_{Y'}+t'S_{Y'})$ is lc and $K_{Y'}+B_{Y'}+t'S_{Y'}\equiv_Z0$. Let $F$ be a general fiber of $Y'\rightarrow Z$, $B_{F}:=B_{Y'}|_F$, and $S_F:=S_{Y'}|_F$. Then $K_{F}+B_{F}+t'S_{F}\equiv0$. By Lemma \ref{lem: low dimension epsilon is larger}, 
$$1-\delta_{\glct}(d,p)<t\leq t'\leq 1-\delta_{\glct}(\dim F,p)\leq 1-\delta_{\glct}(d,p),$$
a contradiction.

\medskip

\noindent\textbf{Case 2}. $-(K_Y+B_Y+S_Y)$ is  pseudo-effective. Since $$-\tau(K_Y+B_Y+S_Y)\sim_{\mathbb R}K_Y+\tilde B_Y$$ and $(Y,\tilde B_Y)$ is dlt, by our assumption, we may run a $(K_Y+\tilde B_Y)$-MMP which terminates with a good minimal model $(Y',\tilde B_{Y'})$ of $(Y,\tilde B_Y)$. More precisely, if $X$ is of Fano type, then the MMP terminates with a good minimal model by Theorem \ref{thm: fano type mds}, and if we assume the existence of good minimal models in dimension $d$, then the MMP terminates with a good minimal model by \cite[Theorem 1.5]{Bir11} and \cite[Remark 2.7]{Bir12}.

We let $B_{Y'},S_{Y'}$ be the images of $B_Y,S_Y$ on $Y'$ respectively, then $-(K_Y+B_{Y'}+S_{Y'})$ is semi-ample. Since $(Y,B_Y+tS_Y)$ is $\Rr$-complementary, $(Y',B_{Y'}+tS_{Y'})$ is $\Rr$-complementary. In particular,  $(Y',B_{Y'}+tS_{Y'})$ is lc. Since $t\in (1-\delta_{\glct}(d,p),1)$, by Lemma \ref{lem: equality glct d and lct d+1}, $(Y',B_{Y'}+S_{Y'})$ is lc. Thus we may pick $D\in |-(K_Y+B_{Y'}+S_{Y'})|_{\mathbb R}$ such that $$(Y',B_{Y'}^+:=B_{Y'}+S_{Y'}+D)$$ is lc. Let $p: W\rightarrow Y$, $q: W\rightarrow Y'$ be a common resolution of the induced birational map $Y\dashrightarrow Y'$, and let $$K_X+B^+:=f_*p_*q^*(K_{Y'}+B_{Y'}^+).$$ Then $(X,B^+)$ is an $\Rr$-complement of $(X,B+S)$. Thus $\Rct(X,T;S)=1$, a contradiction.
\end{proof}

The following proposition is a small improvement of Lemma \ref{lem: rct gap is glct gap}.

\begin{prop}\label{prop: rct gap is glct gap}
Let $d,p$ be two positive integers. Let $(X,B+B')$ be an $\Rr$-complementary pair of dimension $d$ such that $B\in\Phi_p$ and $B'\in (1-\delta_{\glct}(d,p),1)$. Let $S\geq 0$ a non-zero $\Qq$-Cartier Weil divisor on $X$. Assume that
\begin{itemize}
    \item either $X$ is of Fano type, or
    \item the existence of good minimal models holds in dimension $d$.
\end{itemize}
Then $\Rct(X,B+B';S)=1$ or $\Rct(X,B+B';S)\leq 1-\delta_{\glct}(d,p)$.
\end{prop}
\begin{proof}
Possibly replacing $(X,B+B')$ with a dlt model we may assume that $(X,B+B')$ is $\Qq$-factorial dlt. Let $T:=\Supp B'$. Suppose that $\Rct(X,B+B';S)\in (1-\delta_{\glct}(d,p),1)$, then $\Rct(X,B;T+S)\in (1-\delta_{\glct}(d,p),1)$, which contradicts Lemma \ref{lem: rct gap is glct gap}.
\end{proof}

\begin{thm}\label{thm: strong version 1.2 liu22}
Let $d,p$ be two positive integers and $\delta:=\delta_{\glct}(d,p)$. Let $(X,B+B')$ be an exceptional pair of dimension $d$ such that $B\in\Phi_p$ and $B'\in (1-\delta,1)$.  Assume that
\begin{itemize}
    \item either $X$ is of Fano type, or
    \item the existence of good minimal models holds in dimension $d$.
\end{itemize}
Then for any $0\leq G\sim_{\mathbb R}-(K_X+B+B')$, $(X,B+B'+G)$ is $\delta$-lc. In particular, $B'=0$.
\end{thm}
\begin{proof}
Suppose that $\tmld(X,B+B'+G)<\delta$. Possibly taking a small dlt modification, we may assume that $X$ is $\Qq$-factorial. We let $E$ be a prime divisor over $X$ such that $a(E,X,B+B'+G)=\tmld(X,B+B'+G)$. Then either $E$ is on $X$ and we let $f: Y\rightarrow X$ be the identity morphism, or there exists a divisorial contraction $f: Y\rightarrow X$ which extracts $E$, such that $Y$ is $\Qq$-factorial. We let $B_Y:=f^{-1}_*B-f^{-1}_*B\wedge E$ and $B'_Y:=f^{-1}_*B'-f^{-1}_*B'\wedge E$. Then $$(Y,B_Y+B_Y'+(1-a(E,X,B+B'+G))E)$$ is $\Rr$-complementary, so $$\Rct(Y,B_Y+B_Y';E)\geq 1-a(E,X,B+B'+G)>1-\delta.$$ By Proposition \ref{prop: rct gap is glct gap}, $\Rct(Y,B_Y+B_Y';E)=1$. Thus $(Y,B_Y+B_Y'+E)$ has an $\Rr$-complement $(Y,B_Y^+)$. Then $(X,B^+:=f_*B_Y^+)$ is an $\Rr$-complement of $(X,B+B')$, such that $a(E,X,B^+)=0$. This is not possible as $(X,B+B')$ is exceptional.
\end{proof}

\begin{prop}\label{prop: compare exc mld and glct}
Let $d$ and $p$ be two positive integers. Then:
\begin{enumerate}
    \item $\delta_{\mld,\ft}(d,p)\geq\delta_{\glct}(d,p).$
    \item Assume the existence of good minimal models in dimension $d$. Then
    $$\delta_{\mld}(d,p)\geq\delta_{\glct}(d,p).$$
    \item    \begin{align*}
   \delta_{\mld,\ft}(d,p)=\inf\left\{\tmld(X,B+B')\Biggm|
    \begin{array}{r@{}l}
    \dim X=d, B\in\Phi_p, B'\in (1-\delta_{\glct}(d,p),1),\\ X\text{ is of Fano type}, (X,B+B')\text{ is exceptional}
    \end{array}\right\}.
    \end{align*}
    \item Assume the existence of good minimal models in dimension $d$. Then
    \begin{align*}
   \delta_{\mld}(d,p)=\inf\left\{\tmld(X,B+B')\Biggm|
    \begin{array}{r@{}l}
    \dim X=d, B\in\Phi_p, B'\in (1-\delta_{\glct}(d,p),1),\\ (X,B+B')\text{ is exceptional}
    \end{array}\right\}.
    \end{align*}
\end{enumerate}
\end{prop}
\begin{proof}
It immediately follows from Theorem \ref{thm: strong version 1.2 liu22}.
\end{proof}

\subsection{Canonical bundle formulas}

In this subsection we prove some results on  canonical bundle formulas. For basic properties and constructions of the canonical bundle formula, especially the canonical bundle formula and its relationship with generalized pairs, we refer the reader to \cite[3.4]{Bir19}, \cite[Section 4]{HL21}, \cite[2.4]{JLX22}, \cite[2.5]{FS23}, \cite[Part III]{CHLX23} and we will freely use the results therein.

The following proposition is a generalization of \cite[Proposition 6.3]{Bir19}, \cite[Proposition 3.1]{CHL23} and the proofs are similar.

\begin{prop}\label{prop: cbf explicit boundary set}
Let $p$ and $q$ be two positive integers, $\delta\in [0,1)$ a real number, $(X,B)$ a projective lc pair, and $f: X\rightarrow Z$ a contraction, such that
\begin{enumerate}
    \item $\dim Z=1$,
    \item $B\in\Phi_p\cup[1-\delta,1)$ (resp. $\Phi_p\cup (1-\delta,1)$),
    \item $K_X+B\sim_{\mathbb R,Z}0$, and
    \item for any (not necessarily closed point) $z\in Z$, $(X/Z\ni z,B)$ has a monotonic $q$-complement $(X/Z\ni z,B_z^+)$, such that $(X,B_z^+)$ has an lc center $F$ with $f(F)=\bar z$. (Here we consider $X$ as an lc center of $(X,B_z^+)$ as well.)
\end{enumerate}
Then we have an lc g-pair $(Z,B_Z+M_Z)$ induced by the canonical bundle formula for $f: (X,B)\rightarrow Z$, such that
$$q(K_X+B)\sim qf^*(K_Z+B_Z+M_Z),$$
$qM_Z$ is a Weil divisor, and
$$B_Z\in\Phi_{pq}\cup[1-\delta,1)\text{ (resp. }B_Z\in \Phi_{pq}\cup(1-\delta,1)\text{)}.$$
\end{prop}

\begin{proof}

\noindent\textbf{Step 1}. In this step we make a choice of $M_Z$ such that $$q(K_X+B)\sim qf^*(K_Z+B_Z+M_Z).$$

By our assumption, $q(K_X+B)\sim 0$ over the generic point of $Z$. Thus there exists $\alpha\in \mathcal{K}(X)$ such that $q(K_X+B)+\Div(\alpha)=0$ near the generic fiber of $f$. In particular, $q(K_X+B)\sim qL$ and $L$ is vertical$/Z$. Since $L\sim_{\mathbb R,Z}0$, $L=f^*L_Z$ for some $L_Z$ on $Z$. We let $B_Z$ be the discriminant part of the canonical bundle formula for $f: (X,B)\rightarrow Z$ and let $M_Z:=L_Z-K_Z-B_Z$.

By our construction, 
$$q(K_X+B)\sim qL=qf^*L_Z=qf^*(K_Z+B_Z+M_Z).$$

\medskip

\noindent\textbf{Step 2}. In this step we show that $B_Z\in\Phi_{pq}\cup [1-\delta,1)$ (resp. $\Phi_{pq}\cup (1-\delta,1)$).

\begin{claim}\label{claim:lctcmpt}
Let $z\in Z$ be a closed point and let $t_z:=\lct(X/Z\ni z,B;f^*z)$. Then $$(X/Z\ni z,B+t_zf^*z)$$ is a monotonic $q$-complement of $(X/Z\ni z,B)$.
\end{claim}
\begin{proof}[Proof of Claim \ref{claim:lctcmpt}]
By assumption (4), we may let $(X/Z\ni z,B_z^+)$ be a monotonic $q$-complement of $(X/Z\ni z,B)$ such that there exists an lc center $F$ of $(X,B_z^+)$ with $f(F)=z$. Thus $B_z^+-B\sim_{\mathbb R}0$ over a neighborhood of $z$, so $B_z^+-B=sf^*z$ for some real number $s$. Since  there exists an lc center $F$ of $(X,B_z^+)$ with $f(F)=z$, $$(X,B_z^++\epsilon f^*z=B+(s+\epsilon)f^*z)$$ is not lc for any positive real number $\epsilon$, so $s=t_z$. Thus 
$$(X/Z\ni z,B+t_zf^*z)=(X/Z\ni z,B_z^+)$$ is a monotonic $q$-complement of $(X/Z\ni z,B)$.
\end{proof}

\noindent\textit{Proof of Proposition \ref{prop: cbf explicit boundary set} continued}. Let $z\in Z$ be a closed point and let $$t_z:=\lct(X/Z\ni z,B;f^*z).$$ By Claim \ref{claim:lctcmpt}, $$(X/Z\ni z,B_z^+:=B+t_zf^*z)$$ is a monotonic $q$-complement of $(X/Z\ni z,B)$. Let $S$ be a component of $f^*z$, $b:=\mult_SB,\, b^+:=\mult_SB_z^+$, and $m:=\mult_Sf^*z$. Then $m$ is a positive integer, $b^+=b+t_zm$, so $t_z=\frac{b^+-b}{m}$. We let $r:=qb^+$, then $r$ is a non-negative integer and $r\leq q$. We have
$$\mult_zB_Z=1-t_z=1-\frac{r-qb}{qm}\in [0,1].$$
If $b\in [1-\delta,1)$ (resp. $(1-\delta,1)$), then 
$$1\geq 1-\frac{r-qb}{qm}\geq\text{(resp. }>\text{)}1-\frac{r-q+q\delta}{qm}\geq 1-\frac{\delta}{m}\geq 1-\delta,$$
so $\mult_zB_Z\in [1-\delta,1)$ (resp. $(1-\delta,1)$). If $b\in\Phi_p$, then $b=1-\frac{k}{pl}$ for some non-negative integer $k\leq p$ and positive integer $l$. Then 
$$0\leq pql(b^+-b)=rpl+kq-pql=pl(r-q)+kq\leq kq\leq pq,$$ 
so
$$\mult_zB_Z=1-\frac{1}{m}(b^+-b)=1-\frac{1}{ml}\cdot\frac{rpl+kq-pql}{pq}\in\Phi_{pq}.$$

\noindent\textbf{Step 3}. In this step we show that $qM_Z$ is a Weil divisor and conclude the proof. Let $V\subset Z$ be a non-empty open subset such that $\Supp B_Z\subset Z\backslash V$ and $q(K+B)\sim 0$ over $V$. Let
$$\Theta:=B+\sum_{z\in Z\backslash V}t_zf^*z\text{ and }\Theta_Z:=B_Z+\sum_{z\in Z\backslash V}t_zz,$$
where $t_z:=\lct(X/Z\ni z,B;f^*z)$. Then $\Theta_Z$ is a reduced divisor, $(Z,B_Z+M_Z)$ is a g-pair induced by the canonical bundle formula of $f: (X,B)\rightarrow Z$, and $(Z,\Theta_Z+M_Z)$ is a g-pair induced by the canonical bundle formula of $f: (X,\Theta)\rightarrow Z$. By Claim \ref{claim:lctcmpt}, $q(K_X+\Theta)\sim 0$ over any point $z\in Z$. By \cite[Lemma 2.4]{Bir19}, $q(K_X+\Theta)\sim_Z0$. Since
\begin{align*}
    q(K_X+\Theta)&=q(K_X+B)+q(\Theta-B)\\
    &\sim qf^*(K_Z+B_Z+M_Z)+qf^*(\Theta_Z-B_Z)\\
    &=qf^*(K_Z+\Theta_Z+M_Z),
\end{align*}
$q(K_Z+\Theta_Z+M_Z)$ is Cartier. Since $K_Z+\Theta_Z$ is a Weil divisor, $qM_Z$ is a Weil divisor.
\end{proof}

\begin{lem}\label{lem: cls lc center}
    Let $(X,B+M)$ be an lc generalized pair, $f: X\rightarrow Z$ a contraction such that $K_X+B+M\sim_{\mathbb R,Z}0$, and $S^\nu$ the normalization of a component $S$ of $\lfloor B\rfloor$ such that $S$ is horizontal over $Z$. Let $f_S: S^\nu\rightarrow Z$ the induced projective surjective morphism, and
    $$S^\nu\xrightarrow{\pi} T\xrightarrow{\tau}Z$$
    the Stein factorization of $f_S$. We let $(S^\nu,B_S+M_S)$ be the generalized pair induced by adjunction
    $$K_{S^\nu}+B_S+M_S:=(K_X+B+M)|_{S^\nu}.$$
    Then for any lc center $W$ of $(X,B+M)$ any irreducible component $V$ of $\tau^{-1}(f(W))$, there exists an lc center $W_S$ of $(S^\nu,B_S+M_S)$ such that $V=\pi(W_S)$.
\end{lem}
\begin{proof}
    We let $g: Y\rightarrow X$ be a dlt modification of $(X,B+M)$, $K_Y+B_Y+M_Y:=g^*(K_X+B+M)$, and $S_Y:=g^{-1}_*S$. Then the Stein factorization of the induced morphism $S_Y\rightarrow Z$ factorizes through $S^\nu$.  Therefore, we may replace $(X,B+M)$ with $(Y,B_Y+M_Y)$ and $S$ with $S_Y$, and assume that $(X,B+M)$ is dlt. The lemma follows from \cite[Lemma 4.3.11]{CHLX23} (see \cite[Lemma 3.19(3)]{LX23} for the NQC case).
\end{proof}

\begin{thm}\label{thm: cbf preserved under adjunction}
    Let $(X,B+M)$ be an lc generalized pair, $f: X\rightarrow Z$ a contraction such that $K_X+B+M\sim_{\mathbb R,Z}0$, and $S^\nu$ the normalization of a component $S$ of $\lfloor B\rfloor$ such that $S$ is horizontal over $Z$. Let $f_S: S^\nu\rightarrow Z$ the induced projective surjective morphism, and
    $$S^\nu\xrightarrow{\pi} T\xrightarrow{\tau}Z$$
    the Stein factorization of $f_S$. We let $(S^\nu,B_S+M_S)$ be the generalized pair induced by adjunction
    $$K_{S^\nu}+B_S+M_S:=(K_X+B+M)|_{S^\nu},$$
    $(Z,B_Z+M_Z)$ a generalized pair induced by the canonical bundle formula of $f: (X,B+M)\rightarrow Z$, and $(T,B_T+M_T)$ a generalized pair induced by the canonical bundle formula of $\pi: (S^\nu,B_S+M_S)\rightarrow T$.  Let $M^Z$ and $M^T$ be the nef parts of $(Z,B_Z+M_Z)$ and $(T,B_T+M_T)$ respectively, whose traces on $Z$ and $T$ are $M_Z$ and $M_T$ respectively. Let $R$ be the ramification divisor of $\tau$. Then:
    \begin{enumerate}
        \item $B_Z=\frac{1}{\deg\tau}\tau_*(R+B_T)$.
        \item For suitable choices of $M^Z$ and $M^T$,
    $M^Z=\frac{1}{\deg\tau}\tau_*M^T.$
    \item For any $\Rr$-divisor $B_T^+\geq B_T$ on $T$ such that $(T,B_T^++M_T)$ is lc and $K_T+B_T^++M_T\sim_{\mathbb R,Z}0$, $(Z,B_Z+\frac{1}{\deg\tau}(B_T^+-B_T)+M_Z)$ and $(X,B+\frac{1}{\deg\tau}f^*\tau_*(B_T^+-B_T)+M)$ are lc.
    \end{enumerate}
\end{thm}
\begin{proof}
First we prove (1). To prove this, we only need to show that $B_Z=\frac{1}{\deg\tau}\tau_*(R+B_T)$ near the generic point of any prime divisor $D$ on $Z$. Thus possibly fixing a prime divisor $D$ on $Z$ and shrinking $Z$ to a neighborhood of the generic point of $D$, we may assume that $Z$ is smooth, $\Supp B_Z=D$, and 
$\mult_DB_Z=1-t$, where
$$t:=\lct(X,B+M;f^*D).$$
We let $B':=B+tf^*D$ and let
$$K_{S^\nu}+B_S'+M_S:=(K_X+B'+M)|_{S^\nu}.$$
Since $(X,B'+M)$ is lc, $(S^\nu,B_S'+M_S)$ is lc. Moreover,
$$B_S'=B_S+tf^*D|_{S^\nu}=B_S+tf_S^*D=B_S+t\pi^*(\tau^*D).$$
Therefore, $(T,B_T'+M_T)$ is the generalized pair induced by the canonical bundle formula of $\pi: (S^\nu,B_S'+M_S)\rightarrow T$, where $B_T':=B_T+\tau^*D$.

By construction, $(Z,D+M_Z)$ is the lc g-pair induced by the canonical bundle formula of $f: (X,B'+M)\rightarrow Z$, and $D$ is the image of an lc center of $(X,B'+M)$ on $Z$. By Lemma \ref{lem: cls lc center}, any irreducible component of $\tau^{-1}(D)$ is the image of an lc center of $ (S^\nu,B_S'+M_S)$ on $T$. Therefore, any irreducible component of $\tau^{-1}(D)$ is an lc center of $(T,B_T'+M_T)$.

We let $D_1,\dots,D_n$ be the irreducible components of $\tau^{-1}(D)$ and let $r_i$ be the ramification index of $\tau$ along $D_i$ for each $i$. Then $\sum_{i=1}^nr_i\deg(D_i\rightarrow D)=\deg\tau$ and $B_T'=\sum_{i=1}^nD_i$ over the generic point of $D$, so
\begin{align*}
   \frac{1}{\deg\tau}\tau_*(R+B_T')&=\frac{1}{\deg\tau}\tau_*\sum_{i=1}^n((r_i-1)D_i+D_i)=\frac{1}{\deg\tau}\sum_{i=1}^nr_i\tau_*D_i\\
   &=\frac{1}{\deg\tau}\sum_{i=1}^nr_i\deg(D_i\rightarrow D)=D. 
\end{align*}
Since $B_T'=B_T+t\tau^*D$, we have
$$B_Z=D-tD=\frac{1}{\deg\tau}\tau_*(R+B_T')-\frac{1}{\deg\tau}\tau_*(t\tau^*D)=\frac{1}{\deg\tau}\tau_*(R+B_T).$$

Now we prove (2). To prove this, we only need to show that for any birational morphisms $h_Z: Z'\rightarrow Z$, we have
$$M_{Z'}=\frac{1}{\deg\tau}\tau'_*M_{T'},$$
for suitable choices of $M^Z$ and $M^T$, where $T'$ is the main component of $T\times_ZZ'$, $\tau': T'\rightarrow Z'$ is the induced birational morphism, and $M_{Z'}$ and $M_{T'}$ are the traces of $M^Z$ and $M^T$ on $Z'$ and $T'$ respectively.  By \cite[Theorem 2.1 and Proposition 4.4]{AK00} (see also \cite[Theorem B.6]{Hu20}, \cite[Theorem 2.2]{ACSS21}), possibly replacing $Z'$ with a higher model, we may assume that there exist a toroidal g-pair  $(X',\Sigma_{X'})$, a log smooth pair $(Z',\Sigma_{Z'})$, and a commutative diagram
 \begin{center}$\xymatrix{
X'\ar@{->}[r]^{h}\ar@{->}[d]_{f'}& X\ar@{->}[d]^{f}\\
Z'\ar@{->}[r]^{h_Z} & Z\\
}$
\end{center}
satisfying the following.
\begin{itemize}
\item $h$ and $h_Z$ are projective birational morphisms.
\item $f': (X',\Sigma_{X'})\rightarrow (Z',\Sigma_{Z'})$ is a toroidal contraction.
\item $\Supp(h^{-1}_*B)\cup\Supp\Exc(h)$ is contained in $\Supp\Sigma_{X'}$.
\item $X'$ has at most toric quotient singularities.
\item $f'$ is equi-dimensional.
\item $M^X$ descends to $X'$ and $M^Z$ descend to $Z'$. We let $M'$ be the trace of $M^X$ on $X'$.
\item $X'$ is $\Qq$-factorial klt.
\end{itemize}
Let $B':=h^{-1}_*B+\Supp\Exc(h)$. Then $(X,B'+M')$ is $\Qq$-factorial lc,
$$\kappa_{\sigma}(X'/Z',K_{X'}+B'+M')=\kappa_{\sigma}(X'/Z,K_{X'}+B'+M')=\kappa_{\sigma}(X/Z,K_{X}+B+M)=0,$$
and
$$\kappa_{\iota}(X'/Z',K_{X'}+B'+M')=\kappa_{\iota}(X'/Z,K_{X'}+B'+M')=\kappa_{\iota}(X/Z,K_{X}+B+M)=0.$$
Therefore, we have
$$K_{X'}+B'+M'\sim_{\mathbb R,Z'}E^h+E^v$$
where $E^h\geq 0$, any component of $E^h$ is horizontal$/Z'$, and $E^v\geq 0$ is vertical$/Z$. Let $$\nu_P:=1-\sup\{t\mid E^\nu-tf'^*P\geq 0\}$$
for any prime divisor $P$ on $Z'$. Since $Z'$ is smooth, $\nu_P$ is well-defined for any prime divisor $P$ on $Z'$. Since $f'$ is equi-dimensional, possibly replacing $E^v$ with $E-\sum_P\nu_Pf'^*P$, we may assume that $E^v$ is very exceptional over $Z'$.

We run a $(K_{X'}+B'+M')$-MMP$/Z'$ with scaling of an ample divisor $A$. Then after finitely many steps, we get an birational map $X'\dashrightarrow X''$ such that $K_{X''}+B''+M''$ is movable$/Z'$, where $B''$ is the image of $B'$ on $X''$ and $M''$ is the image of $M'$ on $X''$. Thus for any general fiber $F$ of $X''\rightarrow Z'$, $(K_{X''}+B''+M'')|_F$ is movable. Since $\kappa_{\sigma}((K_{X''}+B''+M'')|_F)=0$, $(K_{X''}+B''+M'')|_F\equiv 0$. 

Let $(E^h)''$ and $(E^v)''$ be the images of $E^h$ and $E^v$ on $X''$ respectively. Then $(E^h)''|_F\equiv 0$, so $(E^h)''|_F=0$, hence $(E^h)''=0$. Since $E^v$ is very exceptional over $Z'$, 
$K_{X''}+B''+M''\sim_{\mathbb R,Z'}(E^v)''$ is very exceptional over $Z'$. By \cite[Proposition 3.9]{HL22}, we may assume that the  $(K_{X'}+B'+M')$-MMP$/Z'$ terminates at $X''$ and $K_{X''}+B''+M''\sim_{\mathbb R,Z'}0$. 

By applying the negativity lemma twice, we know that $(X'',B''+M'')$ and $(X,B+M)$ are crepant over the generic point of $Z$. Therefore, $M^Z$ is the moduli part of the canonical bundle formula of $(X'',B''+M'')\rightarrow Z'$.

Let $S':=h^{-1}_*S$ and let $S''$ be the image of $S'$ on $X''$. Since the induced birational map $X'\dashrightarrow X''$ only contract divisors that are contained in $\Supp(E^h+E^v)$, $S''\not=0$, and $X'\dashrightarrow X''$ is an isomorphism near the generic point of $S'$. Therefore, the Stein factorization of the induced morphism $S''^\nu\rightarrow Z'$ factorizes through $T'$, where $S''^\nu$ is the normalization of $S''$ (actually $S''$ is normal, but we do not need this fact). Let
$$K_{S''^\nu}+B_{S''}+M_{S''}:=(K_{X''}+B''+M'')|_{S''^\nu},$$
then $(S''^\nu,B_{S''}+M_{S''})$ and $(S^\nu,B_S+M_S)$ are crepant over the generic point of $Z$, so  $(S''^\nu,B_{S''}+M_{S''})$ and $(S^\nu,B_S+M_S)$ are crepant over the generic point of $T$. Therefore, $M^T$ is the moduli part of the canonical bundle formula of $(S''^\nu,B_{S''}+M_{S''})\rightarrow Z'$.

We let $B_{Z'}$ and $B_{T'}$ be the discriminant parts of $(X'',B''+M'')\rightarrow Z'$ and $(S''^\nu,B_{S''}+M_{S''})\rightarrow Z'$ respectively. By (1), we have
$$B_{Z'}=\frac{1}{\deg\tau'}\tau'_*(R'+B_{T'})$$
where $R'$ is the ramification divisor of $\tau'$. For suitable choices of $M^T$ and $M^{Z}$, we have
$$K_{T'}+B_{T'}+M_{T'}=\tau'^*(K_{Z'}+B_{Z'}+M_{Z'}).$$
Therefore, $M_{Z'}=\frac{1}{\deg\tau}\tau'_*M_{T'}$ and we are done.

Finally, we prove (3). Let $B_Z^+:=B_Z+\frac{1}{\deg\tau}(B_T^+-B_T)$ and $B^+:=B+f^*(B_Z^+-B_Z)$. By (1)(2), $(Z,B_Z^++M_Z)$ is the g-pair induced by the canonical bundle formula of $\tau: (T,B_T^++M_T)\rightarrow Z$, so $(Z,B_Z^++M_Z)$ is lc. Since  $(Z,B_Z^++M_Z)$ is also the g-pair induced by the canonical bundle formula of $f: (X,B^++M)\rightarrow Z$, $(X,B^++M)$ is lc.
\end{proof}

\begin{rem}
In this paper, we will apply Theorem \ref{thm: cbf preserved under adjunction} and Lemma \ref{lem: cls lc center} only for the NQC case. Nevertheless, we prove Theorem \ref{thm: cbf preserved under adjunction} and Lemma \ref{lem: cls lc center} in full generality for potential future applications.
\end{rem}

\subsection{A lemma on  Cartier index}

\begin{lem}\label{lem: cartier index}
    Let $I$ be a positive integer, $(X,B)$ an lc pair such that $B$ is a reduced divisor, $S$ a component of $\lfloor B\rfloor$, and $P$ a prime divisor on $S$. Assume that $I(K_X+B)|_S$ is Cartier near the generic point of $P$. Then:
    \begin{enumerate}
        \item $I(K_X+B)$ or $2I(K_X+B)$ is Cartier near the generic point of $P$.
        \item If $(X,S)$ is plt near the generic point of $P$, then  $I(K_X+B)$ is Cartier near the generic point of $P$.
    \end{enumerate}
\end{lem}
\begin{proof}
Possibly cutting $X$ by general hyperplane sections, we may assume that $X$ is a surface and $P$ is a closed point. Since $(X,B)$ is lc, there are two cases:

\medskip

\noindent\textbf{Case 1}. $S$ is singular near $P$. In this case, near a neighborhood of $P$, $B=S$, $P$ is a node of $S$, and $P$ is a cyclic quotient singularity of $X$. In this case, $K_X+S$ is Cartier near $P$, so $I(K_X+B)$ is Cartier near $P$.

\medskip

\noindent\textbf{Case 2}. $S$ is non-singular near $P$, and $(X,B)$ is plt near $P$. In this case, $B=S$ near $P$, and 
$$(K_X+B)|_S=K_S+\frac{m-1}{m}P$$
near a neighborhood of $P$, where $m$ is the Cartier index of $K_X+S$ at $P$. Since $I(K_X+B)|_S$ is Cartier near $P$, $m\mid I$. Since $m(K_X+B)$ is Cartier near $P$, $I(K_X+B)$ is Cartier near $P$.

\medskip

\noindent\textbf{Case 3}. $S$ is non-singular near $P$, and $(X,B)$ is not plt near $P$. There are two sub-cases:

\smallskip

\noindent\textbf{Case 3.1}. $(X,S)$ is plt near $P$. In this case, $(X,B)$ is log toric near $P$, so $K_X+B$ is Cartier near $P$, hence $I(K_X+B)$ is Cartier near $P$.

\smallskip

\noindent\textbf{Case 3.2}. $(X,S)$ is not plt near $P$. In this case, $P$ is a D-type singularity of $X$ and $B=S$ near $P$. Then $2(K_X+B)$ is Cartier near $P$.
\end{proof}

\section{Curves}\label{sec: curve}

In this section we study explicit bounds of values given in Definition \ref{defn: special explicit values} when $d=1$.

\subsection{Bound of the lc threshold}

\begin{prop}\label{prop: lct 1-gap dim 1}
Let $p$ be a positive integer. Then $\delta_{\lct}(1,p)=\min\{\frac{1}{p},\frac{1}{2}\}$.
\end{prop}
\begin{proof}
First, we consider the triple 
$$\left(X:=\mathbb P^1,B=\left(\min\left\{\frac{1}{p},\frac{1}{2}\right\}\right)P;D:=P\right)$$ where $P$ is a closed point on $\mathbb P^1$. Then $B\in\Phi_p$, $D\in\mathbb N^+$, and 
$$\lct(X,B;D)=\max\left\{\frac{p-1}{p},\frac{1}{2}\right\}.$$
Thus $\delta_{\lct}(1,p)\leq\min\left\{\frac{1}{p},\frac{1}{2}\right\}$.

Suppose that $\delta_{\lct}(1,p)<\min\left\{\frac{1}{p},\frac{1}{2}\right\}$. Then there exists a pair $(X,B)$ of dimension $1$ such that $B\in\Phi_p$, and an an effective Weil divisor $D$ on $X$, such that $$\max\left\{\frac{1}{2},\frac{p-1}{p}\right\}<t:=\lct(X,B;D)<1.$$ Thus there exists a component $C$ of $D$, such that $\mult_CD>0$ and $\mult_C(B+tD)=1$. Therefore,
$$\gamma+td=1$$
for some $\gamma\in\Phi_p$ and $d\in\mathbb N^+$. If $d\geq 2$, then
$$t=\frac{1-\gamma}{d}\leq\frac{1}{2},$$
a contradiction. Thus $d=1$. We write $\gamma=1-\frac{k}{pn}$ where $n\in\mathbb N^+$ and $k\in\mathbb N\cap [0,p]$. Then
$$t=1-\gamma=\frac{k}{pn}.$$
If $n\geq 2$, then 
$$t=\frac{k}{pn}\leq\frac{1}{n}\leq\frac{1}{2},$$ 
a contradiction. Thus  $n=1$, and $t=\frac{k}{p}$. Since $t<1$, $k\leq p-1$, so $t\leq\frac{p-1}{p}$, a contradiction. Therefore, $\delta_{\lct}(1,p)\geq\min\left\{\frac{1}{p},\frac{1}{2}\right\}$ and we are done.
\end{proof}

\subsection{Bound of the global lc threshold and exceptional mlds}

\begin{lem}\label{lem: elementary computation glct 1-gap dim 1}
Let $p$ be a positive integer and $\delta:=\min\left\{\frac{1}{6},\frac{1}{p(p+1)}\right\}$. Then the equation
$$2=\sum_{i=1}^n\gamma_i+\sum_{j=1}^mb_j, \gamma_i\in\Phi_p\backslash\{0\}, b_j\in (1-\delta,1), n\in\mathbb N, m\in\mathbb N^+$$
does not have a solution.
\end{lem}
\begin{proof}
Suppose that the above equation has a solution. Since $\delta<\frac{1}{3}$, $m\leq 2$ and $n\geq 1$.

First suppose that $m=2$. Then
$$\gamma_1<2-2(1-\delta)=2\delta\leq\min\{\frac{1}{3},\leq\frac{2}{p(p+1)}\}.$$
If $p=1$, then $\gamma_1\leq\frac{1}{3}$, which is not possible. If $p>1$, then $\gamma_1<\frac{1}{p}$, which is also not possible. Therefore, $m=1$, and we have
$$\sum_{i=1}^n\gamma_i\in (1,1+\delta).$$
In particular, $N\sum_{i=1}^n\gamma_i$ is not an integer for any $N\leq\max\{6,p(p+1)\}$. We may write $$\gamma_i=\frac{n_i-1+\frac{m_i}{p}}{n_i}$$ for each $i$, where $n_i\in\mathbb N^+,m_i\in\mathbb N$, and $0\leq m_i\leq p$. Possibly reordering indices, we may assume that $n_i\geq n_j$ for any $i\leq j$.
 
 If $n_1=1$, then $n_i=1$ for any $i$. Thus $\gamma_i=\frac{m_i}{p}$ for any $i$. Thus $p\sum_{i=1}^n\gamma_i$ is an integer, which is not possible. If $n_3\geq 2$, then $$\sum_{i=1}^n\gamma_i\geq\frac{3}{2}>1+\delta,$$ which is not possible. Thus $n_i=1$ for any $i\geq 3$.
 
We have the following two cases.

\medskip

\noindent\textbf{Case 1}. $n_2\geq 2$. If $n_1\geq 3$, then $$\sum_{i=1}^n\gamma_i\geq\frac{2}{3}+\frac{1}{2}\geq\frac{7}{6}\geq 1+\delta,$$ which is not possible. Thus $n_1=2$, so $n_2=2$. We have
$$\sum_{i=1}^n\gamma_i=1+\frac{m_1+m_2}{2p}+\frac{\sum_{i\geq 3}m_i}{p},$$
so $2p\sum_{i=1}^n\gamma_i$ is an integer, which is not possible.

\medskip

\noindent\textbf{Case 2}. $n_2=1$. Then
$$\sum_{i=1}^n\gamma_i=\frac{n_1-1+\frac{m_1}{p}}{n_1}+\frac{\sum_{i\geq 2}m_i}{p}\in (1,1+\delta).$$
In particular, $\sum_{i\geq 2}m_i\geq 1$. If $n_1\leq p$, then $n_1p\sum_{i=1}^n\gamma_i$ is an integer, which is not possible. Thus $n_1>p$, so $\gamma_1\geq\frac{p}{p+1}$. Thus 
$$\sum_{i=1}^n\gamma_i\geq\frac{p}{p+1}+\frac{1}{p}=1+\frac{1}{p(p+1)}\geq 1+\delta,$$
which is not possible.
\end{proof}

\begin{prop}\label{prop: glct 1-gap dim 1}
Let $p$ be a positive integer. Then $\delta_{\glct}(1,p)=\min\left\{\frac{1}{6},\frac{1}{p(p+1)}\right\}$.
\end{prop}
\begin{proof}
Let $(X,B+tS)$ be an lc pair such that $t>0$, $\dim X=1,B\in\Phi_p$, $0\not=S\in\mathbb N^+$, and $K_X+B+tS\equiv 0$. We write $B=\sum \gamma_iB_i$, such that $\gamma_i\in\Phi_p$ and $B_i$ are prime divisors. Then either $t=0, B=0$, and $X$ is an elliptic curve, or $X=\mathbb P^1$ and
$$2=\sum_i\gamma_i+t\deg S.$$
By Lemma \ref{lem: elementary computation glct 1-gap dim 1}, $t\leq 1-\min\left\{\frac{1}{6},\frac{1}{p(p+1)}\right\}$. Thus $\delta_{\glct}(1,p)\geq\min\left\{\frac{1}{6},\frac{1}{p(p+1)}\right\}$.
The examples
$$\left(\mathbb P^1,\frac{1}{2}B_1+\frac{2}{3}B_2+\frac{5}{6}S\right)\text{ and }\left(\mathbb P^1,\frac{p}{p+1}B_1+\frac{1}{p}B_2+\frac{p(p+1)-1}{p(p+1)}S\right)$$
show that $\delta_{\glct}(1,p)\leq\min\{\frac{1}{6},\frac{1}{p(p+1)}\}$ and the proposition follows.
\end{proof}

\begin{prop}
Let $p$ be a positive integer. Then $$\delta_{\mld}(1,p)=\delta_{\mld,\ft}(1,p)=\min\left\{\frac{1}{6},\frac{1}{p(p+1)}\right\}.$$
\end{prop}
\begin{proof}
Let $(X,B)$ be an exceptional pair such that $B\in\Phi_p$ and $\dim X=1$. Then either $B=0,\mld(X,B)=1$ and $X$ is an elliptic curve, or $X=\mathbb P^1$. Thus $\delta_{\mld}(1,p)=\delta_{\mld,\ft}(1,p)$. By considering the examples $$\left(\mathbb P^1,\frac{1}{2}B_1+\frac{2}{3}B_2+\frac{5}{6}B_3\right)\text{ and }\left(\mathbb P^1,\frac{p}{p+1}B_1+\frac{p(p+1)-1}{p(p+1)}B_2+\frac{1}{p}B_3\right),$$
the proposition follows from Propositions \ref{prop: compare exc mld and glct} and \ref{prop: glct 1-gap dim 1}.
\end{proof}

\subsection{Bound of complements and indices}

\begin{prop}\label{prop: i11=6}
$I(1,1)=6$.
\end{prop}
\begin{proof}
The lemma follows from the fact that any pair in dimension $1$ with standard coefficients has an $N$-complement for some $N\in\{1,2,3,4,6\}$ (cf. \cite[19.4 Theorem]{Kol+92}).
\end{proof}

\begin{prop}\label{prop: n complemntary bound dim 1}
Let $p$ be a positive integer. Then
$$I(1,p)\leq N(1,p)=N_{\ft}(1,p)\leq N_g(1,p)\leq \max\{6,p\cdot (p(p+1))!\}.$$
\end{prop}
\begin{proof}
By Proposition \ref{prop: i11=6}, we may assume that $p\geq 2$. For any lc pair $(X,B)\in\mathcal{S}_{\rc}(1,p)\backslash\mathcal{S}_{\rc,\ft}(1,p)$, $X$ is an elliptic curve, $B=0$, and $(X,0)$ is a $1$-complement of itself. Thus $N(1,p)=N_{\ft}(1,p)$.

By Lemma \ref{lem: elementary comparison of explicit bound}, we only need to prove that $N_g(1,p)\leq p\cdot (p(p+1))!$. 

For any $(X,B+M)\in\mathcal{S}_{g,\rc}(1,p)$, $X\cong\mathbb P^1$. We write $B=\sum_{i=1}^n b_iB_i$, where $B_i$ are distinct points of $\mathbb P^1$, $b_i\in\Phi_p\backslash\{0\}$, and let $m:=p\deg M$. We let $$D:=\sum_{b_i>1-\frac{1}{p(p+1)}}B_i\text{ and }C:=B-B\wedge D.$$ 
If $(X,C+D+M)$ is not $\Rr$-complementary, then there exists $t\in\left(1-\frac{1}{p(p+1)},1\right)$ such that $$K_X+C+tD+M\equiv 0.$$ This contradicts Proposition \ref{prop: glct 1-gap dim 1}. Thus $(X,C+D+M)$ is $\Rr$-complementary. Possibly replacing $B$ with $C+D$, we may assume that $$B\in\Phi_p\cap\left(\left[0,1-\frac{1}{p(p+1)}\right]\cup\{1\}\right).$$  Since $B\in\Phi_p$, for any $i$, $pIb_i\in\mathbb N$ for some $I\leq p(p+1)$. Thus $$p\cdot (p(p+1))!(K_X+B+M)$$ is Cartier, so $(X,B+M)$ is a monotonic $(p\cdot (p(p+1))!)$-complement of itself. Thus $$N_g(1,p)\leq p\cdot (p(p+1))!$$ and we are done.
\end{proof}

\section{Explicit bounds on surfaces}\label{sec: surface standard coefficient}

In this section we provide explicit bounds of values given in Definition \ref{defn: special explicit values} when $d=2$ and $p=1$. We refer to Section \ref{sec: Non-standard coefficient case} for some results when $d=2$ and $p$ is arbitrary.

\subsection{Bound of the lc threshold}

\begin{prop}\label{prop: lct 1-gap dim 2}
Let $p$ be a positive integer. Then
$\delta_{\lct}(2,p)=\min\left\{\frac{1}{6},\frac{1}{p(p+1)}\right\}$.
\end{prop}
\begin{proof}
It follows from Lemma \ref{lem: equality glct d and lct d+1} and Proposition \ref{prop: glct 1-gap dim 1}.
\end{proof}

\subsection{Bound of complements and indices: the standard coefficient case}

In this subsection, we estimate $I(2,1)$ and $N(2,1)$.
\begin{thm}\label{thm: exceptional 1/42 lc global strong}
$\delta_{\mld,\ft}(2,1)=\delta_{\mld}(2,1)=\delta_{\glct}(2,1)=\frac{1}{42}$.
\end{thm}
\begin{proof}
By considering the example
$$\left(\mathbb P^2,\frac{1}{2}L_1+\frac{2}{3}L_2+\frac{6}{7}L_3+\frac{41}{42}L_4\right)$$
where $L_1,L_2,L_3,L_4$ are four lines on $\mathbb P^2$ in general position, we have $\delta_{\mld,\ft}(2,1)\leq\frac{1}{42}$, $\delta_{\mld}(2,1)\leq\frac{1}{42}$, and $\delta_{\glct}(2,1)\leq\frac{1}{42}$. By Proposition \ref{prop: compare exc mld and glct}, we only need to show that $\delta_{\glct}(2,1)\geq\frac{1}{42}$.

Suppose that $\delta_{\glct}(2,1)<\frac{1}{42}$. Then there exists a projective log Calabi-Yau lc surface pair $(X,B+tS)$ such that $B\in\Phi_1$, $S$ is a non-zero reduced divisor, and $t\in (\frac{41}{42},1)$. Possibly replacing $(X,B+tS)$ with a dlt model we may assume that $(X,B+tS)$ is $\Qq$-factorial dlt. We let $S_1$ be a component of $S$ and run a $(K_X+B+t(S-S_1))$-MMP, which terminates with a Mori fiber space $f: Y\rightarrow Z$. Let $B_Y,S_Y,S_{1,Y}$ be the images of $B,S$ and $S_1$ on $Y$ respectively, then $S_{1,Y}$ is horizontal$/Z$. There are two cases:

\medskip

\noindent\textbf{Case 1}. $\dim Z=1$. In this case, we let $F$ be a general fiber of $f$, $B_F:=B_Y|_F$, and $S_F:=S_Y|_F$. Since $S_{1,Y}|_F\not=0$, $S_F\not=0$. Thus $(F,B_F+tS_F)$ is an lc pair of dimension $1$ such that $B_F\in\Phi_1$, $S_F$ is a non-zero reduced divisor, and $t\in (\frac{41}{42},1)$. Thus $\delta_{\glct}(1,1)<\frac{1}{42}<\frac{1}{6}$, which contradicts Proposition \ref{prop: glct 1-gap dim 1}.

\medskip

\noindent\textbf{Case 2}. $\dim Z=0$. In this case, $\rho(Y)=1$. We let $S_{1,Y},\dots,S_{n,Y}$ be the irreducible components of $S_Y$. Since $K_Y+B_Y+\frac{41}{42}S_Y$ is anti-ample and $K_Y+B_Y+S_Y$ is ample, and since $\rho(Y)=1$, there exists an index $1\leq j\leq n$ and $t'\in (\frac{41}{42},1)$, such that 
$$K_Y+B_Y+\sum_{i=1}^{j-1}S_{i,Y}+t'S_{j,Y}+\sum_{i=j+1}^n\frac{41}{42}S_{i,Y}\equiv 0.$$
We let $B_Y':=B_Y+\sum_{i=1}^{j-1}S_{i,Y}+\sum_{i=j+1}^n\frac{41}{42}S_{i,Y}$ and $S_Y':=S_{j,Y}$. Then $B_Y'\in\Phi_1$, $S_Y'$ is a prime divisor, and $K_Y+B_Y'+t'S_Y'\equiv 0$. Since $(Y,B_Y+tS_Y)$ is lc, by Proposition \ref{prop: lct 1-gap dim 2}, $(Y,B_Y+S_Y)$ is lc. Since $B_Y+S_Y\geq B_Y'+t'S_Y'$, $(Y,B_Y'+t'S_Y')$ is lc. This contradicts \cite[5.3 Theorem]{Kol94}.
\end{proof}

\subsubsection{Non-exceptional case}

\begin{prop}\label{prop: n21 nonexceptional case}
Let $(X/Z\ni z,B+B')$ be an $\Rr$-complementary projective lc surface pair such that $B\in\Phi_1$, $B'\in [\frac{6}{7},1)$, and
\begin{itemize}
    \item either $Z$ is a point and  $(X,B+B')$ is not exceptional, or
    \item $\dim Z\geq 1$.
\end{itemize}
Then $(X/Z\ni z,B+B')$ has a monotonic $N$-complement for some $N\in\{1,2,3,4,6\}$ that is not klt. Moreover, if $(X,B+B')$ has infinitely many lc centers over $z$, then $(X/Z\ni z,B+B')$ has a monotonic $2$-complement.
\end{prop}
\begin{proof}
Possibly shrinking $Z$ to a neighborhood of $z$, we may assume that $(X/Z,B+B')$ is $\Rr$-complementary. Then there exists an $\Rr$-complement $(X/Z,\tilde B:=B+B'+G)$ of $(X/Z,B+B')$ that is not klt. Let $f: Y\rightarrow X$ be a dlt modification of $(X,\tilde B)$, $K_Y+\tilde B_Y:=f^*(K_X+\tilde B)$, $S_Y:=\lfloor \tilde B_Y\rfloor$, $B_Y:=f^{-1}_*B-f^{-1}_*B\wedge S_Y$, $B'_Y:=f^{-1}_*B'-f^{-1}_*B'\wedge S_Y$, and $G_Y:=f^{-1}_*G-f^{-1}_*G\wedge S_Y$. Since
$$-\epsilon(K_Y+S_Y+B_Y+B_Y')\sim_{\mathbb R,Z}K_Y+S_Y+B_Y+B_Y'+(1+\epsilon)G_Y,$$
we may run a $-(K_Y+S_Y+B_Y+B_Y')$-MMP$/Z$ which terminates with a model $Y'$, such that $-(K_{Y'}+S_{Y'}+B_{Y'}+B_{Y'}')$ is semi-ample$/Z$, where $S_{Y'},B_{Y'},B'_{Y'}$ are the strict transforms of $S_Y,B_Y$, $B_{Y'}$ on $Y'$ respectively. By \cite[2.3 Inductive Theorem]{Sho00},  $$(Y'/Z\ni z,S_{Y'}+B_{Y'}+B'_{Y'})$$ has a monotonic $N$-complement $(Y'/Z\ni z,B_{Y'}^+)$ for some $N\in\{1,2,3,4,6\}$. Moreover, if $(X,B+B')$ has infinitely many lc centers over $z$, then by \cite[2.3.2]{Sho00},
$$(Y'/Z\ni z,S_{Y'}+B_{Y'}+B'_{Y'})$$
has a monotonic $2$-complement.

Let $g: Y\rightarrow Y'$ be the induced morphism and let $$K_X+B^+:=f_*g^*(K_{Y'}+B_{Y'}^+).$$ Then $(X/Z\ni z,B^+)$ is a monotonic $N$-complement of $(X/Z\ni z,B+B')$ that is not klt. Moreover, if $(X,B+B')$ has infinitely many lc centers over $z$, then $(X/Z\ni z,B^+)$ is a monotonic $2$-complement of $(X/Z\ni z,B+B')$ that is not klt.
\end{proof}

\subsubsection{Kodaira dimension zero case}

\begin{prop}[cf. {\cite[Proposition 3.5]{ETW22}}]\label{prop: klt standard coefficient surface 66}
$I(2,1)=66$. 
\end{prop}

\begin{prop}\label{prop: n21 kod 0 case}
Let $(X,B+B')$ be an exceptional surface pair such that $B\in\Phi_1$, $B'\in(\frac{41}{42},1)$, and $\kappa_{\iota}(-(K_X+B+B'))=0$. Then $(X,B+B')$ has a monotonic $N$-complement for some integer $N\leq 66$.
\end{prop}
\begin{proof}
 By Theorems \ref{thm: strong version 1.2 liu22} and \ref{thm: exceptional 1/42 lc global strong}, $B'=0$. We may run a $-(K_X+B)$-MMP $g: X\rightarrow Y$ which terminates with a model $Y$ such that $K_Y+B_Y\equiv 0$, where $B_Y:=g_*B$. By Proposition \ref{prop: klt standard coefficient surface 66}, $N(K_Y+B_Y)\sim 0$ for some positive integer $N\leq 66$. Let $$K_X+B^+:=g^*(K_{Y}+B_Y),$$ 
 then $(X,B^+)$ is a monotonic $N$-complement of $(X,B)$.
\end{proof}

\subsubsection{Kodaira dimension one case}

\begin{prop}\label{prop: n21 kod 1 case}
Let $(X,B+B')$ be an exceptional surface pair such that $B\in\Phi_1$, $B'\in (\frac{41}{42},1)$, and $\kappa_{\iota}(-(K_X+B+B'))=1$. Then $(X,B+B')$ has a monotonic $N$-complement for some integer 
$$N\leq\max_{1\leq q\leq 6}\lcm(q,N_g(1,q))<36\cdot (42)!.$$
\end{prop}
\begin{proof}
 By Theorems \ref{thm: strong version 1.2 liu22} and \ref{thm: exceptional 1/42 lc global strong}, $B'=0$. We may run a $-(K_X+B)$-MMP $g: X\rightarrow Y$ which terminates with a model $Y$ such that $-(K_Y+B_Y)$ is semi-ample, where $B_Y:=g_*B$. Possibly replacing $X,B$ with $Y,B_Y$ respectively, we may assume that $-(K_X+B)$ is semi-ample and defines a contraction $f: X\rightarrow Z$. Then $\dim Z=1$. By Propositions \ref{prop: n21 nonexceptional case} and \ref{prop: cbf explicit boundary set}, there exists an lc g-pair $(Z,B_Z+M_Z)$ induced by the canonical bundle formula $f: (X,B)\rightarrow Z$ and an integer $q\in\{1,2,3,4,6\}$, such that
$$q(K_X+B)\sim qf^*(K_Z+B_Z+M_Z),$$
$qM_Z$ is a Weil divisor, and $B_Z\in\Phi_{q}$. Thus  $(Z,B_Z+M_Z)$ has a monotonic $N$-complement $(Z,B_Z^++M_Z)$ for some positive integer $N\leq N_g(1,q)$, so  $(X,B+f^*(B_Z^+-B_Z))$ is a monotonic $\lcm(q,N)$-complement for some positive integer $N\leq N_g(1,q)$. By Proposition \ref{prop: n complemntary bound dim 1}, $N_g(1,q)\leq 6\cdot (42)!$. The proposition follows.
\end{proof}

\subsubsection{Kodaira dimension two case}

\begin{lem}[{\cite[Lemma 2.7]{Liu23}}; cf. {\cite[Proof of Lemma 3.7]{AM04}, \cite[After Theorem A]{Lai16}, \cite[Lemma 2.2]{Bir23}}]\label{lem: index epsilonlc fano surface}
Let $\epsilon\in\left(0,\frac{1}{\sqrt{3}}\right)$ be a real number and $(X,B)$ a projective $\epsilon$-lc surface pair such that $-(K_X+B)$ is ample. Then for any Weil divisor $D$ on $X$, $ID$ is Cartier for some positive integer $I<\left(\frac{2}{\epsilon}\right)^{\lfloor\frac{128}{\epsilon^5}\rfloor}$.
\end{lem}
\begin{proof}
Let $f: Y\rightarrow X$ be the minimal resolution of $X$. By \cite[Corollary 1.10]{AM04}, $\rho(Y)\leq\frac{128}{\epsilon^5}$. Thus the dual graph of $f$, $\mathcal{D}(f)$, is a tree of rational curves with at most $\lfloor\frac{128}{\epsilon^5}\rfloor-1$ vertices. By \cite[Corollary 2.19]{Ale93}, for any prime $f$-exceptional divisor $E$, $-E^2\leq\frac{2}{\epsilon}$. Thus $\det(f)$, the absolute value of the determinant of the intersection matrix of $\mathcal{D}(f)$ is strictly less than $\left(\frac{2}{\epsilon}\right)^{\lfloor\frac{128}{\epsilon^5}\rfloor}$. We may take $I=\det(f)$.
\end{proof}

\begin{prop}\label{prop: n21 kod 2 case}
Let $(X,B+B')$ be an exceptional surface pair such that $B\in\Phi_1$, $B'\in (\frac{41}{42},1)$, and $\kappa_{\iota}(-(K_X+B+B'))=2$. Let $I_0:=(42)!\cdot 84^{128\cdot 42^5}$. Then 
\begin{enumerate}
    \item $I(K_X+B+B')$ is Cartier for some integer $I<I_0$, and
    \item $(X,B+B')$ has a monotonic $N$-complement for some $N<96I_0$.
\end{enumerate}
\end{prop}
\begin{proof}
(1) By Theorems \ref{thm: strong version 1.2 liu22} and \ref{thm: exceptional 1/42 lc global strong}, $B'=0$, and we may pick $0\leq D\in |-(K_X+B)|_{\mathbb Q}$ such that $(X,B+G)$ is $\frac{1}{42}$-lc. In particular, $(X,B)$ is $\frac{1}{42}$-lc, $$B\in\left\{1-\frac{1}{n}\bigm| n\in\mathbb N^+, n\leq 42\right\},$$ and $(42)!B$ is a Weil divisor. 

Since $\kappa_{\iota}(-(K_X+B+B'))=2$, $X$ is of Fano type. Thus there exists a klt pair $(X,\Delta)$ such that $-(K_X+\Delta)$ is ample. Let $0<\epsilon\ll 1$ be a real number.
Then there exists a positive real number $\delta=\delta(\epsilon)$ such that $$(X,\tilde B:=(1-\delta)(B+G)+\delta\Delta)$$ is $(\frac{1}{42}-\epsilon)$-klt, and 
$$-(K_X+\tilde B)\sim_{\mathbb R}-\delta(K_X+\Delta)$$ is ample. By Lemma \ref{lem: index epsilonlc fano surface}, $ID$ is Cartier for any prime divisor $D$ on $X$ for some 
 integer 
 $$I\leq\left(\frac{2}{\frac{1}{42}-\epsilon}\right)^{\frac{128}{\left(\frac{1}{42}-\epsilon\right)^5}}-1<84^{128\cdot 42^5}$$
since $0<\epsilon\ll 1$. This implies (1) as $(42)!B$ is a Weil divisor. 

We prove (2). Since $(X,B)$ is an exceptional surface pair, $X$ is klt. In particular, $X$ is $\Qq$-factorial. Since $\kappa_{\iota}(-(K_X+B))=2$, we may run a $-(K_X+B)$-MMP $g: X\rightarrow Y$ which terminates with a model $Y$ such that $-(K_Y+B_Y)$ is big and semi-ample, where $B_Y:=g_*B$. Then $(Y,B_Y)$ is an exceptional surface pair such that $B_Y\in\Phi_1$ and $\kappa_{\iota}(-(K_Y+B_Y))=2$. By (1), $I(K_Y+B_Y)$ is Cartier for some integer $I<I_0$. By the effective base-point-freeness theorem (\cite[Theorem 1.1, Remark 1.2]{Fuj09}, \cite[1.1 Theorem]{Kol93}), $-96I(K_Y+B_Y)$ is base-point-free. Thus $(Y,B_Y)$ has a monotonic $96I$-complement $(Y,B_Y^+)$. We let $$K_X+B^+:=g^*(K_Y+B_Y^+),$$ then $(X,B^+)$ is a monotonic $96I$-complement of $(X,B)=(X,B+B')$. This implies (2).
\end{proof}

\subsubsection{Conclusion}

\begin{thm}\label{thm: n21}
$N(2,1)<96\cdot (42)!\cdot 84^{128\cdot 42^5}$.
\end{thm}
\begin{proof}
It follows from Propositions \ref{prop: n21 nonexceptional case}, \ref{prop: n21 kod 0 case}, \ref{prop: n21 kod 1 case}, and \ref{prop: n21 kod 2 case}.
\end{proof}

The following result is stronger than Lemma \ref{lem: alternative defn epsilon12}(3) for surfaces.

\begin{lem}\label{lem: n21 with b' 41/42}
Let $(X,B+B')$ be an $\Rr$-complementary surface pair such that $B\in\Phi_1$ and $B'\in (\frac{41}{42},1)$. Then $(X,B+B')$ has a monotonic $N$-complement for some positive integer $N\leq N(2,1)$.
\end{lem}
\begin{proof}
By Proposition \ref{prop: klt standard coefficient surface 66}, $I(2,1)=66$, so $N(2,1)\geq 66$. By Proposition \ref{prop: n21 nonexceptional case}, we may assume that $(X,B+B')$ is exceptional. By Proposition \ref{prop: compare exc mld and glct} and Theorem \ref{thm: exceptional 1/42 lc global strong}, $B'=0$, and the lemma follows.
\end{proof}

\subsection{Lower bound of volumes}

Finally, we recall the following lower bound of surface log pairs of general type with standard coefficients.

\begin{thm}[{cf. \cite[Section 10]{AL19}}]\label{thm: surface standard volume lower bound}
Let $(X,B)$ be a projective lc surface pair, such that $B\in\Phi_1$ and $K_X+B$ is big. Then $$\vol(K_X+B)\geq\frac{1}{42\cdot 84^{128\cdot 42^5+168}}.$$
\end{thm}

\section{Proof of Theorem \ref{thm: explicit bdd nonexc fano threefold intro}}\label{sec: non-exceptional case}

The goal of this section is to prove Theorem \ref{thm: explicit bdd nonexc fano threefold intro}. More accurately, we will prove the following theorem in this section, which will immediately imply Theorem \ref{thm: explicit bdd nonexc fano threefold intro}:

\begin{thm}\label{thm: non-exceptional complement dim 3}
Let $X$ be a non-exceptional Fano type variety of dimension $3$. Then $X$ has an $N$-complement for some positive integer $N<192\cdot (42)!\cdot 84^{128\cdot 42^5}$.
\end{thm}

\subsection{Reduction of the theorem}

\begin{setup}\label{setup: non-exceptional}
$X,B,B',\pi,Z$ satisfy the following:
\begin{enumerate}
    \item $(X,B+B')$ is a $\Qq$-factorial projective lc pair of dimension $3$.
    \item $B$ is a non-zero Weil divisor and $B'\geq 0$,
    \item $X$ is of Fano type,
    \item $-(K_X+B+B')$ is semi-ample, and
    \item $\pi: X\rightarrow Z$ is a contraction such that $K_X+B+B'\sim_{\mathbb R,Z}0$.
\end{enumerate}
\end{setup}

\begin{lem}\label{lem: fano type has horiziontal boundary}
Conditions as in Set-up \ref{setup: non-exceptional}. Suppose that $\dim Z<\dim X$. Then some component of $\Supp B\cup\Supp B'$ is horizontal$/Z$.
\end{lem}
\begin{proof}
Let $F$ be a general fiber of $X\rightarrow Z$, $B_F:=B|_F$, and $B'_F:=B'|_F$. Since $K_X+B+B'\sim_{\mathbb R,Z}0$ and $X$ is of Fano type, $-K_F$ is big, and $K_F+B_F+B_F'\equiv 0$. So $B_F\not=0$ or $B_F'\not=0$. Thus a component of $\Supp B\cup\Supp B'$ is horizontal$/Z$.
\end{proof}

\begin{lem}\label{lem: fano type has horiziontal reduced boundary}
Conditions as in Set-up \ref{setup: non-exceptional}. Suppose that 
\begin{itemize}
    \item either $\dim Z=1$ and $B'\in (\frac{41}{42},1)$, or
    \item $\dim Z=2$ and $B'\in (\frac{5}{6},1)$.
\end{itemize}
Then no component of $\Supp B'$ is horizontal$/Z$, and some component of $\Supp B$ is horizontal$/Z$.
\end{lem}
\begin{proof}
Let $F$ be a general fiber of $X\rightarrow Z$, $B_F:=B|_F$, and $B'_F:=B'|_F$. Since $K_X+B+B'\sim_{\mathbb R,Z}0$ and $(X,B+B')$ is lc, $K_F+B_F+B_F'\equiv 0$ and $(F,B_F+B_F')$ is lc. By Proposition \ref{prop: glct 1-gap dim 1} and Theorem \ref{thm: exceptional 1/42 lc global strong}, $B'\in (1-\delta_{\glct}(\dim F,1),1)$, so $B'_F\in  (1-\delta_{\glct}(\dim F,1),1)$. By Lemma \ref{lem: alternative defn epsilon12}, $B'_F=0$. Thus no component of $\Supp B'$ is horizontal$/Z$. By Lemma \ref{lem: fano type has horiziontal boundary}, some component of $\Supp B$ is horizontal$/Z$.
\end{proof}

\subsection{The plt case}\label{subsec: the plt case}

In this subsection, we deal with the case when $(X,B+B')$ is plt. 

\subsubsection{The big case}

\begin{prop}\label{prop: non-exceptional kod 3 case}
Conditions as in Set-up \ref{setup: non-exceptional}. Assume that $(X,B+B')$ is plt, $$B'\in 
\left[1-\frac{1}{N(2,1)+1},1\right),$$ and $-(K_X+B+B')$ is big. Then $(X,B+B')$ has a monotonic $N$-complement for some positive integer $N\leq N(2,1)$. In particular, if $$B'\in \left[1-\frac{1}{96\cdot (42)!\cdot 84^{128\cdot 42^5}},1\right),$$ 
then $(X,B+B')$ has a monotonic $N$-complement for some positive integer $$N<96\cdot (42)!\cdot 84^{128\cdot 42^5}.$$
\end{prop}
\begin{proof}
It is a special case of \cite[Proposition 6.2]{PS01}. The in particular part follows from Theorem \ref{thm: n21}.
\end{proof}

\subsubsection{Kodaira dimension two case}

\begin{prop}\label{prop: non-exceptional kod 2 case}
Conditions as in Set-up \ref{setup: non-exceptional}. Assume that $B'\in (\frac{41}{42},1)$ and $\dim Z=2$. Then $(X,B+B')$ has a monotonic $N$-complement for some positive integer $$N\leq 2N(2,1)<192\cdot (42)!\cdot 84^{128\cdot 42^5}.$$
\end{prop}
\begin{proof}
By Lemma \ref{lem: fano type has horiziontal reduced boundary}, $B'$ is vertical$/Z$, and there exists a component $S$ of $B$ that is horizontal$/Z$. We let $S\xrightarrow{\pi_S} T\xrightarrow{\tau} Z$ be the Stein factorization of $f|_S$ and 
$$K_S+B_S:=(K_X+B+B')|_S.$$ 
Then $\pi_S$ is a birational morphism and $B_S\in\Phi_1\cup (\frac{41}{42},1)$. Since $-(K_X+B+B')$ is semi-ample, $(S,B_S)$ is $\Rr$-complementary. By Lemma \ref{lem: n21 with b' 41/42}, $(S,B_S)$ has a monotonic $N_1$-complement $(S,B_S^+)$ for some $N_1\leq N(2,1)$. We let $B_T:=(\pi_S)_*B_S$ and $B_T^+:=(\pi_S)_*B_S^+$. Then $(T,B_T^+)$ is a monotonic $N_1$-complement of $(T,B_T)$.
 
 Let $(Z,B_Z+M_Z)$ be a g-pair induced by the canonical bundle formula of $(X,B+B')\rightarrow Z$. By Theorem \ref{thm: cbf preserved under adjunction}, $(Z,B_Z+M_Z)$ is the g-pair induced by the canonical bundle formula of $f|_S: (S,B_S)\rightarrow Z$, 
$$K_T+B_T=\tau^*(K_Z+B_Z+M_Z),$$ 
 and $B_Z=\frac{1}{\deg\tau}\tau_*(R+B_T)$, where $R$ is the ramification divisor of $\tau$. In particular, the nef part $M'$ of $(Z,B_Z+M_Z)$ is $0$. We let $$B_Z^+:=B_Z+\frac{1}{\deg\tau}\tau_*(B_T^+-B_T).$$ By Lemma \ref{lem: plt stein degree 2}, $\deg\tau\leq 2$, so $2N_1(K_Z+B_Z^+)\sim 0$. Thus $(Z,B_Z^+)$ is a monotonic $2N_1$-complement of $(Z,B_Z)$. Since $\dim X-\dim Z=1$ and $S$ is horizontal$/Z$, $K_X+B+B'\sim 0$ over the generic point of $Z$. Thus $K_X+B+B'\sim \pi^*(K_Z+B_Z)$. Let $N:=2N_1$, then $$N\leq 2N_1\leq 2N(2,1),$$ and
 \begin{align*}
 N(K_X+B+B'+\pi^*(B_Z^+-B_Z))&=N(K_X+B+B')+N\pi^*(B_Z^+-B_Z)\\
 &\sim N\pi^*(K_Z+B_Z)+N\pi^*(B_Z^+-B_Z)\\
 &=N\pi^*(K_Z+B_Z^+)\sim 0.
 \end{align*}
 Thus $(X,B+B'+\pi^*(B_Z^+-B_Z))$ is a monotonic $N$-complement of $(X,B+B')$ and $N\leq 2N(2,1)$.

The inequality $2N(2,1)<192\cdot (42)!\cdot 84^{128\cdot 42^5}$ follows from Theorem \ref{thm: n21}.
\end{proof}

\subsubsection{Kodaira dimension one case}

\begin{prop}\label{prop: non-exceptional kod 1 case}
Conditions as in Set-up \ref{setup: non-exceptional}. Assume that $B'\in (\frac{41}{42},1)$ and $\dim Z=1$. Moreover, assume that there exists a component $S$ of $B$ such that $S$ is horizontal$/Z$ and $(X,B+B')$ is plt near $S$ over the generic point of $Z$. Then $(X,B+B')$ has a monotonic $N$-complement for some positive integer
$$N\leq \max_{1\leq q\leq 6}\lcm(q,2N_g(1,q))<72\cdot (42)!.$$
\end{prop}
\begin{proof}
By Lemma \ref{lem: fano type has horiziontal reduced boundary}, no component of $B'$ is horizontal$/Z$ and some component of $B$ is horizontal$/Z$. By Proposition \ref{prop: n21 nonexceptional case}, over the generic point of $Z$, $I(K_X+B+B')\sim 0$ for some integer $I\in\{1,2,3,4,6\}$. We let $S\xrightarrow{\pi_S} T\xrightarrow{\tau} Z$ be the Stein factorization of $f|_S$ and $$K_S+B_S:=(K_X+B+B')|_S.$$ 

By construction, $B_S\in\Phi_1\cup (\frac{41}{42},1)$ and $\dim T=1$. We let $(T,B_T+M_T)$ be a g-pair induced by the canonical bundle formula of $(S,B_S)\rightarrow T$ and $(Z,B_Z+M_Z)$ a g-pair induced by the canonical bundle formula of $(X,B)\rightarrow Z$. By Theorem \ref{thm: cbf preserved under adjunction}, $(Z,B_Z+M_Z)$ is a g-pair induced by the canonical bundle formula of $f|_S: (S,B_S)\rightarrow Z$. In particular, we have
$$K_T+B_T+M_T=\tau^*(K_Z+B_Z+M_Z).$$

If $(S,B_S)$ is non-exceptional, then by Proposition \ref{prop: n21 nonexceptional case}, $(S,B_S)$ has a monotonic $q$-complement $(S,B_S^+)$ for some $q\in\{1,2,3,4,6\}$. Since $K_S+B_S\sim_{\mathbb R,T}0$, $B_S^+-B_S$ is vertical$/T$. Let $(T,B_T^++M_T)$  be a g-pair induced by the canonical bundle formula of $(S,B_S^+)\rightarrow T$, then 
$$0\sim q(K_S+B_S^+)\sim q(K_T+B^+_T+M_T)$$ 
and $B_T^+\geq B_T$. Thus $(T,B^+_T+M_T)$ is a monotonic $q$-complement of $(T,B_T+M_T)$. Let $$B_Z^+:=B_Z+\frac{1}{\deg\tau}(B^+_T-B_T).$$  By Theorem \ref{thm: cbf preserved under adjunction}, $(Z,B_Z^++M_Z)$ is lc. By Lemma \ref{lem: plt stein degree 2}, $\deg\tau\leq 2$, so $(Z,B_Z^++M_Z)$ is a monotonic $2q$-complement of $(Z,B_Z+M_Z)$. Therefore, $$(X,B+B'+\pi^*(B_Z^+-B_Z))$$ is a monotonic $\lcm(I,2q)$-complement of $(X,B+B')$. Since $\lcm(I,2q)\leq 12$, we are done. Thus we may assume that  $(S,B_S)$ is exceptional. By Theorem \ref{thm: exceptional 1/42 lc global strong}, $B'=0$ and $B\in\Phi_1$. 

By Proposition \ref{prop: n21 nonexceptional case}, $(S/T\ni t,B_S)$ has a monotonic $q$-complement for some $q\in\{1,2,3,4,6\}$ for any point $t\in T$. By Proposition \ref{prop: cbf explicit boundary set}, $B_Z\in\Phi_q$, and we may choose $M_T$ so that $qM_T$ is a Weil divisor. Thus $(T,B_T+M_T)$ has a monotonic $N_1$-complement $(T,B_T^++M_T)$ for some positive integer $N_1\leq N_g(1,q)$. By Lemma \ref{lem: plt stein degree 2}, $\deg\tau\leq 2$. By Theorem \ref{thm: cbf preserved under adjunction}, $$\left(Z,B_Z^+:=B_T^++\frac{1}{\deg\tau}(B_T^+-B_T)+M_Z\right)$$ is a monotonic $2N_1$-complement of $(Z,B_Z+M_Z)$. Let $N:=\lcm(I,2N_1)$ and $$B^+:=B+B'+\pi^*(B_Z^+-B_Z),$$ then
\begin{align*}
    N(K_X+B^+)&=N(K_X+B+B'+\pi^*(B_Z^+-B_Z))\\
    &\sim N\pi^*(K_Z+B_Z+M_Z)+N\pi^*(B_Z^+-B_Z)\\
    &=N\pi^*(K_Z+B_Z^++M_Z)\sim 0.
\end{align*}
Therefore, $(X,B^+)$ is a monotonic $N$-complement of $(X,B+B')$. The proposition follows from Proposition \ref{prop: n complemntary bound dim 1}.
\end{proof}

\subsubsection{The numerically trivial case}

\begin{prop}\label{prop: non-exceptional kod 0 case picard 1}
Conditions as in Set-up \ref{setup: non-exceptional}. Assume that $B'\in (\frac{41}{42},1)$, $Z$ is a point, and $\rho(X)=1$. Moreover, suppose that $S$ is a component of $B$ such that $(X,S)$ is plt. Then $I(K_X+B+B')\sim 0$ for some positive integer $I\leq 66$. In particular, $(X,B+B')$ is a monotonic $I$-complement of itself.
\end{prop}
\begin{proof}
Let $$K_S+B_S:=(K_X+B+B')|_S.$$ Then $B_S\in\Phi_1\cup (\frac{41}{42},1)$, and if $B'\not=0$, then at least one coefficient of $B_S$ belongs to $(\frac{41}{42},1)$. Since $K_S+B_S\equiv 0$, by Theorem \ref{thm: exceptional 1/42 lc global strong} and Lemma \ref{lem: alternative defn epsilon12}(2), $B'=0$, and $B_S\in\Phi_1\cap ([0,\frac{41}{42}]\cup\{1\})$. By Proposition \ref{prop: klt standard coefficient surface 66}, $I(K_S+B_S)\sim 0$ for some positive integer $I\leq 66$. 

By Condition (1) of  Set-up \ref{setup: non-exceptional}, $X$ is $\Qq$-factorial, so $S$ is $\Qq$-Cartier. Since $I(K_S+B_S)\sim 0$, by Lemma \ref{lem: cartier index}, $I(K_X+B)$ is Cartier near the generic point of any prime divisor $P$ on $S$. By \cite[Lemma 2.42]{Bir19}, we have the short exact sequence
$$0\rightarrow\mathcal{O}_X(I(K_X+B)-S)\rightarrow\mathcal{O}_X(I(K_X+B))\rightarrow\mathcal{O}_S(I(K_S+B_S))\rightarrow 0$$
which induces the long exact sequence
$$\dots\rightarrow H^0(I(K_X+B))\rightarrow H^0(I(K_S+B_S))\rightarrow H^1(I(K_X+B)-S)\rightarrow\dots.$$
If $S\not=B$, then $$I(K_X+B)-S=K_X+(I(K_X+B)-S-K_X)$$ and $$I(K_X+B)-S-K_X\equiv B-S$$ is ample, so $H^1(I(K_X+B)-S)=0$ by Kawamata-Viehweg vanishing \cite[Theorem 1-2-7]{KMM87}. If $S=B$, then $$I(K_X+B)-S$$
is anti-ample, so $H^1(I(K_X+B)-S)=0$ by a variation of the Kawamata-Viehweg vanishing \cite[Theorem 2.70]{KM98}. Thus $$H^0(I(K_X+B))\rightarrow H^0(I(K_S+B_S))$$ is a surjection, so $|I(K_X+B)|\not=\emptyset$. Thus $I(K_X+B)\sim 0$. 
\end{proof}

\begin{prop}\label{prop: non-exceptional kod 0 case}
Conditions as in Set-up \ref{setup: non-exceptional}. Assume that  $(X,B+B')$ is plt, $B'\in (\frac{41}{42},1)$ and $Z$ is a point. Then $I(K_X+B+B')\sim 0$ for some positive integer $$I\leq\max\{2N(2,1),\max_{1\leq q\leq 6}\lcm(q,2N_g(1,q)),66\}<192\cdot (42)!\cdot 84^{128\cdot 42^5}.$$ In particular, $(X,B+B')$ is a monotonic $I$-complement of itself.
\end{prop}
\begin{proof}
We let $S$ be an irreducible component of $\lfloor B\rfloor$. We run a $(K_X+B+B'-S)$-MMP, which terminates with a Mori fiber space $\pi_Y: Y\rightarrow Z'$.  Let $B_Y,B_Y',S_Y$ be the images of $B,B',S$ on $Y$ respectively.
Since the induced birational map $X\dashrightarrow Y$ is a $(-S)$-MMP, $S_Y\not=0$.

\begin{claim}\label{claim: plt on y}
    $(Y,B_Y+B_Y')$ is plt.
\end{claim}
\begin{proof}
    Since $(X,B+B')$ is plt, there are only finitely many lc places of $(X,B+B')$. Since $(X,B+B')$ and $(Y,B_Y+B_Y')$ are crepant, $(Y,B_Y+B_Y')$ is lc, and there are only finitely many lc centers of $(Y,B_Y+B_Y')$. In particular, all lc centers of $(Y,B_Y+B_Y')$ are isolated. Since $S_Y\not=0$, $(Y,B_Y+B_Y')$ is plt near $S_Y$.

    Since the induced map $X\dashrightarrow Y$ is a $(K_X+B+B'-S)$-MMP, $(Y,B_Y+B_Y'-S_Y)$ is dlt, so $(Y,B_Y+B_Y'-S_Y)$ is plt. Thus $(Y,B_Y+B_Y')$ is plt outside $S_Y$.

    Therefore, $(Y,B_Y+B_Y')$ is plt.
\end{proof}
 By Claim \ref{claim: plt on y},   $(Y,B_Y+B_Y')$ is plt. By Propositions \ref{prop: non-exceptional kod 2 case}, \ref{prop: non-exceptional kod 1 case}, and  \ref{prop: non-exceptional kod 0 case picard 1}, $(Y,B_Y+B_Y')$ is a monotonic $I$-complement of itself for some positive integer 
 $$I\leq\max\{2N(2,1),\max_{1\leq q\leq 6}\lcm(q,2N_g(1,q)),66\}.$$ 
 Since $(X,B+B')$ is crepant to $(Y,B_Y+B_Y')$, the proposition follows.
\end{proof}

\subsection{The non-plt case}\label{subsec: non-plt}

In this section we deal with the case when $(X,B+B')$ is not plt. We shall provide two proofs with two different insights. The first proof is a very straightforward one which relies on \cite{FFMP22}. Since \cite{FFMP22} is an unpublished preprint, we will also explain the details of the result we cite in \cite{FFMP22}. The second proof is more lengthy and more algorithm-like, which is parallel to the different cases in subsection \ref{subsec: the plt case}.

\subsubsection{A quick proof via \cite{FFMP22}}

\begin{thm}\label{thm: dlt but not plt case}
Conditions as in Set-up \ref{setup: non-exceptional}. Assume that $(X,B+B')$ is dlt but not plt and $B'\in [\frac{6}{7},1)$. Then $(X,B+B')$ has a monotonic $N$-complement for some integer $N\in\{1,2,3,4,6\}$.
\end{thm}
\begin{proof}
Possibly replacing $X$ with a small $\Qq$-factorialization, we may assume that $X$ is $\Qq$-factorial. Since $(X,B+B')$ is not plt, $\lfloor B\rfloor$ is not normal. We let $B'':=\frac{6}{7}\Supp B'$. Then $B'\geq B''$. Since $X$ is $\Qq$-factorial, $(X,B+B'')$ is an $\Rr$-complementary lc pair. Since $\lfloor B\rfloor$ is not normal, $(X,B+B'')$ is not plt. Since $(X,B+B'')$ is $\Rr$-complementary and the coefficients of $B+B''$ belong to $\Phi_1$, by \cite[Theorems 4,5]{FFMP22}, $(X,B+B'')$ has a monotonic $N$-complement $(X,B^+)$ for some $N\in\{1,2,3,4,6\}$. Since $B''\in [\frac{6}{7},1)$, $$B^+\geq B+\Supp B''\geq B+B'.$$ Thus $(X,B^+)$ is a monotonic $N$-complement of $(X,B+B')$.
\end{proof}

\begin{rem} By taking a dlt modification, we may reduce any case in Set-up \ref{setup: non-exceptional} either to the plt case when the arguments in Subsection \ref{subsec: the plt case} can be applied, or to the case when Theorem \ref{thm: dlt but not plt case} can be applied. In the following, we will explain how the explicit bounds in \cite[Theorems 4,5]{FFMP22} are derived. 
    
We will only focus on the proof of \cite[Theorem 5]{FFMP22} as the pair case of \cite[Theorem 4]{FFMP22} is alternatively proven in \cite[Theorems B.1, B.2]{ABBDILW23} with a more detailed proof. Henceforth, we can freely reference \cite[Theorem 4]{FFMP22}.

Analyzing the proof of \cite[Proposition 7.4]{FFMP22} in depth, we aim to demonstrate that both $I=I(\Phi_1,1,1)$ from \cite[Theorem 6]{FFMP22} and $q=q(\Phi_1,1,1)$ from \cite[Theorem 8]{FFMP22} are ``essentially" contained within $\{1,2,3,4,6\}$. By ``essentially," we mean that under the conditions of \cite[Theorem 6]{FFMP22} (or \cite[Theorem 8]{FFMP22}), we can consistently select $I$ and $q$ from $\{1,2,3,4,6\}$ such that the relevant arguments remain valid, but the choices of $I$ and $q$ within this set might be varied under different scenarios. By making this small change, the outcome would state ``we can choose $N\in\{1,2,3,4,6\}$" rather than ``we can take $N=12=\lcm(1,2,3,4,6)$".

We first control $q$. The explicit value of $q$ can be traced to the proof of \cite[Proposition 6.3]{FFMP22}, and subsequently to the proofs of \cite[Theorem 6.2]{FFMP22}, and then \cite[Theorem 6.1]{FFMP22}. The first step of \cite[Theorem 6.1]{FFMP22} argues that $q$ can be chosen as the number $N(\lambda,0,1)$ provided in \cite[Theorem 5.3]{FFMP22}. It is important to note that we cannot extend the current line of reasoning since the proof of \cite[Theorem 5.3]{FFMP22} takes a least common multiple to get $N(\Lambda,c,p)$. However, as per \cite[Theorem 4]{FFMP22}, we can simply set $N(\lambda,0,1)=2$, which implies $q=2$.

Now we control $I$. This can be linked to the proof of \cite[Proposition 7.2]{FFMP22}, and subsequently to \cite[Theorem 4.9]{FFMP22}, followed by \cite[Theorem 4.8]{FFMP22}, and eventually to \cite[Lemma 4.6]{FFMP22}. Tracing this lineage, we only need to prove that $I(\Phi_1,1)$ in \cite[Lemma 4.6]{FFMP22} can be ``essentially" chosen as either $4$ or $6$. When $X$ is an elliptic curve, \cite[Remark 3.5]{Fuj01} completes the discussion. If we assume $X=\mathbb P^1$, then we let $(\tilde X,\tilde B)$ be the pair derived from the quotient $(X,B)/\Aut(X,B)$. Given this, $(\tilde X,\tilde B)$ is recognized as a klt log Calabi-Yau pair with standard coefficients, ensuring $N(K_{\tilde X}+\tilde B)\sim 0$ for some $N\in\{1,2,3,4,6\}$. The section in $N(K_{\tilde X}+\tilde B)$ induces an admissible section in $A(X,N(K_X+B))$, allowing $I(\Phi_1,1)$ in \cite[Lemma 4.6]{FFMP22} to be either $4$ or $6$. Thus, we can choose $I$ as $4$ or $6$. It is also important to note that $\Aut(X,B)$ may have order $\geq 12$; therefore, we have the control $\Aut(X,B)$-invariant sections rather than the order of $\Aut(X,B)$ in order the get the explicit bound $N\in\{1,2,3,4,6\}$.
\end{rem}

\begin{rem}
In the following, we shall provide an alternative proof for the non-plt case which is independent of \cite{FFMP22}. For the alternative proof, similar to what we have done before for the plt case, we need to treat four different cases for Set-up \ref{setup: non-exceptional}: the case when $-(K_X+B+B')$ is big, and the cases when $\dim Z=0,1,2$. We remark that the $\dim Z=2$ case is already done by Proposition \ref{prop: non-exceptional kod 2 case} so we will exclude this case in our arguments below.

We also remark that we cannot get the explicit bound $N\in\{1,2,3,4,6\}$ as in \cite[Theorems 4 and 5]{FFMP22}. Nevertheless, we will not get worse bounds: all bounds for the non-plt case are not larger than the bound
$$192\cdot (42)!\cdot 84^{128\cdot 42^5}$$
given in Proposition \ref{prop: non-exceptional kod 2 case}.
\end{rem}

\subsubsection{The big case}

\begin{prop}\label{prop: non-exceptional kod 3 case non-plt}
Conditions as in Set-up \ref{setup: non-exceptional}. Assume that $$B'\in 
\left(1-\frac{1}{N(2,1)+1},1\right),$$ and $-(K_X+B+B')$ is big. Then $(X,B+B')$ has a monotonic $N$-complement for some positive integer $N\leq N(2,1)$. In particular, if $$B'\in \left(1-\frac{1}{96\cdot (42)!\cdot 84^{128\cdot 42^5}},1\right),$$ 
then $(X,B+B')$ has a monotonic $N$-complement for some positive integer $$N<96\cdot (42)!\cdot 84^{128\cdot 42^5}.$$
\end{prop}
\begin{proof}
We let $$0<\epsilon\ll\frac{1}{N(2,1)+1}$$ be a real number such that 
$$(1-\epsilon)B'\geq\left(1-\frac{1}{N(2,1)+1}\right)\Supp B'.$$
Since $-(K_X+B+B')$ is big and semi-ample, possibly replacing $\pi: X\rightarrow Z$, we may assume that $\pi$ is defined by $-(K_X+B+B')$. We let $B_Z$ and $B_Z'$ be the images of $B$ and $B'$ on $Z$ respectively. Then $-(K_Z+B_Z+B_Z')$ is ample. By \cite[Lemma 2.7]{Bir21}, there exists a birational morphism $\phi: Y\rightarrow Z$ and a prime divisor $T$ on $Y$, such that
\begin{itemize}
        \item either $\phi$ is small or it contracts $T$ but no other divisors,
        \item $(Y,T)$ is plt,
        \item $-(K_Y+T)$ is ample$/Z$, and
        \item $a(T,Z,B_Z+B_Z')=0$.
\end{itemize}
We let $B_Y':=\phi^{-1}_*B_Z'$ and $B_Y:=\phi^{-1}_*B_Z$. Let $0<\delta\ll\epsilon$ be a real number, such that
$$-\delta(K_Y+T)-(1-\delta)\phi^*(K_Z+B_Z+B_Z')$$
is ample. We write
$$K_Y+T+B_Y'':=\delta(K_Y+T)+(1-\delta)\phi^*(K_Z+B_Z+B_Z').$$
By our conditions, $$B_Y''\in\left[1-\frac{1}{N(2,1)+1},1\right),$$ $(Y,T+B_Y'')$ is plt, and $-(K_Y+T+B_Y'')$ is ample. We let $\psi: W\rightarrow Y$ be a small $\Qq$-factorialization of $(Y,T+B_Y'')$, and let $B_W:=\psi^{-1}_*T$ and $B_W':=\psi^{-1}_*B_Y''$. By Proposition \ref{prop: non-exceptional kod 3 case}, $(W,B_W+B_W')$ has a monotonic $N$-complement $(W,B_W^+)$ for some positive integer $N\leq N(2,1)$. Since $$B_W'\in \left[1-\frac{1}{N(2,1)+1},1\right),$$ $(W,B_W^+)$ is a monotonic $N$-complement of $(W,B_W+\Supp B_W')$. Let $B_Z^+:=(\phi\circ\psi)_*B_W^+$, then $(Z,B_Z^+)$ is a monotonic $N$-complement of $(Z,B_Z+B_Z')$. Let 
$$K_X+B^+:=\pi^*(K_Z+B_Z^+),$$
then $(X,B^+)$ is a monotonic $N$-complement of $(X,B+B')$. The in particular part follows from Theorem \ref{thm: n21}.
\end{proof}

\subsubsection{Kodaira dimension one case}

To deal with the non-plt case when $\dim Z=1$, we use the $\gamma$-invariant introduced in \cite{LS23}.

\begin{nota}[{\cite[Notation 4.2]{LS23}}]\label{nota: gamma invariant}
Let $X\ni x$ (resp. $X$) be a klt surface germ (resp. a klt surface), $f: Y\rightarrow X$ the minimal resolution of $X\ni x$ (resp. $X$), and
$$K_Y+\sum_{i=1}^m b_iE_i=f^*K_X,$$
where $E_1,\dots,E_m$ are the prime $f$-exceptional divisors and $b_i=1-a(E_i,X,0)$ for each $i$. We put 
$$\gamma(X\ni x)\text{(resp. }\gamma(X)\text{)}:=m-\sum_{i=1}^mb_i(K_Y\cdot E_i).$$
It is clear that
$$\gamma(X)=\sum_{x\in X}\gamma(X\ni x).$$
\end{nota}

\begin{lem}\label{lem: 1/2 lc nonnegative gamma}
Let $X$ be a $\frac{1}{2}$-lc surface. Then $\gamma(X)\geq 0$.
\end{lem}
\begin{proof}
    We only need to show that $\gamma(X\ni x)\geq 0$ for any closed point $x\in X$. We may assume that $X$ is not smooth at $x$. Let  $f: Y\rightarrow X$ be the minimal resolution of $X\ni x$ and let
$$K_Y+\sum_{i=1}^m b_iE_i=f^*K_X.$$
If $E_i^2\leq -4$ for some $i$, then by the classification of surface singularities (cf. \cite[Corollary 2.19]{Ale93}), $m=1$ and $E_1^2=-4$, and we have $\gamma(X\ni x)=0$. If $E_i^2\geq -3$ for each $i$, then $K_Y\cdot E_i\leq 1$, so 
$$\gamma(X\ni x)=m-\sum_{i=1}^mb_i(K_Y\cdot E_i)=\sum_{i=1}^m(1-b_i(K_Y\cdot E_i))>0.$$
\end{proof}

\begin{lem}\label{lem: 1/2lc picard 9}
    Let $X$ be a $\frac{1}{2}$-lc del Pezzo surface. Then $\rho(X)\leq 9$.
\end{lem}
\begin{proof}
    Let  $f: Y\rightarrow X$ be the minimal resolution of $X$. Then $Y$ is rational, so $10=\rho(Y)+K_Y^2$ (cf. \cite[Lemma 3.19]{LS23}). Since $X$ is a $\frac{1}{2}$-lc del Pezzo surface, $K_X^2>0$. By \cite[Lemma 4.3]{LS23} and Lemma \ref{lem: 1/2 lc nonnegative gamma},
    $$0\leq\gamma(X)=\rho(Y/X)+K_Y^2-K_X^2=10-\rho(X)-K_X^2<10-\rho(X).$$
    Thus $\rho(X)\leq 9$.
\end{proof}

\begin{prop}\label{prop: non-exceptional kod 1 case nonplt}
Conditions as in Set-up \ref{setup: non-exceptional}. Assume that $B'\in (\frac{41}{42},1)$ and $\dim Z=1$. Then $(X,B+B')$ has a monotonic $N$-complement for some positive integer
$$N\leq \max\left\{\max_{1\leq q\leq 6}\lcm(q,2N_g(1,q)),2N(2,1)\right\}<192\cdot (42)!\cdot 84^{128\cdot 42^5}.$$
\end{prop}
\begin{proof}
    Possibly replacing $(X,B+B')$ with a dlt model, we may assume that $(X,B+B')$ is dlt. By Proposition \ref{prop: non-exceptional kod 1 case}, we may assume that $(X,B+B')$ is not plt over the generic point of $Z$. By Lemma \ref{lem: fano type has horiziontal reduced boundary}, no component of $B'$ is horizontal$/Z$ and some component of $B$ is horizontal$/Z$. By Proposition \ref{prop: n21 nonexceptional case}, over the generic point of $Z$, $2(K_X+B+B')\sim 0$.

    We let $F$ be a general fiber of $\pi$ and let $B_F:=B|_F$. Then $B'|_F=0$,
    $$K_F+B_F=(K_X+B+B')|_F,$$
    $(F,B_F)$ is dlt but not plt, and $2(K_F+B_F)\sim 0$. Therefore, $F$ is $\frac{1}{2}$-lc, so $X$ is $\frac{1}{2}$-lc over the generic point of $Z$.

    We run a $K_X$-MMP$/Z$ which terminates with a Mori fiber space $\pi_Y: Y\rightarrow Z'$. We let $B_Y$ and $B_Y'$ be the images of $B$ and $B'$ on $Y$ respectively. 
    
    If $\dim Z'>\dim Z$, then $\dim Z'=2$. By Proposition \ref{prop: non-exceptional kod 2 case}, $(Y,B_Y+B_Y')$ has a monotonic $N$-complement for some positive integer
    $$N\leq 2N(2,1)<192\cdot (42)!\cdot 84^{128\cdot 42^5}$$
    and we are done. Thus we may assume that $\dim Z'=\dim Z=1$, so $Z'=Z$ and $-K_Y$ is ample$/Z$. Since $X$ is $\frac{1}{2}$-lc over the generic point of $Z$, $Y$ is $\frac{1}{2}$-lc over the generic point of $Z$.

    Let $G$ be a general fiber of $\pi_Y$ and let $B_G:=B_Y|_G$. Then $G$ is a $\frac{1}{2}$-lc del Pezzo surface,
    $B_Y'|_G=0$, and
    $$K_G+B_G=(K_Y+B_Y+B_Y')|_G\sim_{\mathbb Q}0.$$
    By Lemma \ref{lem: 1/2lc picard 9}, $\rho(G)\leq 9$. Since $B_G$ is a reduced divisor, by \cite[Corollary 1.3]{BMSZ18}, there are at most $11$ components of $B_G$. Moreover, $B_G\not=0$, so there exists a component $S_Y$ of $B_Y$ that is horizontal$/Z$. 
    
    Let $S_Y^\nu$ be the normalization of $S_Y$ and let $V^\nu\subset S_Y^\nu$ be the smallest closed subset such that $S_Y^\nu\rightarrow S_Y$ is an isomorphism outside $V^\nu$. Then the general fibers of $S^\nu\rightarrow Z$ are of pure dimension $\dim S^\nu-\dim Z$. Thus no irreducible component of such fibers is contained in $V^\nu$ as the general fibers of $V^\nu\rightarrow Z$ have smaller dimension. This implies that a general fiber of $S_Y^\nu\rightarrow Z$ is birational to a general fiber of $S_Y\rightarrow Z$. Since $B_G$ has at most $11$ components, $\sdeg(S_Y^\nu/Z)\leq 11$.

    We let $h: W\rightarrow Y$ be a dlt modification of $(Y,B_Y+B_Y')$, $B_W':=h^{-1}_*B_Y'$, $S:=h^{-1}_*S_Y$, and 
    $$K_W+B_W+B_W':=h^*(K_Y+B_Y+B_Y').$$
    Then $\sdeg(S/Z)\leq 11$. 

    If $(W,B_W+B_W')$ is plt over the generic point of $Z$ near $S$, then by Proposition \ref{prop: non-exceptional kod 1 case}, $(W,B_W+B_W')$ has a monotonic $N$-complement for some positive integer
    $$N\leq \max_{1\leq q\leq 6}\lcm(q,2N_g(1,q))<72\cdot (42)!.$$
    Since $(X,B+B')$ and $(W,B_W+B_W')$ are crepant,  $(X,B+B')$ has a monotonic $N$-complement for some positive integer
    $$N\leq \max_{1\leq q\leq 6}\lcm(q,2N_g(1,q))<72\cdot (42)!$$
    and the proposition follows. Thus we may assume that $(W,B_W+B_W')$ is not plt over the generic point of $Z$ near $S$.

    Let $\pi_W:=\pi_Y\circ h$, $S\xrightarrow{\pi_S} T\xrightarrow{\tau} Z$ the Stein factorization of $\pi_W|_S$, and $K_S+B_S:=(K_W+B_W+B_W')|_S$. Let $(T,B_T+M_T)$ and be a g-pair induced by the canonical bundle formula of $\pi_S: (S,B_S)\rightarrow T$, and let $(Z,B_Z+M_Z)$ be a g-pair induced by the canonical bundle formula of $\pi_W: (W,B_W+B_W')\rightarrow Z$. By Theorem \ref{thm: cbf preserved under adjunction}, $(Z,B_Z+M_Z)$ is the g-pair induced by the canonical bundle formula of $\pi_W|_S: (S,B_S)\rightarrow Z$. In particular, we have
    $$K_T+B_T+M_T=\tau^*(K_Z+B_Z+M_Z).$$
    Moreover, since $(X,B+B')$ and $(W,B_W+B_W')$ are crepant, $(Z,B_Z+M_Z)$ is a g-pair induced by the canonical bundle formula of $(X,B+B')\rightarrow Z$.

    By Proposition \ref{prop: n21 nonexceptional case}, $(S,B_S)$ has a monotonic $q$-complement $(S,B_S^+)$ for some $q\in\{1,2,3,4,6\}$. Since $2(K_X+B+B')\sim 0$ over the generic point of $Z$ and $(X,B+B')$ and $(W,B_W+B_W')$ are crepant, $2(K_W+B_W+B_W')\sim 0$ over the generic point of $Z$. Thus
    $2(K_S+B_S^+)\sim 0$ over the generic point of $Z$. Since $K_S+B_S\sim_{\mathbb R,T}0$, $B_S^+-B_S$ is vertical$/T$. Let $(T,B_T^++M_T)$  be a g-pair induced by the canonical bundle formula of $(S,B_S^+)\rightarrow T$, then we have a choice of $M_T$ so that $$0\sim \lcm(2,q)(K_S+B_S^+)\sim \lcm(2,q)\pi_S^*(K_T+B^+_T+M_T)$$ and $B_T^+\geq B_T$. Thus $(T,B^+_T+M_T)$ is a monotonic $\lcm(2,q)$-complement of $(T,B_T+M_T)$. Let $$B_Z^+:=B_Z+\frac{1}{\deg\tau}(B^+_T-B_T).$$ 
    By Theorem \ref{thm: cbf preserved under adjunction}, $(Z,B_Z^++M_Z)$ is lc. Since $\sdeg(S/Z)\leq 11$, $\deg\tau\leq 11$. Thus $(Z,B_Z^++M_Z)$ is a monotonic $(11\lcm(2,q))$-complement of $(Z,B_Z+M_Z)$. Since $2(K_X+B+B')\sim 0$ over the generic point of $Z$, $$(X,B+B'+\pi^*(B_Z^+-B_Z))$$ is a monotonic $(11\lcm(2,q))$-complement of $(X,B+B')$. Since $$11\lcm(2,q)\leq 66=I(2,1)<2N(2,1),$$ the proposition follows.
    \end{proof}

\subsubsection{The numerical trivial case}

First we recall some basic definitions of B-representations, slc pairs, admissible sections, and pre-admissible sections.

\begin{defn}
Let $(X,B)$ and $(X',B')$ be two pairs. We say $f: (X,B)\dashrightarrow (X',B')$ is \emph{B-birational} if there exists a common resolution $\alpha: (Y,B_Y)\rightarrow (X,B)$ and $\beta: (Y,B_Y)\rightarrow (X',B')$, such that
$$K_Y+B_Y=\alpha^*(K_X+B)=\beta^*(K_{X'}+B')$$
and $\beta=f\circ\alpha$. We define
$$\Bir(X,B):=\{f|f: (X,B)\dashrightarrow (X,B) \text{ is B-birational}\}.$$
Let $n$ be a positive integer such that $n(K_X +B)$ is Cartier. Then we define 
$$\rho_n:\Bir(X,B)\rightarrow\Aut(H^0(X,n(K_X+B)))$$
to be the representation of the natural action of $\Bir(X, B)$ on $H^0(X,n(K_X+B))$ by pulling back sections.
\end{defn}

\begin{defn}
    A \emph{demi-normal scheme} is a reduced scheme of pure dimension that is $S_2$ and normal crossing in codimensional $1$. An \emph{slc pair} (resp. \emph{sdlt pair}) consists of a demi-normal scheme $X$ and an $\Rr$-divisor $B\geq 0$ on $X$ satisfying the following. Let $\nu: X'\rightarrow X$ be the normalization of $X$ and $K_{X'}+B'=\nu^*(K_X+B)$. Then
    \begin{enumerate}
        \item  $\Supp B$  does not contain any irreducible components of the conductor of $X$,
        \item $K_X+B$ is $\Rr$-Cartier, and
        \item $(X',B')$ is lc (resp. $(X',B')$ is dlt and $\nu$ is an isomorphism).
    \end{enumerate}
\end{defn}

\begin{defn}
We define \emph{preadmissible section} and \emph{admissible section} inductively in the following way. Let $d$ and $m$ be two positive integers. Let $(X,B)$ be a (possibly disconnected) projective sdlt pair of dimension $d$ such that $m(K_X+B)$ is Cartier. Let $\nu: X'\rightarrow X$ be the normalization of $X$, 
$$K_{X'}+B'+D':=\nu^*(K_X+B),$$
where $B'$ is the divisorial part of the inverse image of $B$ and $D'$ is the conductor on $X'$. We let $X_i$ be the irreducible components of $X$, $B_i:=B'|_{X_i}$, and $D_i:=D'|_{X_i}$ for each $i$. We let $D^\nu$ be the normalization of $D'$ and let
$$K_{D^\nu}+B_{D^\nu}:=(K_{X'}+B'+D')|_{D^\nu}.$$
\begin{enumerate}
    \item We say that $s\in H^0(X,m(K_X+B))$ is \emph{preadmissible} if $$s|_{D^\nu}\in H^0(D^\nu,m(K_{D^\nu}+B_{D^\nu}))$$
    is admissible. We remark that $s\in H^0(X,m(K_X+B))$ is preadmissible if and only if $s|_{D'}$ is admissible (\cite[Remark 5.2]{Gon10}, \cite[Remark 2.4.2.2]{Xu20}).
    \item We say that $s\in H^0(X,m(K_X+B))$ is \emph{admissible} if $s$ is preadmissible, and for any $i,j$ and B-birational map $g: (X_i,B_i+D_i)\dashrightarrow (X_j,B_j+D_j)$, $g^*(s|_{X_j})=s|_{X_i}$.
\end{enumerate}
The set of preadmissible (resp. admissible) sections in $H^0(X,m(K_X+B))$ is denoted by $\PA(X,m(K_X+B))$ (resp. $\A(X,m(K_X+B))$.
\end{defn}

\begin{lem}\label{lem: bir(x,b) order dim 1}
Let $n$ be a positive integer and $(X,B)$ a klt pair of dimension $1$ such that $n(K_X+B)\sim 0$. Then:
   \begin{enumerate}
       \item If $X$ is an elliptic curve, $B=0$, and $|\rho_{12}(\Aut(X))|=1$.
       \item If $X$ is not an elliptic curve, then $|\Bir(X,B)|\leq 6\cdot\binom{2n}{3}$. In particular, $$|\rho_n(\Bir(X,B))|\leq 6\cdot\binom{2n}{3}.$$
   \end{enumerate} 
\end{lem}
\begin{proof}
    This is essentially \cite[Proof of Proposition 3.3.1]{Xu20}.

    (1) It follows from \cite[12.2.9.1]{Kol+92}.

(2) There are at least $3$ components and at most $2n$ components of $\Supp B$, so 
$$|\Bir(X,B)|\leq 6\cdot\binom{2n}{3}.$$
\end{proof}

\begin{lem}\label{lem: admissible section dim 1}
Let $n$ be a positive integer and $(X,B)$ a projective klt pair of dimension $1$ such that $n(K_X+B)\sim 0$. Then there exists an admissible section of $A(X,N(K_X+B))$ for some positive integer
$$N\leq 6n\cdot\binom{2n}{3}.$$
\end{lem}
\begin{proof}
    If $X$ is an elliptic curve, then we let $m=12$ and $G:=\rho_{12}(\Aut(X))$. Otherwise, we let $m=n$ and $G=\rho_n(\Bir(X,B))$. Then there exists a section $s\in H^0(X,m(K_X+B))$. Since
    $$t:=\prod_{\sigma\in G}\sigma^*s\in H^0(X,m|G|(K_X+B))$$
    is $\Bir(X,B)$-invariant, $\A(X,m|G|(K_X+B))$ is non-trivial. Thus we may let $N=m|G|$. By Lemma \ref{lem: bir(x,b) order dim 1}, $N=12$ if $X$ is an elliptic curve, and
    $$N\leq 6n\cdot\binom{2n}{3}$$
    if $X$ is not an elliptic curve.
\end{proof}

\begin{lem}\label{lem: slc surface trivial complement}
    Let $(X,B)$ be a connected projective slc surface pair, such that $(X,B)$ is not lc, $B\in\Phi_1\cup [\frac{6}{7},1)$, and $K_X+B\equiv 0$. Then $B\in\Phi_1$ and $I(K_X+B)\sim 0$ for some even integer
    $$I\leq 6\cdot (7920)!.$$
\end{lem}
\begin{proof}
    We let $n: \tilde X\rightarrow X$ be the normalization of $X$ and let 
    $$K_{\tilde X}+\tilde B:=n^*(K_X+B).$$
    Then $(\tilde X,\tilde B)$ is an lc pair, any connected component of $(\tilde X,\tilde B)$ is not klt, $\tilde B\in\Phi_1\cup[\frac{6}{7},1)$, and $K_{\tilde X}+\tilde B\equiv 0$. By Proposition \ref{prop: n21 nonexceptional case}, there exists a positive integer $N_1\in\{4,6\}$ such that $N_1(K_{\tilde X}+\tilde B)\sim 0$. Thus $\tilde B\in \{1\}\cup [0,\frac{5}{6}]$, so $B\in \{1\}\cup [0,\frac{5}{6}]$. Therefore, $B\in\Phi_1$.
    
    By \cite[Proposition 3.2.12]{Xu20} (see also \cite[Propositions 4.5, 4.7]{Fuj00}) and Lemma \ref{lem: admissible section dim 1},   $$\PA(X,I(K_{\tilde X}+\tilde B))$$
    is non-trivial for some even integer 
    $$I\leq N_1\cdot\left(6N_1\cdot\binom{2N_1}{3}\right)!\leq 6\cdot (7920)!.$$
    By \cite[Proposition 2.4.5]{Xu20}  (see also \cite[Lemma 4.2]{Fuj00}), $I(K_X+B)\sim 0$.
\end{proof}

\begin{prop}\label{prop: non-exceptional kod 0 case picard 1 nonplt}
Conditions as in Set-up \ref{setup: non-exceptional}. Assume that $B'\in (\frac{41}{42},1)$, $Z$ is a point, and $\rho(X)=1$. Then $I(K_X+B+B')\sim 0$ for some positive even integer 
$$I\leq 6\cdot (7920)!.$$ In particular, $(X,B+B')$ is a monotonic $I$-complement of itself.
\end{prop}
\begin{proof}
Let $S$ be a component of $B$. If $(X,S)$ is plt, then the proposition follows from Proposition \ref{prop: non-exceptional kod 0 case picard 1}. Thus we may assume that $(X,S)$ is not plt. Let $$K_S+B_S:=(K_X+B+B')|_S.$$ Then $B_S\in\Phi_1\cup (\frac{41}{42},1)$, and if $B'\not=0$, then at least one coefficient of $B_S$ belongs to $(\frac{41}{42},1)$. By Lemma \ref{lem: slc surface trivial complement}, $B'=0$, $B_S\in\Phi_1$, and $I(K_S+B_S)\sim 0$ for some positive even integer $I\leq 6\cdot (7920)!$. 

By Condition (1) of  Set-up \ref{setup: non-exceptional}, $X$ is $\Qq$-factorial, so $S$ is $\Qq$-Cartier. Since $I(K_S+B_S)\sim 0$, by Lemma \ref{lem: cartier index}, $I(K_X+B+B')$ is Cartier near the generic point of any prime divisor on $S$. By \cite[Lemma 2.42]{Bir19}, we have the short exact sequence
$$0\rightarrow\mathcal{O}_X(I(K_X+B)-S)\rightarrow\mathcal{O}_X(I(K_X+B))\rightarrow\mathcal{O}_S(I(K_S+B_S))\rightarrow 0$$
which induces the long exact sequence
$$\dots\rightarrow H^0(I(K_X+B))\rightarrow H^0(I(K_S+B_S))\rightarrow H^1(I(K_X+B)-S)\rightarrow\dots.$$
If $S\not=B$, then $$I(K_X+B)-S=K_X+(I(K_X+B)-S-K_X)$$ and $$I(K_X+B)-S-K_X\equiv B-S$$ is ample, so $H^1(I(K_X+B)-S)=0$ by Kawamata-Viehweg vanishing \cite[Theorem 1-2-7]{KMM87}. If $S=B$, then $$I(K_X+B)-S$$
is anti-ample, so $H^1(I(K_X+B)-S)=0$ by a variation of the Kawamata-Viehweg vanishing \cite[Theorem 2.70]{KM98}. Thus $$H^0(I(K_X+B))\rightarrow H^0(I(K_S+B_S))$$ is a surjection, so $|I(K_X+B)|\not=\emptyset$. Thus $I(K_X+B)\sim 0$. 
\end{proof}

\begin{prop}\label{prop: non-exceptional kod 0 case nonplt}
Conditions as in Set-up \ref{setup: non-exceptional}. Assume that $B'\in (\frac{41}{42},1)$ and $Z$ is a point. Then $I(K_X+B+B')\sim 0$ for some positive integer $$I\leq\max\{\max_{1\leq q\leq 6}\lcm(q,2N_g(1,q)),2N(2,1),6\cdot (7920)!\}<192\cdot (42)!\cdot 84^{128\cdot 42^5}.$$  In particular, $(X,B+B')$ is a monotonic $I$-complement of itself.
\end{prop}
\begin{proof}
We let $S$ be an irreducible component of $\lfloor B\rfloor$. We run a $(K_X+B+B'-S)$-MMP, which terminates with a Mori fiber space $\pi_Y: Y\rightarrow Z'$.  Let $B_Y,B_Y',S_Y$ be the images of $B,B',S$ on $Y$ respectively.
Since the induced birational map $X\dashrightarrow Y$ is a $(-S)$-MMP, $S_Y\not=0$. By Propositions \ref{prop: non-exceptional kod 2 case}, \ref{prop: non-exceptional kod 1 case nonplt}, and \ref{prop: non-exceptional kod 0 case picard 1 nonplt}, $(Y,B_Y+B_Y')$ is a monotonic $I$-complement of itself for some positive integer 
$$I\leq\max\{\max_{1\leq q\leq 6}\lcm(q,2N_g(1,q)),2N(2,1),6\cdot (7920)!\}<192\cdot (42)!\cdot 84^{128\cdot 42^5}.$$ 
Since $(X,B+B')$ is crepant to $(Y,B_Y+B_Y')$, the proposition follows.
\end{proof}

\subsection{The general case}

\begin{prop}\label{prop: reduction of non-exceptional complement dim 3}
Conditions as in Set-up \ref{setup: non-exceptional}. Assume that $B'\in (1-\frac{1}{N(2,1)+1},1)$. Then $(X,B+B')$ has a monotonic $N$-complement for some positive integer $$N<192\cdot (42)!\cdot 84^{128\cdot 42^5}.$$
\end{prop}
\begin{proof}
It follows from Propositions \ref{prop: non-exceptional kod 2 case}, \ref{prop: non-exceptional kod 3 case non-plt}, \ref{prop: non-exceptional kod 1 case nonplt}, and \ref{prop: non-exceptional kod 0 case nonplt}.
\end{proof}

\begin{thm}\label{thm: non-exceptional complement dim 3 general}
Let $(X,B)$ be a non-exceptional Fano type lc pair of dimension $3$ such that $B\in (1-\frac{1}{N(2,1)+1},1)$. Then $(X,B)$ has a monotonic $N$-complement for some positive integer $$N<192\cdot (42)!\cdot 84^{128\cdot 42^5}.$$
\end{thm}
\begin{proof}
Let $(X,\tilde B)$ be an $\Rr$-complement of $(X,B)$ such that $(X,\tilde B)$ is not klt. Let $f: Y\rightarrow X$ be a dlt modification of $(X,\tilde B)$, $K_Y+\tilde B_Y:=f^*(K_X+\tilde B)$, $S_Y:=\lfloor B_Y\rfloor$, and $B_Y:=f^{-1}_*B-f^{-1}_*B\wedge S_Y$. Then $(Y,S_Y+B_Y)$ is $\Rr$-complementary and $Y$ is of Fano type. We may run a $-(K_Y+S_Y+B_Y)$-MMP $Y\dashrightarrow Y'$ which terminates with a model $Y'$ such that $-(K_{Y'}+S_{Y'}+B_{Y'})$ is semi-ample, where $S_{Y'},B_{Y'}$ are the images of $S_Y$ and $B_Y$ on $Y'$. By Lemma \ref{lem: antimmp does not contract lc place}, no component of $S_Y$ is contracted by $Y\dashrightarrow Y'$. By Proposition \ref{prop: reduction of non-exceptional complement dim 3}, $(Y',S_{Y'}+B_{Y'})$ has a monotonic $N$-complement $(Y',B_{Y'}^+)$ for some positive integer $N<192\cdot (42)!\cdot 84^{128\cdot 42^5}$. We let $p: W\rightarrow Y$ and $q: W\rightarrow Y'$ be a resolution of indeterminacy of $Y\dashrightarrow Y'$, and let $$K_X+B^+:=f_*p_*q^*(K_{Y'}+B_{Y'}^+).$$ Then $(X,B^+)$ is a monotonic $N$-complement of $(X,B)$.
\end{proof}

\begin{proof}[Proof of Theorem \ref{thm: non-exceptional complement dim 3}]
 It is a special case of Theorem \ref{thm: non-exceptional complement dim 3 general} by letting $B=0$.
\end{proof}

\begin{proof}[Proof of Theorem \ref{thm: explicit bdd nonexc fano threefold intro}]
 It is a special case of Theorem \ref{thm: non-exceptional complement dim 3} since any klt Fano variety is of Fano type.
\end{proof}

\begin{rem}
    It will be interesting to ask whether we can get an explicit boundedness of $N$-complements for exceptional Fano threefolds as well. It seems to us that solving this problem requires more complicated methods such as the explicit resolution and the explicit minimal model program for threefolds, so we will not discuss this question in this paper. 

    It is also compelling to consider if explicit boundedness of $N$-complements can be established for varieties with more general coefficients. In fact, even for klt $\mathbb R$-complementary pairs with arbitrary coefficients, it has been amazingly shown that the boundedness of complements still holds for Fano varieties \cite{Sho20} (even for generalized pairs \cite{CHHX23}), surfaces \cite{CHX23} (see also \cite{Zen23}), and threefolds \cite{CHX23}. These results are expected to play a vital role in future birational geometry research. In particular, from the standpoint of explicit birational geometry, the results in \cite{CHX23} led to question if explicit boundedness of $N$-complements for klt $\mathbb R$-complementary threefolds with arbitrary, or at least DCC boundary coefficients, can be achieved. However, this appears to be a more challenging question, even in dimensions 2 and 3, and alternative approaches are required.
    \end{rem}

\section{Proof of the main theorems}\label{sec: proof some main theorems}

In this section, we will first prove Theorem \ref{thm: explicit 1gap glct dim 3}, and then use it to prove Theorems \ref{thm: rct 1gap threefold},  \ref{thm: exc mld 3fold}, \ref{thm: cy mld 3fold}, \ref{thm: lower bound k+s volume threefold}, and Corollary \ref{cor: exc upper vol 3fold}.

\subsection{Proof of Theorem \ref{thm: explicit 1gap glct dim 3}}

\begin{prop}\label{prop: lct31=1/42}
$\delta_{\lct}(3,1)=\frac{1}{42}$.
\end{prop}
\begin{proof}
It follows from Theorem \ref{thm: exceptional 1/42 lc global strong} and Lemma \ref{lem: equality glct d and lct d+1}.
\end{proof}

\begin{prop}\label{prop: threefold 1-gap vol with upper bound}
Let $M>1$ be a real number and $I_0:=42\cdot 84^{128\cdot 42^5+168}$. Then $\delta:=\sqrt{\frac{1}{2I_0M}}$ satisfies the following.  Assume that $(X,B+B')$ is a projective lc pair of dimension $3$ such that $X$ is $\Qq$-factorial, $\rho(X)=1$, $B\in\Phi_1$, $B'\in [1-\delta,1)$, $\vol(-(K_X+B))\leq M$, and $K_X+B+B'\equiv 0$. Then $B'=0$.
\end{prop}
\begin{proof}
It is obvious that $\delta<\frac{1}{42}$. Suppose that $B'\not=0$. Let $T:=\Supp B'$ and $S$ a component of $B'$. By Proposition \ref{prop: lct31=1/42}, $(X,B+T)$ is lc. Since $\rho(X)=1$, $K_X+B+T$ is ample.

Let $S^\nu$ be the normalization of $S$ and
$$K_{S^\nu}+B_{S^\nu}:=(K_X+B+T)|_{S^\nu}.$$
Then $(S^\nu,B_{S^\nu})$ is lc and $K_{S^\nu}+B_{S^\nu}$ is ample. Moreover, $S^\nu\in\Phi_1$. By Theorem \ref{thm: surface standard volume lower bound}, 
\begin{align*}
   \frac{1}{I_0}&\leq\vol(K_{S^\nu}+B_{S^\nu})=(K_{S^\nu}+B_{S^\nu})^2=S\cdot (K_S+B+T)^2\\
   &\leq S\cdot(\delta T)^2\leq\delta^2T^3\leq\delta^2\left(\frac{1}{1-\delta}B'\right)^3=\frac{\delta^2}{(1-\delta)^3}(-K_X-B)^3\leq\frac{\delta^2M}{(1-\delta)^3}<2\delta^2M.
\end{align*}
Thus $\delta>\sqrt{\frac{1}{2I_0M}}$, a contradiction.
\end{proof}

\begin{prop}\label{prop: threefold 1-gap vol with no upper bound}
Let $I_2:=192\cdot (42)!\cdot 84^{128\cdot 42^5}$ and let $I_1>I_2$ be a real number. Assume that $(X,B)$ is a projective lc pair of dimension $3$ such that $X$ is $\Qq$-factorial klt, $\rho(X)=1$, $B\in [1-\frac{I_1-I_2}{I_2(I_1-1)},1]$, $\vol(-K_X)>27I_1^3$, and $K_X+B\equiv 0$. Then $B$ is reduced.
\end{prop}
\begin{proof}
Since $X$ is $\Qq$-factorial klt, $-K_X$ is big, and $\rho(X)=1$, $X$ is Fano. We have
$$\vol(B)=\vol(-K_X)>27I_1^3,$$
So
$$\vol\left(\frac{1}{I_1}B\right)>27.$$
By \cite[Lemma 3.2.2]{HMX14}, there exists an $\Rr$-divisor $0\leq D\sim_{\mathbb R}\frac{1}{I_1}B$, such that $(X,(1-\frac{1}{I_1})B+D)$ is not lc. Thus $(X,(1-\frac{1}{I_1})B)$ is non-exceptional. Since $B\in [1-\frac{I_1-I_2}{I_2(I_1-1)},1]$, the coefficients of $(1-\frac{1}{I_1})B$ belong to $[1-\frac{1}{I_2},1]$. 

Since $X$ is klt Fano, $X$ is of Fano type. By Theorems \ref{thm: n21} and \ref{thm: non-exceptional complement dim 3 general}, $(X,(1-\frac{1}{I_1})B)$ has a monotonic $N$-complement for some positive integer $N<I_2$. Thus $(X,\Supp B)$ has a monotonic $N$-complement. Since $K_X+B\equiv 0$, $B=\Supp B$.
\end{proof}

\begin{prop}\label{prop: glct 1-gap threefold}
Let  $I_0:=42\cdot 84^{128\cdot 42^5+168}$ and $a>I_0$ a real number. Then
$$\delta:=\min\left\{\frac{a-I_0}{I_0(a-1)},\sqrt{\frac{1}{54a^3I_0}}\right\}$$
satisfies the following.  Assume that $(X,B)$ is a projective lc pair of dimension $3$ such that $X$ is $\Qq$-factorial klt, $\rho(X)=1$, $K_X+B\equiv 0$, and $B\in [1-\delta,1]$. Then $B$ is reduced. 

In particular, if $B\in [1-\frac{1}{8I_0^2},1]$, then $B$ is reduced.
\end{prop}
\begin{proof}
We may assume that $B\not=0$. Then $K_X$ is not pseudo-effective. Since $\rho(X)=1$, and $X$ is $\Qq$-factorial klt, $X$ is klt Fano. In particular, $X$ is of Fano type.

We let $I_2:=192\cdot (42)!\cdot 84^{128\cdot 42^5}$. Then $a>I_0>I_2$. Since
$$f_a(x):=\frac{a-x}{x(a-1)}=-\frac{1}{a-1}+\frac{a}{a-1}\cdot\frac{1}{x}$$
is a strictly decreasing function of $x$ when $x>0$ for any fixed $a$, we have 
$B\in [1-\frac{a-I_2}{I_2(a-1)},1].$

Put $M:=27a^3$ and $I_1:=a$. If $\vol(-K_X)\leq M$, then since $\sqrt{\frac{1}{2I_0M}}=\sqrt{\frac{1}{54a^3I_0}}$, the proposition follows from Proposition \ref{prop: threefold 1-gap vol with upper bound}. If $\vol(-K_X)>M$, then since $B\in [1-\frac{I_1-I_2}{I_2(I_1-1)},1]$, the proposition follows from Proposition  \ref{prop: threefold 1-gap vol with no upper bound}. By taking $a=I_0+1$, the in particular part follows.
\end{proof}

\begin{thm}\label{thm: explicit 1gap glct dim 3 with T}
Let $I_0:=2\cdot 84^{256\cdot 42^5+338}$. Let $t\in (0,1)$ be a real number, and let $(X,B:=T+tS)$ be a projective lc pair of dimension $3$ such that $K_X+T+tS\equiv 0$, $T$ is an effective Weil divisor, and $S$ is a non-zero effective Weil divisor. Then $t<1-\frac{1}{I_0}$. In particular, there exists a positive real number $\epsilon$ such that $t\leq 1-\frac{1}{I_0}-\epsilon$.
\end{thm}
\begin{proof}
We may assume that $t>\frac{41}{42}$. Possibly replacing $(X,B)$ with a dlt model, we may assume that $(X,B)$ is $\Qq$-factorial dlt. We run a $(K_X+T)$-MMP which terminates with a Mori fiber space $X'\rightarrow Z$. Let $T',S'$ be the image of $T,S$ on $X'$, $F$ a general fiber of $X'\rightarrow Z$, $T_F:=T'|_F$ and $S_F:=S'|_F$. Then $(X',T'+tS')$ is lc, $S'\not=0$, and $S'$ is horizontal over $Z$.

If $\dim Z\geq 1$, then $\dim F\leq 2$. Since $(F,T_F+tS_F)$ is lc and $K_F+T_F+tS_F\equiv 0$, by Theorem \ref{thm: exceptional 1/42 lc global strong}, $t\leq\frac{41}{42}$, a contradiction. Thus $\dim Z=0$. Therefore, $X'$ is $\Qq$-factorial, $\rho(X')=1$, $K_{X'}+T'+tS'\equiv 0$, $S'\not=0$, and $T'+tS'\in [t,1]$. By Proposition \ref{prop: glct 1-gap threefold}, $t<1-\frac{1}{I_0}$.

The in particular part of the theorem follows from the global ACC \cite[Theorem 1.5]{HMX14}.
\end{proof}

\begin{proof}[Proof of Theorem \ref{thm: explicit 1gap glct dim 3}]
    It is a special case of Theorem \ref{thm: explicit 1gap glct dim 3 with T} by letting $T=0$.
\end{proof}

\subsection{Proof of Theorem \ref{thm: rct 1gap threefold}}

\begin{thm}\label{thm: r complementary 1-gap dim 3}
Let $I_0:=2\cdot 84^{256\cdot 42^5+338}$. Let $(X,T)$ be an $\Rr$-complementary projective threefold pair such that $T$ is an effective Weil divisor, and $S$ is a non-zero effective Weil divisor on $X$. Then $\Rct(X,T;S)=1$ or $\Rct(X,T;S)<1-\frac{1}{I_0}$. 
\end{thm}
\begin{proof}
The proof is very similar to the proof of Lemma \ref{lem: rct gap is glct gap}. For the reader's convenience, we provide a full proof here.

Suppose that the theorem does not hold. Then there exists an $\Rr$-complementary threefold pair $(X,T)$ and a non-zero effective Weil divisor $S$ on $X$, such that $$t:=\Rct(X,T;S)\in [1-\frac{1}{I_0},1).$$ Since $t>\frac{1}{2}$, $S$ is reduced. 

We let $(X,T+tS+G)$ be an $\Rr$-complement of $(X,T+tS)$ and let $f: Y\rightarrow X$ be a $\Qq$-factorial dlt modification of $(X,T+tS+G)$. We let $R:=\lfloor T+tS+G\rfloor$, $G_Y$ the strict transform of $G-G\wedge R$ on $Y$, $S_Y$ the strict transform of $S-S\wedge R$ on $Y$, and
$$K_Y+T_Y+tS_Y+G_Y:=f^*(K_X+tS+G).$$
Then $T_Y$ is a reduced divisor, and $\lfloor T_Y+tS_Y+G_Y\rfloor=T_Y$. In particular, there exists $\tau\in (0,t)$ such that $(Y,T_Y+(1+\tau)(tS_Y+G_Y))$ is dlt, and
$$-\tau(K_Y+T_Y+S_Y)\sim_{\mathbb R}K_Y+T_Y+(1+\tau)(tS_Y+G_Y)-\tau S_Y.$$
 Since $\tau\in (0,t)$, $(1+\tau)t-\tau>0$, so $$(Y,B_Y:=T_Y+(1+\tau)(tS_Y+G_Y)-\tau S_Y)$$ is dlt. 
 
 There are two cases.

\medskip

\noindent\textbf{Case 1}. $-(K_Y+T_Y+S_Y)$ is not pseudo-effective. Since $$-\tau(K_Y+T_Y+S_Y)\sim_{\mathbb R}K_Y+B_Y,$$
we may run a $(K_Y+B_Y)$-MMP $h: Y\dashrightarrow Y'$ with scaling of an ample divisor, which terminates with a Mori fiber space $\pi: Y'\rightarrow Z$. We let $T_{Y'},S_{Y'}$ be the images of $T_Y,S_Y$ on $Y'$ respectively, then $-(K_{Y'}+T_{Y'}+S_{Y'})$ is anti-ample$/Z$. Since $(Y,T_Y+tS_Y)$ is $\Rr$-complementary, $(Y',T_{Y'}+tS_{Y'})$ is $\Rr$-complementary. In particular,  $(Y',T_{Y'}+tS_{Y'})$ is lc. Since $t\in [1-\frac{1}{I_0},1)$, by Proposition \ref{prop: lct31=1/42}, $(Y',T_{Y'}+S_{Y'})$ is lc. Moreover, $-(K_{Y'}+T_{Y'}+tS_{Y'})$ is pseudo-effective, hence $-(K_{Y'}+T_{Y'}+tS_{Y'})$ is nef$/Z$. Thus there exists $t'\in [t,1)$ such that $(Y',T_{Y'}+t'S_{Y'})$ is lc and $K_{Y'}+T_{Y'}+t'S_{Y'}\equiv_Z0$. Let $F$ be a general fiber of $Y'\rightarrow Z$, $T_{F}:=T_{Y'}|_F$, and $S_F:=S_{Y'}|_F$. Then $K_{F}+T_{F}+t'S_{F}\equiv0$. By Theorem \ref{thm: explicit 1gap glct dim 3 with T}, $t'<1-\frac{1}{I_0}$, a contradiction.

\medskip

\noindent\textbf{Case 2}. $-(K_Y+T_Y+S_Y)$ is  pseudo-effective. Since $$-\delta(K_Y+T_Y+S_Y)\sim_{\mathbb R}K_Y+B_Y$$ and $(Y,B_Y)$ is dlt, by the existence of good minimal models in dimension $3$, we may run a $(K_Y+B_Y)$-MMP which terminates with a good minimal model $(Y',B_{Y'})$ of $(Y,B_Y)$.

We let $T_{Y'},S_{Y'}$ be the images of $T_Y,S_Y$ on $Y'$ respectively, then $-(K_Y+T_{Y'}+S_{Y'})$ is semi-ample. Since $(Y,T_Y+tS_Y)$ is $\Rr$-complementary, $(Y',T_{Y'}+tS_{Y'})$ is $\Rr$-complementary. In particular,  $(Y',T_{Y'}+tS_{Y'})$ is lc. Since $t>\frac{41}{42}$, by Proposition \ref{prop: lct31=1/42},  $(Y',T_{Y'}+S_{Y'})$ is lc. Thus we may pick $D\in |-(K_{Y'}+T_{Y'}+S_{Y'})|_{\mathbb R}$ such that $$(Y',B_{Y'}^+:=T_{Y'}+S_{Y'}+D)$$ is lc. Let $p: W\dashrightarrow Y$, $q: W\rightarrow Y'$ be a common resolution of the induced birational map $Y\dashrightarrow Y'$, and let $$K_X+B^+:=f_*p_*q^*(K_{Y'}+B_{Y'}^+).$$ Then $(X,B^+)$ is an $\Rr$-complement of $(X,T+S)$. Thus $\Rct(X,T;S)=1$, a contradiction.
\end{proof}

\begin{proof}[Proof of Theorem \ref{thm: rct 1gap threefold}]
    It is a special case of Theorem \ref{thm: r complementary 1-gap dim 3} by letting $T=0$.
\end{proof}

\subsection{Proof of Theorem \ref{thm: exc mld 3fold}, Theorem \ref{thm: cy mld 3fold}, and Corollary \ref{cor: exc upper vol 3fold}}

\begin{thm}\label{thm: exc mld dim 3}
Let $I_0:=2\cdot 84^{256\cdot 42^5+338}$. Let $(X,B)$ be an exceptional threefold such that $B\in [1-\frac{1}{I_0},1]$. Then for any $0\leq G\sim_{\mathbb R}-(K_X+B)$, $(X,B+G)$ is $\frac{1}{I_0}$-klt. In particular, $B=0$.
\end{thm}
\begin{proof}
let $a:=\tmld(X,B+G)$. Suppose that $a\leq\frac{1}{I_0}$. Possibly taking a small dlt modification, we may assume that $X$ is $\Qq$-factorial. We let $E$ be a prime divisor over $X$ such that $a(E,X,B+G)=a$. Then either $E$ is on $X$ and we let $f: Y\rightarrow X$ be the identity morphism, or there exists a divisorial contraction $f: Y\rightarrow X$ which extracts $E$, such that $Y$ is $\Qq$-factorial. We let $B_Y:=f^{-1}_*B-f^{-1}_*B\wedge E$. Then $(Y,B_Y+(1-a)E)$ is $\Rr$-complementary. Since $B\in [1-\frac{1}{I_0},1]$ and  $1-a\geq 1-\frac{1}{I_0}$, 
$$B_Y+(1-a)E\in \left[1-\frac{1}{I_0},1\right].$$ 
We let $S:=\Supp(B_Y+E)$, then $S$ is a non-zero reduced divisor and $B_Y+(1-a)E\geq (1-\frac{1}{I_0})S$. Since $(Y,B_Y+(1-a)E)$ is $\Rr$-complementary, $(Y, (1-\frac{1}{I_0})S)$ is $\Rr$-complementary. Thus 
$$\Rct(Y,0;S)\geq 1-\frac{1}{I_0}.$$
By Theorem \ref{thm: r complementary 1-gap dim 3}, $\Rct(Y,0;S)=1$. Thus $(Y,B_Y+E)$ has an $\Rr$-complement $(Y,B_Y^+)$. Then $(X,B^+:=f_*B_Y^+)$ is an $\Rr$-complement of $(X,B)$, such that $a(E,X,B^+)=0$. This is not possible as $(X,B)$ is exceptional.
\end{proof}

\begin{proof}[Proof of Theorem \ref{thm: exc mld 3fold}]
It immediately follows from Theorem \ref{thm: exc mld dim 3}.
\end{proof}

\begin{proof}[Proof of Theorem \ref{thm: cy mld 3fold}]
By abundance (cf. \cite[V.4.6 Theorem]{Nak04}), any klt Calabi-Yau threefold is exceptional. The theorem follows from Theorem \ref{thm: exc mld dim 3}.
\end{proof}

\begin{proof}[Proof of Corollary \ref{cor: exc upper vol 3fold}]
We may assume that $\vol(-K_X)>0$, so $-K_X$ is big and $X$ is of Fano type. Possibly replacing $X$ with a small $\Qq$-factorialization, we may assume that $X$ is $\Qq$-factorial. We run a $(-K_X)$-MMP which terminates with a model $X'$ such that $\vol(-K_X)=\vol(-K_{X'})$ and $-K_{X'}$ is ample. Thus $X'$ is an exceptional Fano threefold. The corollary follows from Theorem \ref{thm: exc mld dim 3} and \cite[Theorem 1.1]{JZ23b}.
\end{proof}

\subsection{Proof of Theorem \ref{thm: lower bound k+s volume threefold}}

\begin{defn}\label{defn: pet}
Let $(X,B)$ be a projective lc pair and $D\geq 0$ an $\Rr$-Cartier $\Rr$-divisor on $X$. We define
$$\pet(X,B;D):=\inf\{+\infty,t\geq 0\mid (X,B+tD)\text{ is lc},  K_X+B+tD\text{ is pseudo-effective}\}\in\mathbb R^+\cup\{+\infty\}.$$
\end{defn}

\begin{lem}\label{lem: threefold gap pet qfact}
    Let $I_0:=2\cdot 84^{256\cdot 42^5+338}$. Let $(X,T)$ be a $\Qq$-factorial projective dlt threefold pair such that $T$ is an effective Weil divisor. Let $S$ be a non-zero effective Weil divisor on $X$. Then $\pet(X,T;S)\not\in [1-\frac{1}{I_0},1)$.
\end{lem}
\begin{proof}
    We let $t:=\pet(X,T;S)$. Suppose that $t\in [1-\frac{1}{I_0},1)$. We run a $(K_X+T+(t-\epsilon)S)$-MMP for some $0<\epsilon\ll 1$, which terminates with a Mori fiber space $f: X'\rightarrow Z$. Let $T'$ and $S'$ be the images of $T$ and $S$ on $X'$ respectively. Since $(X',T'+(t-\epsilon)S')$ is lc and $t-\epsilon>\frac{41}{42}$, by Proposition \ref{prop: lct31=1/42},  $(X',T'+S')$ is lc. Since $K_X+T+tS$ is pseudo-effective, $K_{X'}+T'+tS'$ is pseudo-effective, so $K_{X'}+T'+tS'$ is nef$/Z$. Since $K_{X'}+T'+(t-\epsilon)S'$ is anti-ample$/Z$, $S'$ is ample$/Z$, and there exists $t'\in (t-\epsilon,t]$ such that $K_{X'}+T'+t'S'\equiv_Z0$. In particular, $S'$ is horizontal$/Z$. Let $F$ be a general fiber of $X'\rightarrow Z$, $T_F:=T'|_F$, and $S_F:=S'|_F$. Then $(F,T_F+t'S_F)$ is lc, $S_F\not=0$, and $K_F+T_F+t'S_F\equiv 0$. Since 
    $$t-\epsilon<t'\leq t<1$$ 
    and 
    $$t\geq 1-\frac{1}{I_0},$$
    we get a contradiction to Theorem \ref{thm: explicit 1gap glct dim 3 with T}. 
\end{proof}

\begin{lem}\label{lem: pet gap threefold nonqfact}
    Let $I_0:=2\cdot 84^{256\cdot 42^5+338}$. Let $(X,T+S)$ be an lc projective threefold pair such that $K_X+T+S$ is big, and $T$ and $S$ are effective Weil divisors. Then $K_X+T+(1-\frac{1}{I_0})S$ is big.
    
    We remark that $K_X+T+(1-\frac{1}{I_0})S$ may not be $\Qq$-Cartier.
\end{lem}
\begin{proof}
    Let $f: Y\rightarrow X$ be a dlt modification of $(X,T+S)$, $S_Y:=f^{-1}_*S$, and $$K_Y+T_Y+S_Y:=f^*(K_X+T+S).$$ Then $(Y,T_Y+S_Y)$ is $\Qq$-factorial dlt. Let $t:=\pet(Y,T_Y;S_Y)$. By Lemma \ref{lem: threefold gap pet qfact}, $t=1$ or $t<1-\frac{1}{I_0}$. Since $K_Y+T_Y+S_Y$ is big, $t\not=1$, so $t<1-\frac{1}{I_0}$. Thus $$K_X+T+tS=f_*(K_Y+T_Y+tS_Y)$$ is pseudo-effective, so 
    $$K_X+T+\left(1-\frac{1}{I_0}\right)S=\frac{1-\frac{1}{I_0}-t}{1-t}(K_X+T+S)+\frac{1}{I_0(1-t)}(K_X+T+tS)$$
    is big.
\end{proof}

\begin{proof}[Proof of Theorem \ref{thm: lower bound k+s volume threefold}]
Let $I_0:=2\cdot 84^{256\cdot 42^5+338}$ and $I_1:=42\cdot 84^{128\cdot 42^5+168}$. Let $T$ be the normalization of an irreducible component of $S$ and let $K_T+B_T:=(K_X+S)|_T$. Then $(T,B_T)$ is lc, $B_T\in\Phi_1$, and $K_T+B_T$ is ample. By Theorem \ref{thm: surface standard volume lower bound}, $$\vol(K_T+B_T)\geq (K_T+B_T)^2\geq\frac{1}{I_1}.$$ 
Thus 
\begin{align*}
    \vol(K_X+S)&=(K_X+S)^3\\
    &=(K_X+S)^2\left(K_X+(S-T)+\left(1-\frac{1}{I_0}\right)T\right)+\frac{1}{I_0}(K_X+S)^2\cdot T\\ 
    &>\frac{1}{I_0}(K_X+S)^2\cdot T && (\text{Lemma \ref{lem: pet gap threefold nonqfact}})\\
    &=\frac{1}{I_0}(K_T+B_T)^2\geq\frac{1}{I_0}\cdot\frac{1}{I_1}=\frac{1}{84^{384\cdot 42^5+507}}.
\end{align*}
The theorem follows.
\end{proof}

\section{On explicit geometry of surfaces with more general coefficients}\label{sec: Non-standard coefficient case}

In this section, we study the explicit geometry of surface pairs with coefficients contained in the set $\Phi_p$, where $p$ is any positive integer.

\begin{rem}
We briefly explain why we are focusing on the sets of the form $\Phi_p$ rather than arbitrary DCC sets.

In \cite{AM04}, the authors constructed a function $\beta$ to estimate the explicit bounds of algebraic invariants for surfaces for pairs with arbitrary DCC coefficients. More precisely,  \cite{AM04} has found an ``explicitly computable" function 
$$\beta:\{\text{DCC set of real numbers in }[0,1]\}\rightarrow (0,1),$$
such that for any DCC set $\Ii$ and an lc log Calabi-Yau surface pair $(X,B+B')$ such that $K_X+B+B'\equiv 0$, $B\in\Ii$, and $B'\in (1-\beta(\Ii),1)$, we have $B'=0$. 

However, the concept of ``explicitly computable" itself, is actually not so explicit. The reason is that, when constructing $\beta$, we essentially need to consider the function
$$f: (\Ii,q)\rightarrow\min\left\{\gamma>q\mid \gamma=\sum\gamma_i,\gamma_i\in\Ii\right\}-q.$$
This function is, in fact, not so ``explicitly computable". For example, even if we fix $\Ii$ to be the finite set with two elements $$\left\{\frac{\sqrt{2}}{2},\frac{\sqrt{3}}{3}\right\},$$ 
or the hyperstandard set $$\left\{1-\frac{\sqrt{2}}{n}\Biggm| n\in\mathbb N^+\right\}=\Phi(\{\sqrt{2}\}),$$ 
and $q$ is the only variable for $f(\Ii,q)$, it is still very difficult to get an explicit characterization of the function $$f(\Ii,-):\mathbb N^+\rightarrow\mathbb (0,1).$$
Therefore, we probably need to add some conditions to the set $\Ii$. 

The reason why we consider $\Phi_p$ is now becoming very clear. On the one hand, we want to consider the set of all finite sets of rational numbers and we want the set to be preserved under adjunction, so we have to consider the sets of the form $\Phi_p$. On the other hand, instead of considering all DCC set of real numbers, which is obviously uncountable, we want to consider a countable set of DCC sets so that they can be parameterized by integers: once we have the parametrization
$$g: \mathbb N^+\rightarrow \{\text{The DCC set we want to care about}\},$$
$g\circ\beta$ becomes a function from $\mathbb N^+$ to $\mathbb R$, which will be much easier for us to visualize. Since $p\rightarrow\Phi_p$ is a natural choice of the parametrization $g$, at the end of the day, it becomes natural to consider $\Phi_p$.
\end{rem}

\subsection{The beta function in Alexeev-Mori}

In this section, we study the beta function of Alexeev-Mori when estimating the lower bound of global lc thresholds, and provide an explicit version of this function for the sets of the form $\Phi_p$.

\begin{defn}
Let $p,n$ be two positive integers and $r$ a real number. We let $(p\uparrow)^nr$ be the number defined in the following way: $(p\uparrow)^1r=p^r$, $(p\uparrow)^{n}r=p^{(p\uparrow)^{n-1}r}$.
\end{defn}

\begin{defn}
Let $p,q$ be two positive integers and $\epsilon\in (0,2]$ a real number. We define
\begin{itemize}
    \item $$\epsilon_1(p,q):=\min\left\{\gamma\mid\gamma>q, \gamma=\sum\gamma_i,\gamma_i\in\Phi_p\right\}-q,$$
    \item $$\epsilon_2(p,q):=\frac{\epsilon_1(p,q)}{q+\epsilon_1(p,q)},$$
    \item $$M_\epsilon:=\left\lfloor\frac{2}{\epsilon}\right\rfloor^{\left\lfloor\frac{128}{\epsilon^5}\right\rfloor}\times (\left\lfloor\frac{2}{\epsilon}\right\rfloor+2),$$
    \item $$\alpha(p,\epsilon):=\epsilon_2(p,M_\epsilon),$$
    \item $$\beta(p):=\alpha(p,\alpha(p,\alpha(p,\alpha(p,2)))),$$
    \item $$l(p):=\left\lceil\frac{p}{\beta(p)}\right\rceil,$$ and
    \item $$\upsilon(p):=\frac{1}{l(p)\cdot (2l(p))^{128l(p)^5+4l(p)}}.$$
\end{itemize}
\end{defn}

\begin{lem}\label{lem: compare of mf1 and mf2 of am04}
Let $q>q'$ be two positive integers. Then $\epsilon_1(p,q)\leq\epsilon_1(p,q')$ and $\epsilon_2(p,q)<\epsilon_2(p,q')$.
\end{lem}
\begin{proof}
By induction we may assume that $q=q'+1$. Since $1\in\Phi_p$,
\begin{align*}
    \epsilon_1(p,q)&=\min\left\{\gamma\mid\gamma>q, \gamma=\sum\gamma_i,\gamma_i\in\Phi_p\right\}-q\\
    &\leq \min\left\{1+\gamma\mid 1+\gamma>q, \gamma=\sum\gamma_i,\gamma_i\in\Phi_p\right\}-q\\
    &=\min\left\{\gamma\mid\gamma>q-1,\gamma=\sum\gamma_i,\gamma_i\in\Phi_p\right\}-(q-1)=\epsilon_1(p,q').
\end{align*}
So
$$\epsilon_2(p,q)=\frac{\epsilon_1(p,q)}{q+\epsilon_1(p,q)}\leq\frac{\epsilon_1(p,q')}{q+\epsilon_1(p,q')}<\frac{\epsilon_1(p,q')}{q'+\epsilon_1(p,q')}=\epsilon_2(p,q').$$
\end{proof}

\begin{lem}\label{lem: compare of alpha of am04}
Let $p$ be a positive integer and $2\geq\epsilon\geq\epsilon'>0$ two real numbers. Then $\alpha(p,\epsilon)\geq\alpha(p,\epsilon')$.
\end{lem}
\begin{proof}
We have $M_{\epsilon}\leq M_{\epsilon'}$, and the lemma follows from Lemma \ref{lem: compare of mf1 and mf2 of am04}.
\end{proof}

\begin{defn}
Let $n$ be a positive integer. Then $n$-th Sylvester number, denoted by $S_n$ is defined by $S_1=2$ and $S_n^2=S_{n-1}^2-S_{n-1}+1$ for any $n\geq 2$. We have $S_n=\prod_{i=1}^{n-1}S_i+1$ for any $n\geq 2$, and $S_n\leq 2^{2^n}$ for any $n$.
\end{defn}

\begin{thm}[{\cite[Theorem 1]{Cur22}}]\label{thm: sylvester lower bound}
Let $n$ be a positive integer. Then
$$\min\left(\left\{1-\sum_{i=1}^n\frac{1}{m_i}\Biggm| m_i\in\mathbb N^+\right\}\cap (0,1)\right)=\frac{1}{S_{n+1}-1}.$$
\end{thm}

\begin{lem}\label{lem: estimation of epsilon(p,q)}
Let $p,q$ be positive integers such that $p\geq 2$. Then 
$$\epsilon_1(p,q)\geq\frac{1}{S_{(pq+1)p+1}-1}\text{ and }\epsilon_2(p,q)>\frac{1}{S_{(pq+1)p+2}}\geq\frac{1}{S_{4p^2q}}>\frac{1}{2^{2^{4p^2q}}}$$
In particular, $$\alpha(p,2)>\frac{1}{2^{2^{12p^2}}}.$$
\end{lem}
\begin{proof}
Suppose that $$\epsilon_1(p,q)<\frac{1}{S_{(pq+1)p+1}-1}.$$
Then there exist a positive integer $k$ and $\gamma_1,\dots,\gamma_k\in\Phi_p\backslash\{0\}$, such that $$\sum_{i=1}^k\gamma_i\in\left (q,q+\frac{1}{S_{(pq+1)p+1}-1}\right).$$ We have $$\gamma_i=\frac{n_i-1+\frac{m_i}{p}}{n_i}$$ for each $i$, where $n_i\in\mathbb N^+$, $m\in\mathbb N$, and $0\leq m_i\leq p$ for each $i$. Since $\gamma_i>0$ for each $i$ and $p\geq 2$, $\gamma_i\geq\frac{1}{p}$ for each $i$. Thus $q+1\leq k\leq pq+1$. We have
\begin{align*}
    &\sum_{i=1}^k\sum_{j=1}^{p-m_i}\frac{1}{(k-q)n_ip}=\sum_{i=1}^k\frac{p-m_i}{(k-q)n_ip}\\
    =&\frac{1}{k-q}\left(k-\sum_{i=1}^k\gamma_i\right)\in \left(1-\frac{1}{(k-q)(S_{(pq+1)p+1}-1)},1\right)\subset\left(1-\frac{1}{S_{(pq+1)p+1}-1},1\right).\\
\end{align*}
This contradicts Theorem \ref{thm: sylvester lower bound}. Thus $$\epsilon_1(p,q)\geq\frac{1}{S_{(pq+1)p+1}-1}.$$ 
It is clear that $\epsilon_1(p,q)<1$.  We have
$$\epsilon_2(p,q)=\frac{\epsilon_1(p,q)}{q+\epsilon_1(p,q)}>\frac{\epsilon_1(p,q)}{q+1}>\frac{1}{(q+1)(S_{(pq+1)p+1}-1)}>\frac{1}{S_{(pq+1)p+2}}.$$
The last two inequalities are straightforward.
\end{proof}

\begin{lem}\label{lem: estimation of alpha}
Let $p,q$ be two positive integers such that $q\geq 2^{2^{12p^2}}$. Then $$\alpha\left(p,\frac{1}{q}\right)>\frac{1}{(2\uparrow)^4q}.$$
\end{lem}
\begin{proof}
We have $$M_{\frac{1}{q}}=(2q)^{128q^5}(2q+2)<q^{q^6}.$$ So
$$\alpha\left(p,\frac{1}{q}\right)>\frac{1}{2^{2^{4p^2q^{q^6}}}}>\frac{1}{2^{2^{q^{q^7}}}}>\frac{1}{(2\uparrow)^4q}.$$
\end{proof}

\begin{prop}\label{prop: lower bound of beta}
Let $p\geq 2$ be an integer. Then:
\begin{enumerate}
    \item $\beta(p)>\frac{1}{(2\uparrow)^{14}12p^2}>\frac{1}{(2\uparrow)^{17}p}.$
    \item $l(p)<(2\uparrow)^{17}p$.
    \item $\upsilon(p)>\frac{1}{(2\uparrow)^{19}p}$.
\end{enumerate}
\end{prop}
\begin{proof}
The second inequality of (1) is straightforward. This first inequality of (1) follows from Lemmas \ref{lem: estimation of epsilon(p,q)} and \ref{lem: estimation of alpha}. (2)(3) follow from (1).
\end{proof}

\begin{lem}\label{lem: upper bound of beta}
Let $p$ be a positive integer. Then $\beta(p)\leq\min\{\frac{1}{6},\frac{1}{p(p+1)}\}$.
\end{lem}
\begin{proof}
Since $$\frac{1}{p}+\frac{p}{p+1}-1=\frac{1}{p(p+1)}$$ and $$\frac{1}{2}+\frac{1}{3}-1=\frac{1}{6},$$
we have 
$$\epsilon_1(p,1)\leq\min\left\{\frac{1}{6},\frac{1}{p(p+1)}\right\}.$$ Thus $$\beta(p)\leq\min\left\{\frac{1}{6},\frac{1}{p(p+1)}\right\}.$$
\end{proof}

\subsection{Bound of the global lc threshold and exceptional mlds}
\begin{prop}\label{prop: glct 1-gap dim 2}
Let $p$ be a positive integer. Then $$\delta_{\glct}(2,p)\geq\beta(p)>\frac{1}{(2\uparrow)^{17}p}.$$
\end{prop}
\begin{proof}
By Theorem \ref{thm: exceptional 1/42 lc global strong} we may assume that $p\geq 2$. Suppose that $\delta_{\glct}(2,p)<\beta(p)$. Then there exists an an lc pair $(X,B+tS)$ such that $t\in(1-\beta(p),1)$, $\dim X=2,B\in\Phi_p$, $0\not=S\in\mathbb N^+$, and $K_X+B+tS\equiv 0$. We write $B=\sum \gamma_iB_i$ such that $\gamma_i\in\Phi_p$ and $B_i$ are prime divisors. Possibly replacing $(X,B+tS)$ with a dlt model, we may assume that $(X,B+tS)$ is $\Qq$-factorial dlt. We run a $(K_X+B)$-MMP which terminates with a Mori fiber space $X'\rightarrow Z$. Let $B',S'$ be the images of $B,S$ on $X'$ respectively. Then since $X\rightarrow X'$ is $S$-positive, $S'\not=0$. We let $F$ be a general fiber of $X'\rightarrow Z$ and let $B_F:=B'|_F,S_F:=S'|_F$. Then $(F,B_F+tS_F)$ is lc and $K_F+B_F+tS_F\equiv 0$.

If $\dim F=1$, then by Proposition \ref{prop: glct 1-gap dim 1} and Lemma \ref{lem: upper bound of beta}, $$t<1-\min\left\{\frac{1}{6},\frac{1}{p(p+1)}\right\}\leq1-\beta(p),$$
a contradiction. If $\dim F=2$, Then $F=X'$ and $\rho(X')=1$. This contradicts \cite[Theorem 3.2]{AM04}. 
\end{proof}

\begin{prop}
Let $p$ be a positive integer. Then $\delta_{\mld,\ft}(2,p)\geq\beta(p)>\frac{1}{(2\uparrow)^{17}p}$.
\end{prop}
\begin{proof}
By Theorem \ref{thm: exceptional 1/42 lc global strong}, we may assume that $p\geq 2$. The proposition follows from Propositions \ref{prop: compare exc mld and glct} and \ref{prop: glct 1-gap dim 1}.
\end{proof}

\begin{thm}\label{thm: surface volume lower bound}
Let $p$ be a positive integer and $(X,B)$ is a projective lc surface pair, such that $B\in\Phi_p$ and $K_X+B$ is big. Then $\vol(K_X+B)\geq\upsilon(p)>\frac{1}{(2\uparrow)^{19}p}$.
\end{thm}
\begin{proof}
    By Theorem \ref{thm: surface standard volume lower bound}, we may assume that $p\geq 2$. The theorem follows from \cite[Theorem 4.8]{AM04}.
\end{proof}

\end{document}